%% file: Chow.tex
\documentclass{article}

\input{MacrosChow.tex}

\begin{document}

\input IntroduccionChow


\input PreliminaresChow


\input RepresentacionChow


\input ComputacionChow


\input AplicacionesChow


\input ReferenciasChow

\input DireccionesChow

\end{document}

%% file: MacrosChow.tex
\usepackage{amssymb,amsmath,latexsym,epsfig,euscript}

\def\Id {{\rm \mbox{Id}}}
\def\Ch {{{\cal C}{\it h}}}

\def\Conv {{\rm \mbox{Conv}}}
\def\Vol {{\mbox{\rm Vol}}}
\def\vol {{\rm \mbox{vol\,}}}
\def\mod {{\rm \mbox{mod\, }}}

\def\No {{\mbox{\rm N}}}
\def\Norm {\mbox{\rm Norm}}

\def\Supp {{\rm \mbox{Supp\,}}}

\def\ord {\mbox{\rm ord\,}}

\def\Res {\mbox{\rm Res}}
\def\aff {{\rm \mbox{\scriptsize aff}}}

\def\PowerSeries {\mbox{\rm PowerSeries}}

\def\Expand {\mbox{\rm Expand}}
\def\GradedParts {\mbox{\rm GradedParts}}
\def\PolynomialDivision {\mbox{\rm PolynomialDivision}}
\def\Norm {\mbox{\rm Norm}}
\def\Newton {{\cN}}

\def\ChowForm {\mbox{\rm ChowForm}}
\def\Eval {\mbox{\rm Eval}}

\def\Ar { \mbox{\rm ar}}
\def\size { \mbox{\rm size}}
\def\a {{\rm \mbox{\scriptsize a}}}
\def\e {{\rm \mbox{\scriptsize e}}}
\def\m {{\rm \mbox{\scriptsize m}}}
\def\Prob { \mbox{\rm Prob}}
\def\CompanionMatrix { \mbox{\rm CompanionMatrix}}
\def\lin {{\rm \mbox{\scriptsize lin}}}

\def\Homog {\mbox{\rm Homog}}
\def\Adjoint {\mbox{\rm Adjoint}}
\def\JacobianMatrix {\mbox{\rm JacobianMatrix}}

\def\NumDenNewton {\mbox{\rm NumDenNewton}}

\def\discr {\mbox{\rm discr}}

\def\e {{ \mbox{\rm \scriptsize e}}}
\def\h {{ \mbox{\rm \scriptsize h}}}
\def\red {{ \mbox{\rm \scriptsize red}}}
\def\Equidim {\mbox{\rm Equidim}}
\def\Sep {\mbox{\rm Sep}}
\def\GCD {\mbox{\rm GCD}}

\def\Inter {\mbox{\rm Intersection}}
\def\GeomRes {\mbox{\rm GeomRes}}

\def\Random {\mbox{\rm Random}}

\def\A{\mathbb{A}}
\def\C{\mathbb{C}}
\def \N{\mathbb{N}}
\def \P{\mathbb{P}}
\def \Q{\mathbb{Q}}
\def \R{\mathbb{R}}
\def \Z{\mathbb{Z}}

\def\cA {{\cal A}}
\def\cB {{\cal B}}

\def\cF {{\cal F}}
\def\cI {{\cal I}}
\def\cJ {{\cal J}}
\def\cK {{\cal K}}
\def\cL {{\cal L}}
\def\cM {{\cal M}}
\def\cN {{\cal N}}
\def\cO {{\cal O}}
\def\cP {{\cal P}}
\def\cQ {{\cal Q}}

\def\cS {{\cal S}}

\def\cU {{\cal U}}
\def\cV {{\cal V}}

\def\cZ {{\cal Z}}

\newtheorem{defn}{Definition}[section]
\newtheorem{lem}[defn]{Lemma}
\newtheorem{mainlem}[defn]{Main Lemma}
\newtheorem{prop}[defn]{Proposition}

\newtheorem{cor}[defn]{Corollary}

\newtheorem{rem}[defn]{Remark}
\newtheorem{exmpl}[defn]{Example}
\newtheorem{assumption}[defn]{Assumption}

\newenvironment{undef}[1]%
           {\vspace{3.3mm}
           \noindent{\bf #1}\it}%
           {\vspace{3.3mm}}

\newenvironment{proof}[1]{
  \trivlist \item[\hskip \labelsep{\it #1}]}{\hfill\mbox{$\square$}
  \endtrivlist}

\newenvironment{algo}[4]
{\begin{table}[!t]

\caption{#1}
\label{#2}

\hrulefill

\vspace{2mm}

\hspace*{4.1mm}{\bf procedure} \ #3

\vspace{3mm}

#4 }{

\vspace{2mm}

\hspace*{4.1mm}{\bf end.}

\hrulefill

\end{table}
}


%% file: IntroduccionChow.tex
\centerline{\LARGE{\bf The computational complexity of the Chow
form}}

\medskip 

\bigskip 

\centerline{\large Gabriela Jeronimo\footnote{Partially supported
by UBACyT EX-X198 (Argentina).}, Teresa Krick\footnotemark[1],
Juan Sabia\footnotemark[1], and Mart{\'\i}n
Sombra$^{1,}$\footnote{Also partially supported by grant UNLP
X-272 (Argentina),
 and by a
Marie Curie Post-doctoral  
Fellowship of the European Community
program {\em Improving Human
Research Potential and the Socio-economic Knowledge Base}, contract
n\textordmasculine \ HPMFCT-2000-00709.}
}

\bigskip

\bigskip


\noindent {\small {\bf Abstract.} We present a bounded probability
 algorithm for the computation of the Chow forms of
the equidimensional components of an algebraic  variety. Its
complexity is  {\em polynomial }
 in the length and in the
geometric degree of the input  equation system defining the
variety. 
In particular, it provides an alternative algorithm for
the equidimensional decomposition of a variety.

As an  application we obtain an algorithm for the computation of a
subclass of sparse resultants, whose complexity is
{\em polynomial}   in the dimension and the volume of the input  
set of exponents. 
As a further application, we derive an algorithm for the 
computation of the (unique) 
solution of a generic over-determined equation system.
}

\vspace{1mm}


\vspace{1mm}

\noindent {\small {\bf MSC 2000.} {\it Primary:}  14Q15, 
{\it Secondary:} 68W30. }


\typeout{contenido}

\tableofcontents

\vspace{\fill}

\pagebreak


\typeout{Introduccion}

\section*{Introduction}

\addcontentsline{toc}{section}{Introduction}

\vspace{2mm}

The Chow form of an  equidimensional quasi-projective variety is 
one of
the basic objects of algebraic geometry and plays a central role
in  elimination theory, both from the theoretical and practical
points of view.

\smallskip

Let $V \subset \P^n$ be an equidimensional quasi-projective
variety  of dimension $r$ and degree $D$ defined over $\Q$. Its
{\em Chow form } $\cF_{V}$  is a polynomial with rational
coefficients ---unique up to a scalar factor--- which
characterizes the set of over-determined linear systems over the
projective closure $\overline{V}$. 
More precisely, let
$U_0,\dots, U_r$ denote $r+1$ groups of $n+1$ variables each, and 
set
$ L_i:=U_{i  0}\, x_0+\cdots + U_{i
n}\, x_n $  for the linear form associated to the group $U_i$ for 
$0\le i \le r$. 
Then  $\cF_{V}\in
\Q[U_0, \dots, U_r]$ is the unique ---up to a scalar factor---
squarefree polynomial such that 
$$ {\cal F}_{V}(u_0,\dots,u_r)=0
\iff \overline{V} \cap \{ L_0(u_0,x)=0,\dots,L_r(u_r,x)=0\} \ne
\emptyset $$ 
for $u_0, \dots, u_r \in \C^{n+1}$. This is a
  multihomogeneous polynomial of degree $D$ in each group of variables
$U_i$. Thus we can directly read the dimension and the degree of $V $ from $\cF_V$.
In case $V$ is closed in $\P^n$, 
its Chow
form completely characterizes it, and it is possible to
 derive a complete set of  
{equations} for $V$ from $\cF_V$.

\bigskip

The main result of this work
is that  the computation of the Chow forms of
all  
the equidimensional components of a quasi-projective variety 
defined   by means of
a given set  of polynomials,  has a 
   {\em
polynomial} complexity 
in terms of the number of variables, 
the degree, and also   the  {\em length} and the
{\em geometric
degree} of
the input system. 

\smallskip

The complexity of
all known general algorithms in algebraic geometry  is (at least)
{\em exponential} in the worst-case when the considered input parameters are
just 
the number of variables and the degree of the input system, and
there is strong evidence that this exponential behavior is
unavoidable (see \cite{HeMaPaWa98}).
However, it has been observed that there are many particular instances
which are much more tractable than the general case.
This has motivated the introduction of parameters
associated to the
input system that  identify these particular cases, and the design of
algorithms whose performance is correlated to these parameters.

In this spirit,  the  {\em straight-line program}
representation  of polynomials (slp for short) was introduced in
the polynomial
 equation solving frame  as an alternative data structure
 (see e.g. \cite{GiHe93},
\cite{GiHeSa93}, \cite{GiHeMoPa95}) and its {\em length} is 
now considered
as a meaningful  parameter measuring the input (see
Section \ref{Data and algorithms structure} below for the
definition of these notions and \cite{Strassen72}, \cite{Heintz89},
\cite{BuClSh97} for a broader background).

Soon afterwards, the notion of 
{\em geometric degree} of the input polynomial system  
appeared naturally as
another
 parameter classifying tractable problems. 
Roughly
speaking, this is a parameter which measures 
the degree of the varieties successively cut out by
the  input polynomials
(we refer to Section \ref{Equations in general position} below for 
its  precise definition).
For $n$ polynomials of degree $d$, the geometric degree is 
always bounded by the B{\'e}zout
number $d^{n}$. 
However, there are many situations in which the geometric degree is
much
smaller than this B{\'e}zout 
number (see  \cite[Section 4]{KrSaSo97} for a particular
example or
\cite[Prop. 2.12]{KrPaSo99} for an analysis of the sparse case).

In   \cite{GiHeMoMoPa98} and  \cite{GiHaHeMoPaMo97}, 
J. Heintz,  M. Giusti, L.M. Pardo and their collaborators succeeded in
classifying   the tractability   of
 polynomial equation solving  {\em in the  zero-dimensional case} in
 terms of the length and the geometric degree   of
  the input system. They presented an algorithm whose complexity 
is {\em polynomial } in  the number of
variables,  the degree, the length and the geometric degree   of
  the input system. 
Our main theorem can be seen as an extension of   their result  
to arbitrary  dimensions:

\begin{undef}{Theorem 1}
\label{Theorem}
Let
$f_1, \dots, f_s , \, g \, \in {\Q} [x_0, \dots, x_n]$
be homogeneous polynomials  of degree bounded by $d$
encoded by straight-line programs of length bounded by $L$.
Set $V :=  V(f_1, \dots, f_s) \setminus V(g) \subset \P^n $ for
the quasi-projective variety
$\{ f_1 = 0, \dots, f_s = 0, g\ne 0\}$
and let
$ V = V_0 \cup \dots \cup V_{n} $
be  its   minimal equidimensional
decomposition.
Set $\delta:= \delta( f_1,\dots, f_s; \, g)$ for the geometric degree
of the input polynomial system.

Then there is a bounded probability algorithm which computes
(slp's for)
the Chow forms $\cF_{V_0}, \dots, \cF_{V_n}$
 within (expected) complexity
$
s( n\, d\, \delta )^{\cO(1)} L$. Its worst case complexity is
 $s( n\, d^n )^{\cO(1)} L$.
\end{undef}

Let us precise the formal  frame for our computations.
The support of our algorithms
is the
model of  bounded error probability {\em  Blum-Shub-Smale machine}  (BSS machine
for short) over $\Q$: 
our algorithms are probabilistic
BSS machines that
 manipulate slp's.
A probabilistic
BSS machine is the algebraic 
analogue of the classical notion of probabilistic Turing
machine, in which  the bit operations are replaced by the arithmetic
operations $\{ + , - , \cdot, /\}$  of $\Q$.
It allows  to implement {\em uniform} procedures while
``programming'' using the basic operations of $\Q$.
This model is well suited to control the algebraic complexity ---that
 is the number of arithmetic operations--- performed by the algorithm.


By  {\em bounded error probability } we mean that the error
probability of the machine is {\em uniformly} bounded by $1/4$.
For us, the (expected) complexity is then  the expectation of the
complexity seen as a random variable, and {\em not}  its 
worst-case. 
The choice of the constant $1/4$ as error probability
is not restrictive: 
for any given $N \in \N$
we can easily modify our machine (by running it $\cO (\log N)$
times) so that the final error probability is bounded by \ $1/N$
(see Proposition \ref{Probabilidad de Ms} and Corollary
\ref{remarkrepeticion} below).

We refer to Section \ref{Data and algorithms structure} for a
detailed description and discussion of the data structure and
computational model.

\bigskip

We note that in our situation, it is 
unavoidable to 
 consider a data structure different 
from the dense representation 
 if we look  for
a polynomial time algorithm. 
If
  $V$ is an equidimensional  variety of dimension $r$
defined as the zero set of a family of polynomials of
degrees bounded by $d$, its Chow form $\cF_{V} $ is a polynomial
of degree   $(r+1) \, \deg V$   in   $(r+1) \, (n+1)$
variables. So the dense representation of $\cF_{V} $ (i.e. the
vector of its coefficients)  has
$$
{ (r+1)\, (n+1) +  (r+1) \, \deg V \choose (r+1) \, (n+1)  }\ \ge
\
 \frac{(\deg V)^{(r+1)(n+1)}}{((r+1) (n+1))!}
$$
entries, and hence 
it is not polynomial in $\deg V$  (which in the worst case  equals $
d^{n-r}$).
In fact, Corollary \ref{L(Ch)} below   shows that in the above situation
 the slp representation of  $\cF_{V} $ has length $L(\cF_V)= (n\,d
 \deg V)^{\cO(1)}L$.

\bigskip

The previous algorithms for the computation of the Chow forms of 
the equidimensional components of 
a
positive-dimensional variety
(\cite{Krick90}, \cite{Caniglia90}, \cite{GiHe91}, \cite{PuSa98})
have an essentially  worse complexity than ours, with the
exception of the
one due to
G. Jeronimo, S. Puddu and J. Sabia (\cite{JePuSa01}),
which computes the Chow form of the component of maximal dimension
of an algebraic variety within  complexity $(s\, d^n)^{\cO(1)}$.
Here, not only we  compute the Chow form of  {\em all} of the
equidimensional components but we also replace
the B{\'e}zout number $d^n$ by $d\, \delta$, 
where $\delta$ denotes  the geometric degree.

Note that our algorithm also provides an effective geometric
{\em equidimensional
decomposition}, as each equidimensional component is characterized 
by
its Chow form.
 
We can easily derive equations or a geometric
resolution of a variety from its Chow form (see 
\ref{GeomResFromChow} below). 
Hence our result contains and improves --- maybe up to a polynomial 
---  
all the previous  symbolic general  algorithms for
this task \cite{ChGr83}, \cite{GiHe91}, 
\cite{ElMo99b}, \cite{JeSa02}, \cite{Lecerf00}, \cite{Lecerf01}. 

Of these, the only ones whose  complexity is comparable to ours  are the
algorithms of Jeronimo-Sabia \cite{JeSa02} --- which computes
equations for each equidimensional component --- and G. Lecerf
\cite{Lecerf00}, \cite{Lecerf01} ---which computes  geometric 
resolutions. 
With respect to the first one, we improve the complexity by replacing
the quantity  $d^n$ by the geometric degree $\delta$. 
With respect to the second one, we 
substantially improve both the
probability and the uniformity aspects. 
Thus our work might also be seen as a 
unification and clarification  of 
these equidimensional decomposition algorithms.

\smallskip

On the contrary, it is by no means evident 
how to compute the Chow form 
from  a geometric resolution within an admissible complexity.
The difficulty lies in that the involved morphism is not finite but
just {\em  dominant} (see Remark \ref{dominante}). 
We exhibit a deterministic algorithm 
which performs this task within complexity $(s\, n\,  \deg V)^{\cO(1)}
L$ (Main Lemma \ref{fibra}). 
In fact this is the key in our algorithm, and probably the main
technical contribution of the present work.

\bigskip

As a first application of our algorithm,
 we   compute a subclass of sparse resultants.

 \smallskip

The classical resultant $\Res_{n,d}$ of a system of $n+1$ generic
homogeneous polynomials in $n+1$ variables is a polynomial in the
indeterminate coefficients of the polynomials, that characterizes
for which coefficients the system has a non-trivial solution. 
We show that 
$\Res_{n,d}$ can be deterministically computed within 
complexity 
$ (n d^n )^{\cO(1)}$ (Corollary \ref{resdensa})

\smallskip

The sparse resultant  $\Res_\cA$ 
---{\em the} basic object of sparse  elimination
theory---  
has extensively been used 
as a tool for the resolution of polynomial equation
systems
(see for instance
\cite{Sturmfels93}, \cite{Rojas97}, \cite{ElMo99}). 
Several
effective procedures were 
proposed 
to compute  it 
(see e.g. \cite{Sturmfels93}, \cite{CaEm93},
\cite{CaEm00}). 
Recently, C. D'Andrea 
has obtained an explicit determinantal formula 
which extends 
Macaulay's formula to the sparse case (\cite{DAndrea01}).

{}From the algorithmic point of view, the basic point of sparse
elimination theory is that computations
should be substantially faster when  the input polynomials are
sparse (in the sense that their Newton polytopes are restricted). 
Basically,
the parameters which control the sparsity  are the number of
variables $n$ and the normalized volume $\Vol(\cA)$ of the convex
hull of the set $\cA$ of exponents. None of the previous
algorithms computing sparse resultants is completely
satisfactory, as their predicted complexity is exponential in all
or some of these parameters (see \cite[Cor. 12.8]{CaEm00}).

\smallskip

The precise definition of $\Res_\cA$ is as follows:

Let $\cA = \{\alpha_0,
\dots, \alpha_N\} \subset \Z^n$ be a finite set of integer
vectors. We assume here that $\Z^n$ is generated by the differences of
elements in $\cA$.
For $i=0, \dots, n$
let  $U_i$ be a group of variables indexed by the elements of
$\cA$,  and set
$$
F_i := \sum_{\alpha \in \cA} U_{i \alpha }  \, x^{\alpha} \ \in
\Q[U_i][x_1^{\pm 1}, \dots, x_n^{\pm 1}]
$$
for the generic Laurent polynomial with support equal to 
$\cA$. 
Let $
W_\cA \subset (\P^{N})^{n+1} \times (\C^*)^n  $ be the incidence
variety of $ F_0, \dots, F_n$ in $ (\C^*)^n $, that is
$$
W_\cA= \{ (\nu_0, \dots, \nu_n; \, \xi) ; \ \  F_i (\nu_i, \xi)=0
\ \ \forall\, 0\le i\le n\},
$$
and let $ \ \pi : (\P^{N})^{n+1} \times (\C^*)^n \to
(\P^{N})^{n+1} \ $ be the canonical projection.

The sparse {\it $\cA$-resultant} \ $\Res_\cA$ is defined as the unique
---up to a sign--- irreducible polynomial in
 $ \Z[U_0, \dots, U_n]$
which defines
$ \overline{\pi (W_\cA )} $, which
is an irreducible variety of codimension 1 (see \cite[Chapter 8,
Prop.-Defn. 1.1]{GeKaZe94}).
It coincides with  the Chow form of the toric variety associated to
the input set $\cA$ (see Section
\ref{The Chow form of a quasi-projective variety} below).

We show that
the computation of $\Res_\cA$ in  case $\cA \subset (\N_0)^n$
and $\cA$ contains 
$ 0, e_1, \dots, e_n $  
---the vertices of the
standard simplex of $\R^n$--- 
is an instance of our main 
algorithm (see Subsection \ref{Sparse resultants}).
We thus obtain:

\begin{undef}{Corollary 2}
Let $\cA \subset (\N_0)^n $ be a finite set which contains $\{0 ,
e_1, \dots, e_n \}$. Then there is a bounded probability
algorithm which computes (a slp  for) the   $\cA$-resultant
$\Res_\cA$
 within (expected) complexity
$ (n \, \Vol(\cA))^{\cO(1)}$. Its worst case complexity 
 $(n d^{n})^{\cO(1)}$, where $d:= \max \{ |\alpha | \, ; \,
\alpha \in \cA\}$.
\end{undef}

Hence our result represents a significant 
improvement in the theoretical
complexity of computing the $\cA$-resultant as we show it is {\em
polynomial} in $n$ and $\Vol(\cA)$.
We remark that to achieve this result, we had to abandon all 
matrix
formulations. 
In fact, this polynomial behavior of the complexity is out of reach 
of the known matrix formulations, as in all of  them the involved
matrices have an exponential size.

It would be desirable to extend this algorithm in order to compute 
general {\em mixed} resultant, a
point will be a subject of our future research.

\bigskip

As a further application,  we compute the
unique solution of a generic over-determined system over 
an equidimensional variety
$V\subset \P^n $:

\smallskip

Set $r:= \dim V$.
Let  $u= (u_0, \dots, u_r ) \in (\C^{n+1})^{r+1} $ and  set
 $
\ell_i := L_i(u_i, x) = u_{i   0} \, x_0 + \cdots + u_{i  n} \,
x_n$ for $0\le
i\le r$.

There is a
non-empty open subset of the set of coefficients $u$ satisfying
the condition that the linear forms  $\ell_0,
\dots, \ell_r $ have at least a  common root in $V$, for which the
common root is unique. 
It turns out that this unique solution 
$\xi(u)$ is   a rational function
of the coefficient $u$, and this function can be easily
computed using the Chow form $\cF_V$. 
We can successfully apply our algorithm
to that situation
 (see Section \ref{Generic over-determined
systems} for the details).

\bigskip

The results presented
above suggest the following
new research line: 

It is usual to  associate to  a family of $s$ polynomials of
degrees $d_1,\dots,d_s$ in $n$ variables such that $d_1\ge \dots
\ge d_s$, the {\em B{\'e}zout number} $D:=d_1 \cdots
d_{\min\{n,s\}}$.
 The main property of this B{\'e}zout number is that it
  bounds the geometric degree of the variety
defined by the input  polynomials. 

We claim that  the
precise definition of such a B{\'e}zout number $D$ should depend
intimately  on  the {\em encoding}  
of the input polynomials: for
polynomials of degree $d$ encoded in  {\em dense}
representation, $D:=d^n$ is a good choice, while for {\em sparse}
polynomials with support in $\cA$, $D:=\Vol(\cA)$ seems to be the
right notion of B{\'e}zout number, as this quantity also controls the
degree of the variety. 

This digression is motivated by the following  crucial
observation: in our computation of the resultant, the
length of the input $L$ together with the associated B{\'e}zout
number $D$  and the number of variables $n$ controls the
complexity. 
In case of dense
 representation of the input and   $d\ge 2$, 
we have $L=\cO(n d^n)$ and  $D=d^n$
 while for the  sparse representation, we have 
$L\ge 1$ and  $D=\Vol(\cA)$. 
In both cases,   
the complexity of computing the resultant is   
$(n D)^{\cO(1)}L$.
We think that    
the 
optimal
complexity estimate should in fact be 
{\em linear} 
 in $  D$ as well. 
On the other hand, 
it is not clear what should be the exact  
dependence on 
$n$. 
In the linear case, that is for 
  $n+1$ dense linear
  forms, we have $L=\cO(n^2)$ and  $D=1$, 
and the resultant equals  the determinant, which is
  conjectured ---but still
  not proved--- 
to be computable in $\cO(n^2)$.

For the more general 
problem of computing the Chow form of polynomials
encoded by slp's, the notion of
geometric degree $\delta$ plays the
role of the B{\'e}zout number $D$, and our main result
 states that the complexity of computing the Chow form is 
 $s\, (ndD)^{\cO(1)}L$. 
Again, we think that the optimal estimate should be linear in 
$D$, and that the number of polynomials $s$ 
and the  maximum degree $d$ should not occur.

In an even more general framework, we conjecture that 
the computation of  (a slp representation of) any geometric 
object associated to a family of
 polynomials in $n$ variables represented in a given encoding, 
with associated
 B{\'e}zout number $D$ and associated length of the input $L$,
 should be linear in both  $D$ and $L$, and (possibly) quadratic in 
 $n$.

\bigskip

Now we briefly sketch our main algorithm. 

\smallskip

In a first step, we prepare the input data:  
We  take
$n+1$ random 
linear combinations of the input polynomials so that 
---with high probability---
these new polynomials define the same
variety $V$ and   behave properly with respect to the dimensions  
and 
radicality of certain ideals they successively define.
We also take a 
random change of variables to 
ensure
good conditions for the considered projections.

After this preparatory step, we compute  recursively the Chow
forms of the components of a non-minimal equidimensional decomposition of $V$. For
$0\le r\le n-1$,  the algorithm deals with an equidimensional
subvariety $W'_r$  of the variety defined by  the
first $n-r$ polynomials.

The recursive step is as follows:

{}From a {geometric resolution}  of a
zero-dimensional fiber of $W'_{r+1}$, we compute the Chow form of
the variety obtained by intersecting  $W'_{r+1}$ with the set of zeroes
of the next polynomial.  From this Chow form, we obtain the
Chow form of an equidimensional variety of dimension $r$ which is
a subset of $V$ and contains the equidimensional component of
dimension $r$ of $V$ together with a
geometric resolution of the zero-dimensional fiber of $W'_r$ that
is needed for the next recursive step.

The recursion yields the Chow forms of the components of a non-minimal equidimensional decomposition of $V$.
The required {\it minimality} of
the equidimensional decomposition imposes a third step 
in which we remove the 
spurious components.

\bigskip

Finally, a word with respect to   practical implementations: 
there
is a Magma package called Kronecker written by G. Lecerf
(see \cite{Lecerf00b}) which implements
---with
remarkable success---
the algorithm for polynomial equation solving in
\cite{GiHeMoMoPa98}, \cite{GiHaHeMoPaMo97} (in the version of \cite{GiLeSa01}).
Our algorithm is closely related to this one 
and so it seems possible to implement it using this package as support.
However, the lack of a suitable tool for {\em manipulating slp's}  is
still a considerable obstacle for a
successful implementation which would reflect the good theoretical
behavior of our algorithm.

\bigskip

The outline of the paper is the following:

\smallskip

In Section 1 we recall the definition and basic properties 
of the  Chow form, and
we precise    
the data structure and the computational model.
We also describe 
some basic subroutines that we need in the sequel, 
and we estimate 
their complexities.

In Section 2 we present 
a deterministic algorithm for the computation of the Chow form
of an equidimensional variety from a particular zero-dimensional
fiber,
provided some
genericity conditions are
fulfilled.

In Section 3 we describe the algorithm underlying 
Theorem 1, and we estimate its complexity.  
First we establish  
the relationship between geometric resolutions and
Chow forms, and then we present subroutines 
for computing Chow forms of
intersections and of components outside hypersurfaces. 
Combined with the algorithm in Section 2, this 
yields the desired algorithm.

In Section 4, 
we apply the main algorithm to the computation of a 
sparse resultants, 
and to the resolution of generic over-determined equation systems.

\bigskip

\bigskip

\noindent {\bf Acknowledgments.} The authors wish to thank Agnes
Sz{\'a}nt{\'o}  and Gregoire Lecerf for several 
helpful discussions, and to Joos
Heintz for sharing with us his ideas concerning  B{\'e}zout numbers
and complexity.
We also thank Mike Shub for suggesting us the problem of solving
generic over-determined systems.


%% file: PreliminaresChow.tex

\typeout{Preliminaries}

\section{Preliminaries} \label{Preliminaries}


\vspace{2mm}

Throughout this paper
$\Q$ denotes the field of rational numbers,
$\Z$ the ring of rational integers, $\R$ the field of real numbers,
and $\C$ the field of
complex numbers.
We denote by $\N$ the set of positive rational integers, and we also denote
by $\N_0$ the set of non-negative integers.

\smallskip

We denote by  $\A^n$ and $\P^n$ the $n$-dimensional affine space and  projective
space  over $\C$, respectively, equipped with the Zarisky topology definable over $\C$.
A quasi-projective variety $V$ is an open set of a closed
projective (non necessarily irreducible) variety.
We denote by $\overline{V} \subset \P^n $ the projective closure of $V$, that
is the minimal closed projective variety which contains it.

\smallskip

If $f_1, \dots, f_s$ are polynomials in $\Q[x_0, \dots, x_n]$, $V(f_1, \dots, f_s)$ will denote the set of common zeros of $f_1, \dots, f_s$ in $\P^n$. This notation will also be used in the affine case.

\smallskip

Let $V$ be a quasi-projective variety and let $g \in \Q[x_0, \dots, x_n]$ be a
homogeneous
polynomial.
Then we denote by $V_g$ the basic open set  $V \setminus V(g)$ of $V$.

\smallskip

We adopt the usual notion of  degree of  an irreducible projective variety.
  The degree of  an arbitrary
projective variety is here  the sum of the degrees of its irreducible components.
If the variety is quasi-projective,   its degree is defined as the degree of
its projective closure.

\smallskip

We only  consider polynomials and  rational functions  with coefficients in $\Q$
and varieties defined by polynomials with coefficients in
$\Q$ unless otherwise explicitly stated.

\smallskip
The determinant of a square matrix $M$ is denoted by $| M |$.

%


\typeout{The Chow form of a quasi-projective variety}

\subsection{The Chow form of a quasi-projective variety}
\label{The Chow form of a quasi-projective variety}

\vspace{2mm}

We gather in this subsection some definitions and basic facts about Chow
forms.
For a more detailed account we refer to
\cite[Section I.6.5]{Shafarevich72}, \cite[Chapter 3]{GeKaZe94} and
\cite{DaSt95}.

\medskip

First we define the notion of Chow form of an equidimensional
 quasi-projective variety:

\smallskip Let $V \subset \P^n$ be an equidimensional
quasi-projective variety of dimension $r$.


For $i=0, \dots, r$  let
$
{U_i}= ( U_{i  0}, U_{i  1}, \dots, U_{i  n}) $ be a group of
$n+1$ variables and   set $ {U}:= (U_0 , \dots , U_n)$. Then set
$$ L_i:= U_{i  0} \,x_0 +  \cdots + U_{i  n}\, x_n \  \in
\Q[{U}][{x}]$$ for the associated generic linear form, where $x$
denotes the group of variables $(x_0,\dots,x_n)$. Let $$ \Phi_V =
\{ ({u_0}, \dots, {u_r}; \,  { \xi}) \in (\P^n)^{r+1} \times \P^n
\, ; \, \xi \in V , \ L_0 ({u_0},  \xi) = 0, \dots, L_r ({u_r},
 \xi) = 0 \}
\ \subset (\P^n)^{r+1} \times\P^n
$$
be the incidence variety of these linear forms in $V$, and let
$\pi: (\P^n)^{r+1} \times \P^n
\to (\P^n)^{r+1}$ be the projection $({u}, {\xi})
\mapsto {u}$.



\begin{lem}
\label{defchow}
Under the previous assumptions and notations, $\overline{\pi(\Phi_V)}
= \pi(\Phi_{\overline{V}})$.
\end{lem}

\begin{proof}{Proof.--}

Let $V = \cup_C \, C$ be the irreducible decomposition of $V$.
{}From the definition above we deduce that $\Phi_V = \cup_C \,
\Phi_{C}$ and so $\pi(\Phi_V) = \cup_C \, \pi( \Phi_{C})$.

We also have that $\overline{V} = \cup_C \, \overline{C}$ is the
irreducible decomposition of $\overline{V}$. Then
$\pi(\Phi_{\overline{V}}) = \cup_C \, \pi( \Phi_{\overline{C}})$
and so,  without loss of generality, we can assume that $V$ is
irreducible.

The map $\Phi_V \to V$ defined by $({u} , {\xi}) \mapsto {\xi}$
makes $ \Phi_V$ a fiber bundle over $V$ with fiber
$(\P^{n-1})^{r+1}$. Then $\Phi_V$ is an irreducible variety of
codimension  $n+1$, and the same is true for
$\Phi_{\overline{V}}$.

\smallskip

As $ \Phi_{\overline{V}}$ is a closed set,  $\overline{\Phi_V}
 \subset \Phi_{\overline{V}}$. These are irreducible
projective
varieties of the same dimension and, therefore,  they are equal.
 The fact that $\pi$ is a closed map
 implies that
 $\overline{\pi(\Phi_V)} = \pi(\Phi_{\overline{V}})$.
\end{proof}

Then
$\overline{\pi(\Phi_V )} \subset (\P^n)^{r+1} $ is a closed
hypersurface \cite[p.66]{Shafarevich72}.  We define a {\em Chow form of } $V$
  as any  squarefree defining equation $\cF_{V}\in \Q[ {U_0}, \dots, {U_r}]$
of
 the Zariski closure
$\overline{\pi(\Phi_V )} \subset (\P^n)^{r+1} $.  Note that the
Chow form of an equidimensional variety is uniquely determined up
to a scalar factor. We extend this
 to  dimension $-1$ defining a Chow form of the empty variety as any non-zero
 constant in $\Q$.

This definition extends the usual notion of Chow form of closed projective
equidimensional varieties. In fact,
Lemma \ref{defchow} states that a Chow form of an equidimensional
 quasi-projective variety is a Chow form of its
projective closure.

\smallskip

{}From this definition, we see that any Chow form of $V$
characterizes the sets of over-determined linear systems over the
variety $\overline{V}$ which intersect it: for $u_0, \dots, u_r
\in \C^{n+1}$  we have $$ \cF_{V}({u_0}, \dots, {u_r}) = 0 \
\Leftrightarrow
\
\overline{V} \cap \{
L_0({u_0},  {x}) = 0\} \cap \cdots \cap \{
L_r({u_r},  {x}) = 0 \} \ne \emptyset.
$$

\smallskip

A Chow form $\cF_V$ is a multihomogeneous polynomial of degree
$\deg V$ in each group of variables ${U_i} \ (0 \le i \le r)$.
The variety $\overline{V}$ is uniquely determined by a Chow form
of $V$ (\cite[p. 66]{Shafarevich72}). Moreover, it is possible to
derive equations for the variety $\overline{V}$ from a Chow form
of $V$ (\cite[Chapter 3, Cor. 2.6]{GeKaZe94}).

In case
$V$ is irreducible, $\cF_{V}$ is an
irreducible polynomial and, in the general case, a Chow form of $V$ is
the product of Chow forms of its
irreducible components.

\medskip

Following \cite{KrPaSo99} we avoid the indeterminacy of $\cF_V$ by
fixing one of its coefficients under the following assumption on
the equidimensional quasi-projective variety $V$:

\begin{assumption} \label{assumption}
If $\dim V = 0$, we assume $V \subset \{x_0 \ne 0\}$. If  $\dim V
 = r > 0$, we assume that the projection $\pi_V: V  \dashrightarrow  \P^r$ defined by
 ${x} \mapsto (x_0: \dots :x_r)$ verifies $\# \pi_V^{-1}((1:0: \cdots :
 0))
= \deg V $.
\end{assumption}

This assumption implies that
$ \overline{V} \cap
\{x_1=0\} \cap \dots \cap \{x_r=0\} $ is a 0-dimensional variety
lying in the affine chart $\{ x_0 \ne 0\}$.
In particular $V$
has no components contained in the hyperplane
$\{x_0 =0\}$.

We also note  that, in case $V$ is a closed affine variety, the
hypothesis $ \#  (V \cap \{x_1=0\} \cap \dots \cap \{x_r=0\}) =
\deg V$ implies that the map $\pi_V: V \to \A^r$  defined by $x
\mapsto (x_1, \dots, x_r) $ is finite;  that is, the variables
$x_1,\dots, x_r$ are in Noether normal position
 with respect to $V$ (\cite[Lemma 2.14]{KrPaSo99}).

\smallskip

Set ${e_i}$ for the $(i+1)$-vector of the canonical basis of
$\Q^{n+1}$ and $D:=\deg V$. Then, under Assumption
\ref{assumption}, $\cF_V({e_0},\dots,{e_r})$ ---that is, the
coefficient of the monomial $U_{0\,0}^D \cdots U_{r\,r}^D  $---
is non-zero. Then, we  define {\em the  (normalized) Chow form }
$\Ch_V$  of $V$ by fixing the choice  of $\cF_V$ through the
condition $$ \Ch_V ({e_0},\dots,{e_r}) =1 . $$ Note that if $V$
satisfies Assumption \ref{assumption} then each of its irreducible
components also does. Therefore,  the normalized Chow form of $V$
equals the product of the normalized Chow forms of its irreducible
components. The normalized Chow form of the empty variety equals
the polynomial $1$.

\smallskip

Here there are some examples of Chow forms:

\begin{itemize}

\item In
case $\dim  V=0$ we have $$ \cF_V({U_0})= \prod_{\xi \in V}
L_0(U_0,\xi ) \  \in \Q[{U_0}]. $$ Furthermore, if $V$ satisfies
Assumption \ref{assumption},
 $\Ch_V $ is equal to the above expression provided  we choose homogeneous
  coordinates of the type
$\xi:= (1: \xi' )  \in \P^n$ for each point in $V$.

\item In case $V$ is a hypersurface of degree $d$,
then $V= V(F) \subset \P^n$ where $F \in \Q[x_0, \dots, x_n] $ is
a squarefree homogeneous polynomial of degree $d$. We consider
the $n \times (n+1)$-matrix $M:= (U_{i\, j})_{\substack {1 \le i
\le n \\ 0 \le j \le n}}$, and, for $0\le j=0\le  n$, we set $M_j
$ for the maximal minor obtained by deleting its $(j+1)$ column.
Then
$$
\cF_V  = F(M_0, -M_1,  \dots, (-1)^n \, M_n) \ \in \Q[U_0, \dots, U_n].
$$
In this case, Assumption \ref{assumption} is equivalent to the fact that
$f:= F(1, 0, \dots, 0, t)$ is a squarefree polynomial of degree $d$ in $t$.
Therefore,
$\Ch_V $ is equal to the above expression if  we choose
 $F$ such that the coefficient  of the monomial
$x_n^d$ is 1.

\item The  sparse resultant provides an important family of
examples:
let $\cA = \{\alpha_0, \dots, \alpha_N\} \subset \Z^n$ be a finite set
of integer vectors, such that  the
differences of elements in $\cA$
generate  $\Z^n$.
Consider the map
$$
\varphi_\cA: (\C^*)^n \to \P^{N} \quad \quad , \quad \quad \xi
\mapsto (\xi^{\alpha_0}: \dots : \xi^{\alpha_N}) .
$$
This is always well-defined as $\xi_i \ne 0$ $(1 \le i \le n)$ for all $\xi \in
(\C^*)^n$.
The Zariski closure of the image of this map
$ X_\cA:= \overline{\varphi_\cA ((\C^*)^n)} \subset \P^{N} $
 is the {\em toric
  variety} associated to the set $\cA$.
This is an irreducible variety of dimension $n$ and degree
$\Vol(\cA)$ (the normalized volume of the convex hull of $\cA$).

The $\cA$-resultant equals the Chow form of this variety
\cite[Chapter 8, Prop. 2.1]{GeKaZe94}, that is: $$ \cF_{X_\cA} =
\Res_\cA. $$ We refer to
 \cite{GeKaZe94} and to
\cite[Chapter 7]{CoLiOs98} for a broader background on sparse
resultants and toric varieties.

\end{itemize}


\typeout{Data structure and algorithmic model }

\subsection{Data and algorithms structure}

\label{Data and algorithms structure}

\vspace{2mm}

First we specify our data structure:

  The objects we deal with are
polynomials with rational coefficients. The data structure we
adopt to represent them concretely is the {\em straight-line
program} encoding (slp's for short). The input, output and
intermediate objects computed by our algorithms are polynomials
codified through slp's. We insist on the fact that in the present
work  the crucial feature of slp's is their role as data
structures, rather than the more traditional functional role as
programs without branchings for the evaluation of polynomials at
given points.

\smallskip

For the standard terminology of slp's, see  \cite[Defn.
4.2]{BuClSh97}. In this paper all slp's are defined over $\Q$,
without divisions and expecting the variables $x_1, \dots, x_n$ as
input.

\smallskip

For purpose of completeness we restate the definition in our
particular case:

Let $n \in \N$. We denote by $\{ +, -, \cdot\}$ the addition,
substraction and multiplication in the $\Q$-algebra $\Q[x_1,
\dots, x_n]$.   We consider apart the addition and multiplication
by scalars, that is for $\lambda \in \Q$ and $f \in \Q[x_1, \dots,
x_n]$ we set ${\lambda^\a} (f):= f+\lambda $ and ${\lambda^\m} (f)
:= \lambda \cdot f$.   We denote by ${\Q^\a}$ and ${\Q^\m}$ the
set of all scalar additions and multiplications for $\lambda \in
\Q$, respectively.

We set $\Omega_n:= \Q^\a \cup \Q^\m \cup \{ + , - , \cdot \}$ and
denote by  $\Ar(\omega)$ the arity of an operation $\omega \in
\Omega_n$: that is $1$ if it is a scalar operation and $2$ if it
is a binary one.

A {\em straight-line program } $\gamma$ (over $\Q$ and  expecting
$x_1, \dots, x_n$ as input) is a sequence $\gamma:= (\gamma_1,
\dots, \gamma_L)$ of {\em instructions} $$ \gamma_i = \left\{
\begin{array}{lcl}(\omega_i; \, k_{i1}) & \mbox{if}& \Ar(\omega_i)=1\\
(\omega_i; \, k_{i1},k_{i2}) & \mbox{if}& \Ar(\omega_i)=2
\end{array}\right.  $$ where each $\omega_i \in \Omega_n  $ is
 an operation and   for every
$j$, $ k_{i j} \in \Z $ satisfies $-n+1 \le k_{i j} \le i-1$ and
represents a  choice of a previous index. The number of
instructions $L$ in $\gamma$ is called the {\em length} of
$\gamma$ and is denoted by $L(\gamma)$. This is,  in the standard
terminology, the complexity induced by the cost function which
charges 1 to each operation in $\Omega_n$ (see \cite[Defn.
4.7]{BuClSh97}).

\smallskip

Given a slp $\gamma=(\gamma_1,\dots,\gamma_L)$, its {\em result
sequence } $(f_{-n+1} , \dots, f_L) $ is classically defined as
$$f_{-n+1 } := x_1,\dots, f_{0}:=x_n \quad \mbox{and for } 1\le i
\le L, \quad f_i:= \left\{\begin{array}{lcl} \omega_i(f_{k_{i1 }})
& \mbox{if}& \Ar(\omega_i)=1\\   \omega_i(f_{k_{i1 }},f_{k_{i 2}})
& \mbox{if}& \Ar(\omega_i)=2\end{array}\right. .$$

Here we make a slight modification with respect to the standard
terminology. According to the data structure role played by slp's
we  consider only the final result of the slp $\gamma$, that is
the final polynomial $f_L \in \Q[x_1, \dots, x_n]$. We call it the
{\em result} of $\gamma$. Here is the precise definition:
 Let $\Gamma_\Q[x_1,\dots,x_n]$ denote the set of slp's over
$\Q$ expecting $x_1, \dots, x_n$  as input. Then there is a
well-defined surjective function $$ \Eval:
\Gamma_\Q[x_1,\dots,x_n] \to \Q[x_1, \dots, x_n] \quad \quad ,
\quad \quad \gamma \mapsto f_L \qquad \mbox{where} \
L:=L(\gamma).$$ In this way each slp defines precisely one
polynomial (and not a finite set). We say that $\gamma\in
\Gamma_\Q[x_1,\dots,x_n]$ {\em encodes} $f\in \Q[x_1,\dots,x_n]$
if $f$ is the result of $\gamma$. Given a polynomial  $f \in
\Q[x_1, \dots, x_n]$ we define its {\em length } $L(f) $ as the
minimal length of a slp which encodes $f$. (We always have $\deg f
\le 2^{L({f})}$.) For a finite set $\cP \subset \Q[x_1, \dots,
x_n]$ we define naively its {\em length } as $L(\cP):=\sum_{f\in
\cP} L(f)$.

\smallskip

{}From the dense representation $\sum_\alpha a_\alpha \, x^\alpha$
of a polynomial $f\in \Q[x_1,\dots,x_n]$ we   obtain
straight-forward  a slp for $f$.  First, it is easy to show
inductively that for any $r\in \N$, there is a slp of length
bounded by $n+r\choose r$ whose result sequence is the  set of all
monomials $x^\alpha$ of degree $|\alpha|\le r$. This is   due to
the fact that once one has a list of all monomials of degree
bounded by $r-1$, each one of the $n+r-1\choose r$ homogeneous
monomials of degree $r$ is simply obtained from one of the list
multiplying by a single variable. Now set $d:=\deg f$. We multiply
all monomials of degree bounded by $d$  by their coefficients and
 add them up, that is we add $2\,{n+d \choose d}$ instructions
to the slp, in order  to obtain a slp which encodes $f$. Hence $$
L(f) \le 3\,{n+d\choose d} \le 3\,(d+1)^n.$$

We call this particular slp  the {\em dense slp} of $f$. The
previous computation shows that in the worst case, the length
$L(f)$ of  a polynomial $f$ of degree $d$ is essentially its
number of monomials.

\medskip
We can operate with the data structure slp, extending directly the
operations in $\Omega_n$: for instance for $*\in \{+,-,\cdot\}$,
given two slp's $\gamma, \delta \in \Gamma_\Q[x_1,\dots,x_n]$ we
obtain the  new slp $$\gamma * \delta := *(\gamma,\delta)  :=
(\gamma_1, \dots, \gamma_{L(\gamma)},{\delta}_1 , \dots,
{\delta}_{L(\delta)}, ( * ; L(\gamma), L(\gamma)+L(\delta))), $$
where the choice of previous indexes for $\delta$ are suitable
modified. This slp obviously encodes the $*$ of the two
polynomials encoded by $\gamma$ and $\delta$, and its length is
$L(\gamma) + L(\delta) + 1$.

\smallskip

 More
generally, for $\gamma \in \Gamma_\Q [y_1, \dots, y_m] $ and
$\delta_1, \dots , \delta_m \in \Gamma_\Q [x_1, \dots, x_n] $, we
can define the {\em composition slp }
$
\gamma \circ \delta:=\gamma \circ (\delta_1,\dots,\delta_m) \in
\Gamma_\Q [x_1, \dots, x_n] $. We have $$ L(\gamma \circ \delta )
= L(\delta_1)+\cdots + L(\delta_m)  + L(\gamma). $$

This operation is  compatible with the map $\Eval$, that is
$\Eval(\gamma \circ \delta) = \Eval(\gamma) \circ \Eval (\delta)$.
Hence for $f \in  \Q [y_1, \dots, y_m] $ and $g_1, \dots, g_n \in
\Q[x_1, \dots, x_n]$ we have that $L(f(g_1, \dots, g_n)  ) \le
L(g_1)+\cdots + L( g_n)+L(f)$.

\bigskip
Now we specify the computational model that produces and
manipulates our data structure: it is  the {\em Blum-Shub-Smale
machine}
 over $\Q$
(BSS machine for short), which captures the informal notion of
uniform algorithm over $\R$. We refer to \cite[Chapters 3 and
4]{BlCuShSm98} for the definition, terminology  and basic
properties. However, there are again some slight
 modifications in our definition (restrictions on
 the operations --only over rational numbers-- and the branchs
 --only equality of numbers to zero), and we restate it for
  purpose of
 completeness:

\smallskip

We recall that a BSS machine $\cM$ over $\Q$ has five types of
nodes:
 input,
computation, branch, shift, and output.

Set $$ \Q^\infty := \bigsqcup_{n \ge 0} \Q^n $$ for the disjoint
union of the $n$-dimensional spaces $\Q^n$, i.e. the natural space
to represent problem instances  of arbitrarily large dimension.
For $a\in \Q^\ell \subset  \Q^\infty$ we call $\ell$ the {\em
size} of $a$, and we denote it  by $\size(a) $.

On the other hand, let $$ \Q_\infty:= \bigoplus_{m\in \Z} \Q $$ be
the bi-infinite direct sum space over $\Q$. The elements  $b\in
\Q_\infty$ are of the form $$ b:= (\dots, b_{-2}, b_{-1}, b_0 \,
.\,  b_1, b_2, \dots), $$ where $b_i =0 $ for $|i| \gg 0$. The dot
between $b_0$ and $b_1$ is a distinguished marker which allows to
visualize the position of the coordinates of $b$.

\smallskip

Now we define the computation maps: For each $\omega \in \{+, -,
\cdot, / \}$ and $i, j, k \in \N$ there is a map $$ \Q_\infty \to
\Q_\infty \quad \quad , \quad \quad b \mapsto (\dots, b_{k-2},
b_{k-1}, \omega(b_i, b_j) , b_{k+1}, \dots)$$ (observe that unlike
in the case of our data structure, here we also allow division).
Also for each $\lambda^a \in \Q^\a$ or $\lambda^m \in \Q^\m$ and
$i , k \in \N$ there is in an analogous way a map $$\Q_\infty \to
\Q_\infty \quad \quad , \quad \quad b \mapsto (\dots, b_{k-2},
b_{k-1}, \lambda^{a,m} (b_i) , b_{k+1}, \dots). $$  These are our
computation nodes.

The only  branch node we allow is the one  which checks the
equality $b_1=0$. In other words its associated map is $\Q_\infty
\to \Q \quad, \quad b \mapsto b_1$.

The shift nodes are  of two types: shifting to the left or to the
right, associated with the two  maps $\Q_\infty \to \Q_\infty
\quad, \quad b \mapsto \sigma_l (b)_i=b_{i+1} $ or $\sigma_r
(b)_i=b_{i-1}$.

\smallskip
The machine $\cM$ over $\Q$ is a finite connected directed graph
containing these five types of nodes (input, computation, branch,
shift and output). The space $\Q^\infty$ is both the  input space
$\cI_\cM$ and the output space $\cO_\cM$, and $\Q_\infty$ is the
state space  $\cS_M$, that is the ``working'' space of the
machine. The dimension  $K_\cM$ of the machine $\cM$ is the
maximum dimension of the computation maps, which, under our
assumptions, coincides with the maximum of the natural numbers $i,
j$ and $k$ involved in the computations.

\medskip

We are interested in the {\em algebraic } complexity of these
machines. We assume that  the cost  of each computation and branch
node is 1, and that of the shift nodes is 0. In other words, we
assume that $\cM$  can do arithmetic operations and equality
questions in $\Q$ at unit cost, and shifts  at no cost. Hence the
{\em complexity } $C_\cM(a)$ of the machine $\cM$
 on an input $a$ is just the
number of computation and branch nodes of the graph,  from input
to output. Note that each computation node performs exactly one
operation, and hence $C_\cM(a)$ is just the number of basic
arithmetic operations and equality questions performed on the
input $a$.

\smallskip
Observe that any slp $\gamma \in \Gamma_\Q[x_1,\dots,x_n]$ is an
example of a (finite-dimensional) BSS machine $M_\gamma$ without
branchs or shift nodes for computing $f:=\Eval(\gamma) \in
\Q[x_1,\dots,x_n]$ at any input point $a\in \Q^n$. The dimension
of this machine is $n+L(\gamma)$ and its complexity is
$L(\gamma)$.

\smallskip

Given $\ell \in \N$ we consider the {\em complexity }
 $C_\cM(\ell) $ of the machine  on inputs of size bounded by
$\ell$, that is $$ C_\cM(\ell ):= \sup \, \{ C_\cM(a)\, ; \,
\size(a) \le \ell \}. $$

\smallskip

This computational model is a natural algebraic analogue of the
notion of Turing machine. In particular it provides us a support
for the implementation of {\em uniform} procedures.

The difference with the Turing model is that we replace
 bit operations by arithmetic ones, and that
we do not count the cost of the shifts operators. The
reason not to count the  shifts  is not fundamental,
and is just that the available literature on the algebraic
complexity of basic polynomial manipulation only takes into
account the number of arithmetic operations and comparisons (see Subsection
\ref{Complexity of basic computations}, see also \cite{BuClSh97}).

Since all  the involved computations are done over the rational
field, the machine $\cM$ can be effectively transformed into a
classical Turing machine. However our complexity counting does not
provide any reasonable control on the complexity of the resulting
Turing machine.

\bigskip

As we have already anticipated, our algorithms are BSS machines
that manipulate slp's.  A machine $\cM$ receives as input a finite
family of slp's  $\gamma \in \Gamma_\Q[x_1,\dots,x_n]^M$ and gives
back a finite family of slp's  $\cM(\gamma) \in
\Gamma_\Q[y_1,\dots,y_m]^{M'}$.

 A finite family of slp's $\gamma\in
\Gamma_\Q[x_1,\dots,x_n]^M$ can be easily codified as an input
element in $\cI_\cM=\Q^\infty$, in fact $\gamma$ can be identified
with  a vector in  $ \Q^{M+3\, L(\gamma)}$ in the following way:

 The first coordinate is for the dimension
$n$, that is the number of variables. Then each instruction of the
first slp $\gamma_1$ is codified as a triple:  we enumerate the
operations in $\Omega_n $ with numbers from 2 to 6, 2 and 3
corresponding to the operations in $\Q^a$ and $\Q^m$, and 4 to 6
to $+, -$ and $\cdot$. For operations in $\Q^\a \cup \Q^\m$ we
write the operation number in the first coordinate, the
corresponding coefficient in the second one, and the position to
which it applies in the third one. The binary operations are
codified in a similar way, by writing first the operation number,
and then the position of the two elements to which it applies. The
positions are numerated from $1-n$ to $L(\gamma)$ according to the
definition of result sequence. For instance the vector
$(2,(3,5,-1),(4,0,1),(6,2,2))$ codifies the slp $x_1,x_2, 5x_1,
x_2+5x_1,(x_2+5x_1)^2$. The instruction to separate two
consecutive slp's is an empty cell, that is a $0$. The second slp
$\gamma_2$ is now codified exactly as the first one. Therefore,
$\gamma:=(\gamma_1,\dots,\gamma_M) \in \Gamma_\Q[x_1,\dots,x_n]^M$
is codified as a vector in $\Q^\infty$, in fact in
$\Q^{M+3\,L(\gamma)}$ since we need to add $M-1$ ``0" to separate
two consecutive slp's.

 \smallskip
The machine $\cM$ manipulates this input, the finite family of
slp's $\gamma \in \Gamma_\Q[x_1,\dots,x_n]^M$: it  operates these
slp's and gives as the output an element of $\cO_\cM$
corresponding to a finite family of slp's in
$\Gamma_\Q[y_1,\dots,y_m]^{M'}$. As we've just seen, the input and
output size is (essentially) the length of each of these families.
Thus, we   speak here of a finite family of slp's $\gamma$ as the
input of $\cM$ and we simply  denote by $\cM(\gamma)$ its output
in $\Gamma_\Q[y_1,\dots,y_m]^{M'}$.

\smallskip

\begin{rem} \label{L(gamma)}
Let $\gamma \in \Gamma_\Q[x_1, \dots, x_n]^M$ be the  input slp
family of a BSS machine $\cM$ and let $\cM(\gamma) \in
\Gamma_\Q[y_1, \dots, y_m]^{M'}$ be its output. Then $$
 L (\cM(\gamma)) \le 3\, L(\gamma) + C_\cM(\gamma).
$$
\end{rem}

\begin{proof}{Proof.--} As we do not know how the machine $\cM$ operates on
$\gamma$, the only bound for  $L(\cM(\gamma))$ is the number of
operations labelled from $2$ to $6$ of the representation of
$\cM(\gamma)$ in $\cM$, which is bounded by the number of non-zero
cells of this representation minus $1$ (since the first cell of
the output corresponds to the number of variables $m$ of the
output). This is bounded by $1 + 3\,L(\gamma) + C_\cM(\gamma) -
1$, that is, the size of $\gamma$ as an input of  $\cM$ (excepting
the $M-1$ zero-cells separating different input slp's) plus the
number of computation nodes $C_\cM(\gamma)$ minus  $1$.
\end{proof}

\bigskip

Our main algorithms are in fact probabilistic. For this reason we
implement them in the model of {\em probabilistic BSS machine}
over $\Q$ \cite[Section 17.1]{BlCuShSm98}. This is a BSS machine
$\cM$  with an additional kind of  node, called probabilistic.
These are nodes that  have two next nodes and no associated map
and that ``flip coins'', that is when a computation reaches a
probabilistic node, it randomly chooses the next node between the
two possible ones with probability 1/2 for each of them.

\smallskip

 In this probabilistic setting, each run ---on the same
given input $\gamma$---
 of the machine $\cM$ may lead to a
different path computation. In our case, for any given input, the
number of probabilistic nodes traversed is finite, that is, the
number of possible paths is finite. We treat the probalistic nodes
as branches and charge cost 1 for each of them.

As every path $\cP$ of $\cM$ corresponds to a BSS machine of
complexity $C_\cP(\gamma)$, the algebraic complexity
$C_\cM(\gamma)$ of the machine $\cM$ on the input $\gamma$ turns
out to be a random variable, with finite sample set. Moreover,
again in our context, every path is finite: it may happen that a
branch ends in an error message but in any case the complexity of
any path is bounded. Thus the random variable $C_\cM(\gamma)$
satisfies

 $$\mbox{Prob}\,(C_\cM(\gamma)=C):= \sum
 \mbox{Prob}(\cP; \ \cP \ \mbox{path   such that}\  {C_\cP(\gamma)=  C} ).$$

 We are interested in the {\em worst case complexity}  $C^{\max}_\cM(\gamma)$, the maximum
complexity of
 the  paths of $\cM$ on $\gamma$, and
the  {\em expected complexity} $E_\cM(\gamma)$, defined as
 the  (finite) expectation of this random variable, that is
 $$E_\cM(\gamma):= E(C_\cM(\gamma))= \sum_{C\in \N} C\cdot\mbox{Prob}(C_\cM(\gamma)=C).$$

Observe that $C^{\max}_\cM(\gamma)\ge E_\cM(\gamma) $ always holds.

\smallskip

As before, we also consider the function $E_\cM: \N \to \N$: $$
E_\cM(\ell) := \sup \, \{ E_\cM(\gamma) \, ; n , M \in \N, \,
\gamma \in \Gamma_\Q[x_1,\dots,x_n]^M \  \mbox{and}   \
M+3\,L(\gamma) \le \ell \}. $$

\medskip
 We define now the error probability of the machine on a
given input. Again, there is here a modification with respect to
traditional probabilistic BSS machines. Keeping in mind that for
any run  of the probabilistic machine $\cM$ on the input $\gamma
\in \Gamma_\Q[x_1,\dots,x_n]^M$, its output (independently from
the path randomly taken)  encodes a finite family of polynomials
$f:=(f_1,\dots,f_M) \in \Q[y_1,\dots,y_m]^{M'}$ we define:

\begin{defn} (Bounded Probability Algorithm) \label{error}
 Given $\gamma  \in \Gamma_\Q[x_1,\dots,x_n]^M$ and given
 a set of polynomials $f:=(f_1,\dots,f_{M'})\in \Q[y_1,\dots,y_m]^{M'}$,   the
{\em error probability} \ $e_\cM(\gamma, f)$ \ that  $\cM$
computes $f$ on the given input $\gamma$ is the probability that
the output of $\cM$ on $\gamma$ does  {\em not} encode $f$; that
is the probability that the computation finishes with an error
message, or that it outputs $\delta \in
\Gamma_\Q[y_1,\dots,y_m]^{M'}$ which does not encode $f$.

We say that
 $\cM$ {\em computes } $f$ if
\ $e_\cM (\gamma, f) \le 1/4$. As this happens at most for one
$f$, when it happens, we set $e_\cM (\gamma):=e_\cM (\gamma, f)$.

When the latter happens for {\em every} input, we say that $\cM$
is a {\em bounded probability  machine} for polynomial slp's, and
we speak of {\em bounded probability algorithm}.
\end{defn}

Observe that our probabilistic machine  is a little unusual since
in fact, as different slp's may encode the same polynomial,  the
polynomial $f$ computed by the probabilistic machine $\cM$
corresponds to   an equivalence class of outputs  rather than a
single one. In this paper, all machines are bounded probability
machines for polynomial slp's in the sense of this definition.

\bigskip

In our setting, probability is introduced by choosing a random
element with equidistributed probability in a set
$[0,\ell)^n:=\{0,\dots,\ell-1\}^n$ for given natural numbers $\ell
$ and $n$. Since probabilistic machines flip coins to decide
binary digits, each of these random choices can be easily
simulated with a machine with complexity $\cO( n \, \lceil\log
\ell\rceil)$, where here and in the sequel, $\log $ denotes
logarithm in base $2$. This machine is denoted  by ${\rm
Random}(n,\ell)$.

In this work, in each case, there is a non-zero  polynomial $F\in
\Q[x_1, \dots, x_n] \setminus \{0\}$ such that a random choice $a$
is good
---that is, leads to the computation of the desired output--- if  $F(a) \ne 0$. The
error probability of this random choice is then estimated by the
Zippel-Schwartz's  zero test (\cite{Zippel79}, \cite[Lemma
1]{Schwartz80}): $$
 \Prob(F(a) = 0) \le \frac{\deg F  }{\ell}.
$$

\smallskip
The choice of $1/4$ as a bound for the error probability is not
restrictive and we can easily modify it in order to reduce the
error probability as much as desired. The usual procedure is to
run the machine $\cM$ many times  and to declare that the
polynomial family $f$ is computed by $\cM$ if it is the output of
more than half the times. There is a slight difficulty here,
appearing from the fact that our machine computes slp's instead of
polynomials, and two different runs may lead to different
encodings of the same polynomials. That is why we need here to be
more careful in our definition. We define it  in the following
way:

Given the bounded probability  machine $\cM$ which on input
$\gamma \in \Gamma_\Q[x_1,\dots,x_n]^M$ computes $f\in
\Q[y_1,\dots,y_m]^{M'}$, and given $s\in \N$, the machine $\cM_s$
is the machine which performs the following:

\begin{enumerate}
\item
$\cM_s$ runs  $s$ times the machine $\cM$ on the  given input
$\gamma$: for $1\le i\le s$ it   obtains the   output slp family
$\delta_i$ together with the complexity $C_i$ of the path followed
to compute $\delta_i$.
\item Then $\cM_s$ chooses randomly $a\in [0, M'\, 2^{s+ 3\,L(\gamma)
+ C_1+\dots + C_s})^m$
and computes $\delta_i(a),\ 1\le i \le s$.
\item For $j=1$ to $\lceil s/2 \rceil$
\begin{itemize}
\item
  it computes $\delta_j(a)-\delta_k(a), j<k\le s$, and compares the results to $0$.
  \item if   $0$ is obtained for strictly more than $s/2$ values of $k$,
  it returns the polynomial family encoded by $\delta_j$ as the
  output and ends.
  \item if not, it goes to $j+1$
  \end{itemize}
 If for no $j\le \lceil s/2 \rceil$ the algorithm obtains 0 enough
 times,
 it outputs error and ends.
 \end{enumerate}

\begin{prop} \label{Probabilidad de Ms} Given  a bounded probability machine $\cM$
which on $\gamma \in \Gamma_\Q[x_1,\dots,x_n]^M$ computes $f \in
\Q[y_1,\dots,y_m]^{M'}$  and given $s\in \N$,
 the worst case complexity, the expected
complexity and the error probability of the machine $\cM_s$ on
$\gamma$ verify the following bounds:
\begin{eqnarray*}
C^{\max}_{\cM_s}(\gamma) & = & \cO\Big((m+1)\,s\,( L(\gamma)  +
C^{\max}_{\cM}(\gamma)) + m\,\log M' \Big) + M' \,{s\choose 2}
\Big) ,\\ E_{\cM_s}(\gamma) & = & \cO\Big((m+1)\,s\,( L(\gamma) +
E_\cM(\gamma)) + m\,\log M'  + M' \,{s\choose 2} \Big) ,\\
e_{\cM_s}(\gamma) &\le & 2\,(3/4)^{s/2}.
\end{eqnarray*}
\end{prop}

\begin{proof}{Proof.--}
Let us first describe the algebraic complexity $C$ of a given run
of the machine   $C_{\cM_s}$ on $\gamma$  in terms of the
complexities $C_i$ of the paths followed by the machine $\cM$ on
$\gamma$  on the $i$-th run.

1. has complexity $C_1 + \cdots + C_s$.

2. Producing the random choice $a$ costs $\cO(m\,(\log M' +
s+L(\gamma) + C_1 + \cdots + C_s))$,  and, from Remark
\ref{L(gamma)}, computing $\delta_1(a),\dots ,\delta_s(a)$ costs
$3 s\, L(\gamma) + C_1 + \cdots + C_s$.

3. As $\delta_j(a)\in \Q^{M'}$, to compute all
$\delta_j(a)-\delta_k(a)$ and compare them to $0$ costs $2\, M'
\,{s\choose 2}$.

Hence, the worst-case complexity of the machine ${\cM_s}$ on $
\gamma$ is
$$C^{\max}_{\cM_s}(\gamma) = \cO\Big((m+1)\, s\,(L(\gamma) +
C^{\max}_{\cM}(\gamma)) + m\,\log M'   + M' \,{s\choose 2}
\Big),$$
while, as the complexity is an affine combination of the $s$ independent
random variables $C_1,\dots,C_s$, its expectation verifies
$$E_{\cM_s}(\gamma)= \cO\Big((m+1)\,s\,( L(\gamma) +
E_\cM(\gamma)) + m\,\log M'  + M' \,{s\choose 2} \Big).$$

\smallskip

The error is bounded by the probability that there is no group of
more than $ s/2 $ vectors which coincide plus the probability that
$\delta_j(a)=\delta_k(a)$ but the two polynomial families encoded
by $\delta_j$ and $\delta_k$ do not coincide.

The first error is bounded by $(3/4)^{s/2}$ as in \cite[Sect.17.2,
Lemma 1]{BlCuShSm98}. To estimate the second error we apply
Schwartz' lemma: for $1\le i\le s$ the output $\delta_i$ codifies
$f\in \Q[y_1,\dots,y_m]^{M'}$ where the degree of each component
is bounded by $2^{3\,L(\gamma) + C_i}$. Thus the error of one test
is bounded by $(M'\, 2^{3\,L(\gamma) + C_i})/(M'\,2^{s+
3\,L(\gamma) + C_1 + \cdots + C_s})\le (1/2)^s$. As there are at
most ${s\choose 2}$   such independent tests, the total error
verifies $$e_{\cM_s}(\gamma) \le  (3/4)^{s/2} +{s\choose
2}\,(1/2)^{s}\le 2\,(3/4)^{s/2}$$ for $s\ge 2$.
\end{proof}

 \begin{cor} \label{remarkrepeticion} Given  a bounded probability machine $\cM$
which on $\gamma \in \Gamma_\Q[x_1,\dots,x_n]^M$ computes $f \in
\Q[y_1,\dots,y_m]^{M'}$  and given
  $N\in \N$, $N\ge 4$,
 the error probability of the machine $\cM_s$ on $\gamma$
 for $s:= \lceil {6\, (\log N  + 1)} \rceil$ is
 bounded by $1/N$ while  its  worst-case  complexity  is of order
 $$ \cO\Big((m+1)\,\log N\,( L(\gamma) +
C^{\max}_{\cM}(\gamma))   + m\,\log M'    + M' \,\log^2N\Big) .$$
 \end{cor}

 \begin{proof}{Proof.--} As  $ (3/4)^3<(1/2)$,
$$e_{\cM_s}(\gamma) \le  2\,(3/4)^{3\,( \log N+1)} \le (3/4)^{3\,
\log N } \le  (1/2)^{\log N} = 1/N.$$
 \end{proof}

  Proposition \ref{remarkrepeticion} will be
 used to decrease the error probability of intermediate
 subroutines of our main algorithm  and keep  control of the
 complexity  in order that the error probability of the latter is
 bounded by $1/4$. Observe that the length of the output slp
 is of the same order that the length of the slp obtained when
 running any of the repetitions of the algorithm.

\smallskip
Given a bounded probability machine $\cM$, any time we want to
obtain the output of $\cM$ for a slp input family $\gamma$ with
error probability bounded by $1/N$, we run Subroutine
\ref{svecesM} which gives a new probability machine $\cM(\gamma;
N)$ doing so. Any time we run $\cM$  for the input family
$\gamma$, we will denote by ${\rm Complexity}(\cM(\gamma))$ the
complexity of doing it this time.

\newpage

\begin{algo}{Decreasing error probability of $\cM$}{svecesM}
{$\cM (\gamma ; N)$}

\hspace*{4.1mm}\# $\gamma$ is a slp input family for $\cM$.

\hspace*{4.1mm}{\# $N\in \N, N \ge 4$.

\smallskip

\hspace*{4.1mm}\# The procedure returns the output of $\cM$ with
error probability bounded by $1/N$.

\begin{enumerate}

\item $s:= \left\lceil 6\,(\log N +1)\right\rceil$;\vspace{-1mm}

\item {\bf for } $i$ { \bf from } $1$ { \bf to } $s$  { \bf do}\vspace{-2mm}

\item \hspace*{4.1mm} $(\delta_i , C_i):= (\cM(\gamma), {\rm Complexity}(\cM(\gamma)))$;
\vspace{-2mm}

\item {\bf od};\vspace{-2mm}

\item $a:= \Random(m, M' \, 2^{s + 3 L(\gamma) + C_1 + \cdots +C_s})$;\vspace{-2mm}

\item $(\Delta_1,\dots, \Delta_s):= (\delta_1(a), \dots, \delta_s(a))$;\vspace{-2mm}

\item $j: = 1$;\vspace{-2mm}

\item {\bf while } $j\le \lceil s/2\rceil$ { \bf do }\vspace{-2mm}

\item \hspace*{4.1mm} $k:= j+1$;\vspace{-2mm}

\item \hspace*{4.1mm} $t:= 0$;\vspace{-2mm}

\item \hspace*{4.1mm} {\bf while } $k\le s$  {\bf do}\vspace{-2mm}

\item\hspace*{8.2mm} {\bf if } $\Delta_j - \Delta_k = 0$ { \bf then}\vspace{-2mm}

\item\hspace*{12.3mm} $t:= t+1$;\vspace{-2mm}

\item\hspace*{8.2mm} {\bf fi};\vspace{-2mm}

\item\hspace*{8.2mm} $k:= k+1$;\vspace{-2mm}

\item\hspace*{4.1mm} {\bf od};\vspace{-2mm}

\item\hspace*{4.1mm} {\bf if } $t\ge s/2$ { \bf then}\vspace{-2mm}

\item\hspace*{8.2mm} {\bf return}$(\delta_j)$;\vspace{-2mm}

\item\hspace*{4.1mm} {\bf else}\vspace{-2mm}

\item \hspace*{8.2mm} $j:= j+1$;\vspace{-2mm}

\item\hspace*{4.1mm} {\bf fi};\vspace{-2mm}

\item {\bf od};\vspace{-2mm}

\item {\bf return}(``error''); \vspace{-2mm}

\end{enumerate}}
\end{algo}



\typeout{Complexity of basic computations. }

\subsection{Complexity of basic computations}

\label{Complexity of basic computations}

\vspace{2mm}

We summarize the complexity of the basic operations on polynomials
and matrices  our algorithms rely on. We refer to  \cite{BuClSh97}
for a rather complete account of the subject or to \cite{GiLeSa01}
for a brief survey of the existing literature.

\smallskip

Let $R$ denote a commutative $\Q$-algebra and let $d \in \N$.

The multiplication of $d\times d$-matrices with coefficients in
$R$ can be done with $\cO(d^{2.39})$  operations of $R$ and no
branches (\cite[Cor. 15.45]{BuClSh97}).  The computation of the
coefficients of the
 characteristic polynomial of a $d\times d$-matrix
---and in particular the computation of the  adjoint and
the  determinant of this matrix--- can be done with  $\cO(d^4)$
arithmetic
 operations and no branches, the same bounds hold for the inversion
 of an invertible matrix
(\cite{Berkowitz84}, \cite{Abdeljaoued97}).

\smallskip

The quantity $$ M (d):= d \, \log d \, \log \log d $$ controls the
complexity of the basic arithmetic operations (addition,
multiplication, and division with remainder) for univariate
polynomials with coefficients in $R$ of degree bounded by $d$ in
dense representation: Addition of univariate polynomials can be
done in $d+1$ arithmetic operations, while the
Sch{\"o}nhage-Strassen polynomial multiplication algorithm takes
$\cO(M(d))$ arithmetic operations and has no branches (\cite[Thm.
2.13]{BuClSh97}). Division with remainder
---provided the divisor is a monic  polynomial--- has also
complexity $\cO(M(d))$ and no branches (\cite[Cor.
2.26]{BuClSh97}). The greatest common divisor can be computed
through subresultants with $\cO (d)$ branches (computing the
degree of the gcd corresponds to checking the vanishing of the
determinant of submatrices of the Sylvester matrix) and complexity
$\cO(d^5)$ (solving a linear system) (\cite {Collins67}, \cite
{BrTr71}). Alternatively, the Knuth-Sch{\"o}nhage algorithm could
be used to compute the greatest common divisor with complexity
$\cO(M(d) \, \log d )$ and $\cO( d^3)$ branches (\cite[Cor.
3.14]{BuClSh97}).

\medskip
Now we are going to consider some procedures involving polynomials
encoded by slp's.

\smallskip

First, given a   slp $ \gamma $ which computes  $f \in \Q[x_1,
\dots, x_n]$  and given  $a\in \Q^n$, we can compute   $f(a) \in
\Q $ within complexity $L(\gamma)$ and so, we can also check $f(a)
=0$ within the same complexity.
The derivative of the polynomial $f$ with respect to one of its variables will be computed by means of the Baur-Strassen's algorithm (see \cite{BaSt83}) within complexity $O(L(\gamma))$.

\smallskip
For a group of variables $y:=(y_1,\dots,y_m)$ and $a\in \Q^m$, we
will denote by $\Expand(f,y,a,d)$ the subroutine
 which, given a  slp $\gamma$ which encodes a multivariate polynomial $f$,
  computes as intermediate results slp's for  the homogeneous
  components centered at $a$
  and of  degree bounded by $d$ of
  the
 polynomial $f$
with respect to the given group of variables $y$: In
\cite[Lemma
13]{KrPa96}, \cite[Lemma 21.25]{BuClSh97} are given  slp's  of
length $\cO(d^2 \, L(\gamma))$ in which  all the homogeneous
components of $f$ of degree bounded by $d$ appear as intermediate
  computations. These procedures  can be
easily modified within the same complexity to compute the
homogeneous components centered at $a$ and  up to degree $d$ of a
polynomial with respect to the given group of variables. In
particular, if $y$ consists of a single variable and $a=0$, this
procedure computes the coefficients of the given polynomial with
respect to $y$.

\smallskip

Quite frequently we use a mixed representation of $f$:   instead
of encoding it by means of a single slp, we consider $f$ as a
polynomial in a distinguished variable, and if $d$ is a bound for
the degree of $f$ in this variable, we  give a $(d+1)$-uple of
slp's, which  encode the  coefficients $f_0,\dots,f_d$ of $f$
with respect to the
 variable. The length of this  mixed encoding does not essentially
 differ from
 the length of
 $f$; denote by $L'(f)$ the length
  of the mixed encoding and by $L(f)$ the length of $f$,
  we have:
$$L(f)= \cO (d + L'(f)) \quad {\rm and} \quad L'(f) = \cO (d^2\,L(f)).$$

\medskip

Sometimes we  need   to compute the exact degree of a polynomial
with respect to a particular variable. We will call ${\rm Deg}(f,
d)$ the procedure which computes the degree of the univariate
polynomial $f$ given by its dense representation, where $d$ is a
bound for its degree. This computation is done by simply
comparing the coefficients of $f$ with 0. This procedure can be
adapted to obtain a probabilistic algorithm ${\rm
Deg}(f_1,\dots,f_s, x, d; N)$ which computes, with error
probability bounded by $1/N$, the total degrees of the
polynomials $f_1,\dots, f_s$ in the group of variables $x$, from
slp's encoding $f_1,\dots, f_s$ and an upper bound $d$ for their
degrees in the variables $x$. To do so, first we apply subroutine
${\rm Expand}(f_i, x, 0, d)$ for $1\le i \le s$, to obtain the
homogeneous components of $f_i$. Then by choosing a random point
in $[0,1,\dots,  sdN )^n$ we decide probabilistically which is
the component of greatest degree different from zero of each
polynomial $f_1,\dots, f_s$. If the given polynomials are encoded
by slp's of length bounded by $L$, the worst-case complexity of
this procedure is of order $\cO( s d^2 L + n \log(sdN))$.


\typeout{Effective division procedures}

\subsection{Effective division procedures}

\label{Effective division procedures}

\vspace{2mm}

Here, we gather the division subroutines  we will need. Basically,
they compute the division of multivariate polynomials and power
series, and the greatest common divisor of multivariate
polynomials. In all cases, the objects will be {\em multivariate}
polynomials {\em encoded by slp's} and power series, whose known
graded parts will be also {\em encoded by slp's}. The proposed
procedure for multivariate power series division is new and plays
an important role in Subroutine \ref{algo-fibra}, which in turn is
the key step of our main algorithm.

\medskip

The following subroutine  is the well-known  Strassen's Vermeidung
von Divisionen (division avoiding) algorithm (\cite{Strassen73}).
We re-prove it briefly in order to estimate its complexity.

\begin{algo}{Polynomial Division }{PolynomialDivision}
{$\PolynomialDivision (f,g,d,a)$ }

\hspace*{4.1mm}{\# $f,g \in \Q[x_1,\dots, x_n]$ such that $f$
divides $g$,

\hspace*{4.1mm}\# $d \in \N$ an upper bound for the degree of the
quotient $g/f$,

\hspace*{4.1mm}\# $a \in \Q^n$  such that $f(a) \ne 0$.

\smallskip

\hspace*{4.1mm}\# The procedure returns $h := g/f$.

\begin{enumerate}

\item $\alpha:= f(a)$;\vspace{-2mm}

\item $v:=
\frac{1}{\alpha}\sum_{i=0}^d\left(\frac{t}{\alpha}\right)^i$;\vspace{-2mm}

\item $H = g \cdot v(\alpha-f)$;\vspace{-2mm}

\item $(H_0, \dots, H_d):= \Expand(H, x, a , d)$;\vspace{-2mm}

\item $h:= \sum_{m=0}^d H_m $;\vspace{-2mm}

\end{enumerate}

\hspace*{7.7mm} {\bf return}($h$); }

\end{algo}

\begin{lem} \label{vermeidung}
Let $f,g \in \Q[x_1, \dots, x_n]$ be polynomials encoded by slp's
of length bounded by $L$
 such that
$f$ divides $g$. Let  $d \in \N$  be  such that $\deg ( g/f ) \le
d$, and
 $a \in \Q^n $  such that $f(a) \ne 0$.

Then $\PolynomialDivision$ (Subroutine \ref{PolynomialDivision})
computes  $g/f$  within complexity $\cO(d^2\, (d+L))$.
\end{lem}

\begin{proof}{Proof.--}
The quotient polynomial  $h:=g/f \in \Q[x_1, \dots, x_n]$ can also
be seen as a power series in $\Q[[x-a]]$. For $\alpha:=f(a)$, we
have $$ h = \frac{g}{f}= g \,{\alpha}^{-1} \, \left(1-
\frac{\alpha- f}{\alpha} \right)^{-1} = g \, \sum_{i = 0}^\infty
\frac{(\alpha -f)^i}{\alpha^{i+1}} \ \in \Q[[x-a]]. $$ For $$H:= g
\, \sum_{i= 0}^d {(\alpha-f)^i}/ {\alpha^{i+1}} \in \Q[x_1, \dots,
x_n]$$ we have $ h \ \equiv \  H \ \mod{(x-a)^{d+1}}.$ Thus, if $
(H_m)_{m\le d}$ are the homogeneous components of $H$ centered at
$a$ and of degree bounded by $d$ , we conclude $h = \sum_{m=0}^d
H_m $.

\smallskip

The stated complexity is obtained as follows: We compute the
univariate polynomial $v$  with $\cO(d+L)$ operations. Hence we
compute $H$ within complexity $\cO(d + L)$. We compute its
homogeneous components in $x-a$ up to degree $d$   within
complexity \ $\cO(d^2(L(H)))=\cO(d^2 (d+L))$. Finally we obtain
$h $ as  $ \sum_{m=0}^d H_m $ within the same complexity bound.
\end{proof}

\smallskip

 Observe that the same procedure can be used to compute
the graded parts centered at $a$ and of
 a certain  bounded degree
of the rational function $g/f$, even in case $f$ does not divide $
g$. We  denote this subroutine by $\GradedParts(f,g,D,a)$, where
the argument $D$ corresponds to the bound  for the degree of the
graded parts to be computed. Its complexity is  of order $ \cO(D^2
(D +L))$.

\medskip
Subroutine \ref{PolynomialDivision} converts slp's {\em with}
divisions computing polynomials in $\Q[x_1,\dots,x_n]$  into
ordinary slp's: Slp's with divisions are defined
 as ordinary  slp's,
but with the set of basic operations  enlarged to include the
division,  which we denote by the bar $/$.   A further requirement
is that all divisions should be {\em well-defined}, that is, no
intermediate denominator can be zero. In general, the result of a
slp with divisions is a rational function in $\Q(x_1,\dots,x_n)$.

\smallskip

Observe that, given a slp with divisions $\gamma$ which encodes a
rational function $h$, we can easily compute
 {\em separately}
a numerator $g$ and a denominator $f$ by means of two slp's
$\zeta, \eta$ without divisions: for instance, for each {\em
addition }
$
h_i := h_{j} + h_{k}
$
 in the result sequence of $\gamma$, if $h_j := h_{j 1} / h_{j 2}$ and $h_k : =h_{k 1}/
 h_{k 2}$,
we set $ g_k:= h_{j 1}\,h_{k 2} + h_{ j 2}\,h_{k 1} $ and $ h_k:=
h_{j 2}\,h_{k 2}$ for the corresponding result sequence in $\zeta$
and $\eta $  respectively. We proceed analogously for the other
operations in $\Omega_n' := \Omega_n \cup \{ /\}$.

We have $$ h:=\Eval (\gamma) = \frac{\Eval (\zeta)}{\Eval(\eta)}.
$$ Furthermore the slp's $\zeta$ and $\eta$
  can be computed within complexity
 $L( \zeta) \le  3 \, L(\gamma)$ and $ L( \eta ) \le  L(\gamma)$.
In particular, given $a \in \Q^n$, we can check if  $\gamma$ is
well-defined at $a$ and, if that is the  case, if $h(a) =0$ within
complexity $\cO(L(\gamma))$.

In case $h$ is a {\em polynomial} of degree bounded by $d$, the
previous considerations together with Lemma \ref{vermeidung} show
that we can compute a slp {\em without divisions} for $h$ with
complexity $\cO(d^2 ( d  + L(\gamma)))$.

\bigskip

Now follows a bounded probability algorithm (in the sense of
Definition \ref{error}) to compute the greatest common divisor of
two multivariate polynomials encoded by slp's (for another
algorithm solving this task see  \cite{Kaltofen88}). Herein, ${\rm
GCD1}(F, G, d, e)$ is the subresultant algorithm which computes a
greatest common divisor of two univariate polynomials $F$ and $G$
of degrees $d$ and $e$ respectively with coefficients in a ring
$A$. The output $(q, Q)$ of GCD1 is the  multiple $Q$  modulo the
factor $q\in A - \{ 0 \}$ of the monic greatest common divisor of
$F$ and $G$ over the fraction field of $A$.

\begin{algo}{Greatest Common Divisor }
{GreatestCommonDivisor} {$ \GCD ( f,  g, x, d)$ }

\hspace*{4.1mm}{\# $f, g \in \Q[x_1,\dots, x_n]$  of degrees
bounded by $d$;

\hspace*{4.1mm}\# $x:= (x_1,\dots, x_n)$;

\smallskip
\hspace*{4.1mm}\# The procedure returns $ h  := \gcd(f,g)$.

\begin{enumerate}

\item  $a:= {\rm Random}( n, 8\,d(d+1))$;\vspace{-2mm}

\item {\bf if} $f(a) = 0$ {\bf then}\vspace{-2mm}

\item \hspace*{4.1mm} {\bf return}{(``error'')};\vspace{-2mm}

\item {\bf else}\vspace{-2mm}

\item \hspace*{4.1mm} $(f_0,\dots, f_d):= \Expand(f, x, a, d)$;\vspace{-2mm}

\item \hspace*{4.1mm} $(g_0, \dots, g_d):=\Expand(g,x,a,d)$;\vspace{-2mm}

\item \hspace*{4.1mm} $e:= 0$;\vspace{-2mm}

\item \hspace*{4.1mm} {\bf while} $g_e(a) = 0$ and $e \le d$ {\bf do}\vspace{-2mm}

\item \hspace*{8.2mm} $e:= e+1$;\vspace{-2mm}

\item \hspace*{4.1mm} {\bf od};\vspace{-2mm}

\item \hspace*{4.1mm} {\bf if} $e = d+1$ {\bf then}\vspace{-2mm}

\item \hspace*{8.2mm} {\bf return}$(f)$;\vspace{-2mm}

\item \hspace*{4.1mm} {\bf else}\vspace{-2mm}

\item \hspace*{8.2mm} $F:= \sum_{k=0}^d f_{k}\, t^{d-k}$ and
 $G:= \sum_{k=0}^e g_{k}\, t^{e-k}$;\vspace{-2mm}

\item \hspace*{8.2mm} $(q,Q):={\rm GCD1}(F,G,d, e)$;\label{GCD1}\vspace{-2mm}

\item \hspace*{8.2mm} $h:= \PolynomialDivision(q(x),Q(x,1),d,q(a))$;\vspace{-2mm}

\item \hspace*{8.2mm} {\bf return}$(h)$;\vspace{-2mm}

\item \hspace*{4.1mm} {\bf fi;}\vspace{-2mm}

\item {\bf fi};\vspace{-2mm}

\end{enumerate}}

\end{algo}

\begin{lem}  \label{mcd}
Let $f, g \in \Q[x_1, \dots, x_n]$ be   polynomials of degree
bounded by $d$  encoded by slp's of length  bounded by $L$.

Then $\GCD$ (Subroutine \ref{GreatestCommonDivisor}) is a bounded
probability algorithm  which computes (a slp for) the greatest
common divisor between $f$ and $g$. Its (worst-case) complexity is
of order \ $ \cO (n\,\log d + d^4 (d^2 + L))$.
\end{lem}

\begin{proof}{Proof.--}

For $a \in \Q^n$   such that $f(a) \ne 0$ and $t$   an additional
variable, we set $$ F (x,t):=  t^ d \, f( \frac{x-a}{t}+a )\ , \
G(x,t):= t^d \, g( \frac{x-a}{t}+a ) \ \in \Q[x][t] . $$ Since
$f(a) \ne 0$,  $F$ is monic
---up to the  scalar factor $f(a)$--- of degree $d$ in $t$.

Set $H$ for the $\gcd$ of $F$ and $G$ in $\Q(x)[t]$.  Since $F$ is
monic in $t$,   $H$ belongs to $ \Q[x,t]$, and it  is easy to
check that $ \gcd(f,g) = H(x,1)$ (up to a scalar factor).

\smallskip
The procedure runs as follows: First we observe that if
$f=\sum_\alpha f_\alpha (x-a)^\alpha$, then $$F=\sum_{0\le k\le d}
\big( \sum_{|\alpha| = k} f_\alpha (x-a)^\alpha \big) t^{d-k},$$
(and the same holds with $g$ and $G$) so the  homogeneous
components of $f$ and $g$ centered at $a$ turn out to be the
coefficients of the monomial expansion of $F$ and $G$ with respect
to $t$. Then, we apply the subresultant algorithm GCD1 to compute
the multiple $Q \in \Q[x,t]$ and the superflous factor $ q$ in
$\Q[x]$ of their gcd $H$ in $\Q(x) [t]$. Finally, we apply
Subroutine \ref{PolynomialDivision} to eliminate divisions in the
expression $\gcd(f,g):= H(x,1) = Q(x,1)/q(x)$.

\smallskip
Let us decide the size of the sets of points we have to take to
insure that the algorithm has an error probability bounded by
$1/4$:

We are going to choose randomly a point $a \in \Q^n$. This same
point $a$ will be used in each step we need a random point.

The first condition the point $a$ must satisfy so that the
algorithm computes a greatest common divisor of $f$ and $g$ is
that $f(a) \ne 0$.

Then we use the point $a$ to compute the degree of $G$ in $t$.
Finally it is used in the subresultant algorithm to compute the
degree of the gcd (by deciding whether certain determinants are
zero or not). Checking the degree of $G$ involves testing an
$n$-variate polynomial of degree bounded by $d$ (the coefficients
of $G$ as a polynomial in
 $\Q[x][t]$) while
checking the degree of the gcd involves testing  $n$-variate
polynomials of degree bounded by $2d^2$.

Thus, applying  Schwartz bound for the set $[0,\ell)^n$, the
conditional probability $p$ of success verifies $$ p \ge
\big(1-\frac{d}{\ell} \big)\,\big(1-\frac{d}{\ell}
\big)\,\big(1-\frac{2\,d^2}{\ell}\big)\ge
1-\frac{d+d+2\,d^2}{\ell} = 1-\frac{2\,d\,(d+1)}{\ell}.$$
Therefore, taking $\ell:= 8\, d \, (d+1) $  insures that the error
probability is bounded by $1/4$.

 \smallskip
 Now let us  compute the worst-case complexity of the machine:

 The cost of simulating the random choices here is $\cO (n\,\log d)$.
 Computing the homogeneous components
 of $f$ and $g$ centered at $a$ and checking the exact degree of $G$
 (that is finding the first non-zero coefficient of $G$ with
respect to $t$) can be done within complexity $\cO(d^2 (d+L))$. In
Algorithm GCD1, to compute the degree of the gcd involves
computing at most $d+1$ determinants of Sylvester-type matrices of
size at most $2d \times 2d$, that is adds at most $
(d+1)\,\cO(d^4) $ operations. Once we know this degree, computing
the gcd by means of an adjoint adds $ \cO(d^4 )$ steps. That is,
the complexity of computing $Q(x,t)$ (and $q(x)$ which is the
non-vanishing determinant) is of order $\cO(d^2(d^3+L))$ while
$L(Q(x,t), q(x))= \cO(d^2(d^2+L))$ since the computation of the
degree does not intervene in the length.
Applying Subroutine \ref{PolynomialDivision} at $q(a)$ which is
different from $0$ (if not, subroutine GCD1 in line \ref{GCD1}
would have returned error) we obtain a final complexity of order
$\cO (n\,\log d + d^2 (d + L(Q(x,t), q(x))))= \cO (n\,\log d + d^4
(d^2 + L)).$
\end{proof}

\medskip
The following procedure (Subroutine \ref{PowerSeriesDivision})
computes the quotient
---provided  it is a polynomial of bounded degree---
 of two multivariate
power series from their graded components up to a certain bound.

\smallskip

Let
\ $
\varphi =  \sum_\alpha a_\alpha \, x^\alpha \ \in \Q[[x_1, \dots, x_n ]]
$ \
be a power series.
For $i \in \N_0$ we denote by \ $\varphi_i :=
\sum_{|\alpha| = i}  a_\alpha \, x^\alpha \in \Q[x_1, \dots, x_n]$ \ the  $i$-graded component of
$\varphi$.
Also we denote by  $\ord \varphi$ the {\em order}  of $\varphi$,
that is   the least
$i$ such that $\varphi_i \ne 0$.

\smallskip

\begin{algo}{Power Series Division }
{PowerSeriesDivision} {$ \PowerSeries ( n,  m, d ,  \, \varphi_m,
\dots, \varphi_{m+d}, \, \psi_m, \dots,
  \psi_{m+d})$ }

\hspace*{4.1mm}{\# $n \in \N$ is the number of variables,

\hspace*{4.1mm}\# $m \in \N_0$ is the order of the denominator $\varphi\in \Q[[x]]$,

\hspace*{4.1mm}\# $d \in \N$ is the degree of the quotient $h:= \psi / \varphi \in \Q[x]$,

\hspace*{4.1mm}\#  the $\varphi_i$'s and $\psi_i$'s
are the
graded  parts of the
power series $\varphi$ and $\psi$
respectively.

\smallskip

\hspace*{4.1mm}\# The procedure returns $ q  := \varphi_m^{d+1}\, h \in \Q[x]$.

\begin{enumerate}

\item  \label{v} $v:= \sum_{i=0}^d y^{d-i} z^i \in
\Q[y,z]$;\vspace{-2mm}

\item \label{P} $P:=  (\sum_{i=0}^{d} \psi_{m+i} \, t^i) \ v  (\varphi_m,
         -\sum_{j=1}^{d} \varphi_{m+j} \, t^j)$;\vspace{-2mm}

\item \label{P0-Pd} $(P_0, \dots, P_d):= \Expand(P, \, t,0,
d)$;\vspace{-2mm}

\item  \label{q} $q:= \sum_{i=0}^d P_i $;\vspace{-2mm}

\end{enumerate}

\hspace*{7.7mm} {\bf return}($q$); }
\end{algo}

\begin{prop}\label{dos polinomios}
Let   $\varphi, \psi \in \Q[[x_1, \dots, x_n]]$ be power series
such that $h:= \psi / \varphi \in \Q[x_1, \dots, x_n]$. Assume we
are given $m := \ord \varphi$, $d \ge \deg h$, and that the
$i$-graded parts of $\varphi$ and $\psi $ for $i =m, \dots, m+d$
are encoded by slp's of lengths bounded by $L$.

Then $\PowerSeries$ (Subroutine \ref{PowerSeriesDivision})
computes $q: =\varphi_m^{d+1} \, h $  within complexity \ $\cO(d^3
L)$.
\end{prop}

\begin{proof}{Proof.--}

Set
$$
\Phi(x,t):= \varphi( t x ) = \sum_{i=0}^\infty \varphi_i (x) \, t^i
 ,
\ \ \ \Psi(x,t):= \psi(t  x ) = \sum_{i=0}^\infty \psi_i (x) \,
t^i \  \ \in \Q[x][[t]] \hookrightarrow \Q(x)[[t]]. $$  Also set
$H:= h(t  x) \in \Q(x)[t]$.

We first observe that $\ord \Phi = m$, and thus
 $\ord \Psi \ge m$ as
$ \Psi/ \Phi = H \in \Q(x)[t] $ \ is a polynomial.  Hence the
following identity holds in $\Q(x)[[t]]$: $$ H =\frac{\Psi}{\Phi}
= \frac{\Psi}{t^m} \, \frac{1}{\varphi_m} \left( 1 -
\frac{\varphi_m - \Phi/t^m}{\varphi_m} \right)^{-1} =
\frac{\Psi}{t^m} \, \sum_{i=0}^\infty \frac{ (\varphi_m -
\Phi/t^m)^{i}}{\varphi_m^{i+1}} . $$

Thus, for  $$ P:= \left(\sum_{i=0}^{d} \psi_{m+i} \, t^i \right)
\left( \sum_{i=0}^d \varphi_m^{d-i} (-\sum_{j=1}^{d} \varphi_{m+j}
\, t^j)^i \right) \ \in \Q[x][t] $$ we have that \ $
\varphi_m^{d+1} \, H \equiv P \ \  \pmod{t^{d+1}}$. Let $P=
\sum_i P_i \, t^i $ be the monomial expansion of $P$. Then
$\varphi_m^{d+1} \, H =  \sum_{i=0}^d P_i \, t^i $, as the degree
of $H$ with respect to $t$ is bounded by $d$. Hence
$\varphi_m^{d+1} \, h  = \sum_{i=0}^d P_i $.

\smallskip

The stated complexity is obtained as follows:

We compute a slp encoding of $v:= \sum_{i=0}^d y^{d-i} z^i$ within
complexity $\cO(d )$. We compute $P$ as $\sum_{i=0}^{d} \psi_{m+i}
\, t^i $ times $v  (\varphi_m , -\sum_{j=1}^{d} \varphi_{m+j} \,
t^j)$ within complexity $\cO( d \, L)$. We compute the expansion
of $P$ with respect to
  $t$ up to degree $d$  within complexity $\cO(d^3 \, L)$.
Finally we compute $q$ as $\sum_{i=0}^d P_i$. The total complexity
is of order $\cO(d^3 \, L)$.
\end{proof}

\begin{rem}
In case that, in addition, we are given $b \in \Q^n$ such that $
\varphi_m(b) \ne 0$, we can directly apply  Subroutine
\ref{PolynomialDivision} to  compute the quotient polynomial $h$
within total complexity $ \cO(d^5 L ).$
\end{rem}


%% file: RepresentacionChow.tex

\typeout{The representation of the Chow form}

\section{The representation of the Chow form} \label{The representation of the Chow form}


\vspace{2mm}

This section presents an algorithm for the computation of the Chow
form of an equidimensional variety from a 0-dimensional fiber and
a set of local  equations at a neighborhood of  this fiber. This
is the key step in our general algorithm  (see Section \ref{The
computation of the Chow form}), although it   has independent
interest: it shows that
 the Chow form and the geometric resolution are ---up to a polynomial time computation---
equivalent representations of a variety (see Subsection
\ref{Geometric resolutions}). As a further application, we give  a
non-trivial upper bound for the length of a slp representation of
the Chow form (Corollary \ref{L(Ch)}).

\medskip

In order to state the result, we need the following definitions:

\begin{defn}
Let $V \subset \P^n$ be an equidimensional variety of dimension
$r$.

We say that $f_{r+1}, \dots, f_n  \in I(V)$ is a {\em system of
local equations} at  $\xi \in V$ if the polynomials
$f_{r+1},\dots, f_n$ generate $I(V)$ at some  neighborhood of
$\xi$, i.e. $I(V)_\xi =(f_{r+1}, \dots, f_n)_\xi$ (where the
subscript $\xi$ denotes localization at the ideal of the point
$\xi$).

If $Z$ is a subset of $V$, we say that $f_{r+1}, \dots, f_n\in
I(V)$ is a {\em system of local equations (of $V$) at $Z$} if it
is a system of local equations at every $\xi \in Z$.
\end{defn}


The existence of a system of local equations at a point $\xi \in
V$ implies that $(\C[x]/I(V))_\xi$ is Cohen-Macaulay and thus, by
\cite[Thm. 18.15]{Eisenbud95}, for $f_{r+1}, \dots, f_n \in I(V)$
to be a system of local equations at $\xi$ is equivalent to the
fact that the Jacobian matrix of this system has maximal rank
$n-r$ at $\xi$.

\begin{defn} Let $Z \subset \A^n$ be a 0-dimensional variety of
cardinality $D$. A
 {\em geometric resolution } of $Z$
 consists of an affine linear form  $\ell=c_0 + c_1 x_1 + \cdots + c_n x_n \in \Q[x_1,\dots,x_n]$
 and of polynomials
$p  \in \Q[t]$ and $v = (v_1, \dots, v_n)  \in \Q[t]^n$  such
that:

\begin{itemize}

\item  The affine linear form
$\/\ell $ is a {\em primitive element} of $Z$, that is $\ell (\xi)
\not= \ell (\xi')$ for all  $\xi\ne  \xi' $ in $ Z$.

\item  The polynomial $p$ is monic of degree $D$ and $p (\ell (\xi)) =0 $
for all $ \xi \in Z$;
that is, $p$ is the minimal polynomial of $\ell$ over $Z$.

\item  $\deg v_i \le D-1$,  $1\le i \le n$, and
$Z =
  \left\{
 v(\eta) \, ; \ \eta\in \C  \ / \ p(\eta) = 0
\right\}$; that is, $v$ parametrizes  $Z$ by  the zeroes of $p$.
\end{itemize}
\end{defn}

Observe that the minimal polynomial $p$ and the parametrization
$v$ are uniquely determined by the variety $Z$ and the affine
linear form $\ell$. We say that $(p,v)$ is {\em the geometric
resolution of $Z$ associated to $\ell$}.

\smallskip
In case $Z\subset \P^n$ is a zero-dimensional projective variety
which satisfies  that none of its points lie in the hyperplane
$\{x_0=0\}$, $Z $ can be identified to  a 0-dimensional {\em
affine} variety $Z^\aff$, the image of $Z$ under the rational map
$\P^n \dashrightarrow \A^n$ defined by $(x_0: \dots: x_n) \mapsto
(x_1/x_0, \dots, x_n/x_0)$. By a geometric resolution  of $Z$ we
then understand a geometric resolution
---as defined before---
of the affine variety $Z^\aff \subset \A^n$. In homogeneous
coordinates, the definition  of  geometric resolution states that
the homogeneized linear form $\ell^h$ satisfies $\ell^h/x_0(\xi)
\ne \ell^h/x_0 (\xi') $ for all $\xi \ne \xi'$ in $Z$. The
polynomial $p$ is then the minimal monic polynomial of $\ell/x_0$
over $Z^{\aff}$. On the other hand, $v$ defines a  parametrization
$ V(p) \to Z$, \ $ \eta \mapsto (1: v_1(\eta): \cdots :
v_n(\eta))$.

\medskip

Now, we are able to state the lemma:

\begin{mainlem} \label{fibra}
Let $V \subset \P^n$ be an equidimensional variety of dimension
$r$ and degree $D$ which satisfies Assumption \ref{assumption}.
Set     $ Z :=V\cap V(x_1,\dots,x_r)$, and let    $p\in  \Q[t] $
and $v \in \Q[t]^n$  be a given geometric resolution of $Z$
associated to a given affine linear form $\ell  \in
\Q[x_1,\dots,x_n]$. Let $ f_{r+1}, \dots, f_{n}\in I(V)$ be a
system of local equations at $Z$. Assume that $f_{r+1}, \dots,
f_n$ have degrees bounded by $d$ and are encoded by slp's of
length bounded by $L$.

Then there is a deterministic algorithm  (Procedure $\ChowForm$
(Subroutine \ref{algo-fibra}) below)  which computes $ \Ch_V$
within complexity $\cO(r^8 \,
  \log_2 (r   D) \, n^7 \, d^2 \, D^{11} L) $.
\end{mainlem}

In Subsection \ref{The algorithm} we present   the complete proof
of the correctness of the algorithm and its   complexity estimate.
The algorithm  is essentially based on a new Poisson type product
formula for the Chow form (see Proposition \ref{forma de chow}
below), which describes the Chow form as a quotient of products of
norms of certain polynomials. We interpret this expression as a
quotient of two power series, which can be approximated with the
aid of a symbolic version of Newton's algorithm. Finally we apply
Procedure PowerSeries (Subroutine \ref{PowerSeriesDivision} above)
to compute the Chow form from the approximate quotient.


\typeout{ Newton's algorithm }

\subsection{Newton's algorithm } \label{Newton's algorithm}


In this subsection we present a symbolic version of Newton's
algorithm for the approximation of roots of equations.
Newton's algorithm is a widely used tool for polynomial equation
solving. The situation in the present work is not much different
from that in e.g. \cite{HeKrPuSaWa00},  \cite{GiLeSa01}. Hence we
just describe this procedure in order to adapt it to our setting
and to estimate its complexity; its correctness  follows directly
from \cite[Section 2]{HeKrPuSaWa00} and the arguments therein.

\medskip
First, we   state the situation in which Newton's
algorithm is applied.


\smallskip

Let $W \subset \A^r \times \A^n$ be an equidimensional variety of
dimension $r$ such that
the projection map $\pi: W \to \A^r$ is dominant, that is, the image
$\pi(W)$ is a Zariski dense set.

Set $A:= \Q[t_1, \dots, t_r] = \Q[\A^r ]$ and let $K $ be its fraction
field.
Also let $B:= \Q[W]$ and set  $L:= K \otimes_A B$.
Then $L$ is a finite $K$-algebra, and its dimension $D:= [L:K]$  ---that is the degree of $\pi$---
equals
the maximum cardinality of the {\em finite}  fibers of $\pi$
(\cite[Prop. 1]{Heintz83}).

The
 {\em  norm
  }  $\No_\pi (h) \in K $ of a
 polynomial  $h \in A[x_1, \dots, x_n]$
is defined as  the determinant of the $K$-linear map $L \to L$ defined
 by
$f \mapsto h \, f$.

Let $I(W)^\e $ denote the extension of the ideal $I(W)$ to the
polynomial ring $K[x_1, \dots, x_n]$, and
set \ $W^e:= V(I(W)^\e) \subset \A^n(\overline{K})$, which is
a 0-dimensional variety of degree $D$.
Then
$$
\No_\pi(h) = \prod_{\gamma \in W^e} h(\gamma).
$$

We  also denote this norm  by $\No_{W^e} (h) $ when the projection map is
clear from the context.


\smallskip

In different steps, we will be given a polynomial $h$ and an
equidimensional variety $W$  and our aim will be to compute an
approximation of $\No_{W^e} (h)$. The input of the procedure will
be the polynomial $h$, a geometric resolution of a 0-dimensional
fiber of $\pi$ and local equations at this fiber.


\smallskip

Let $ F_{1} , \dots, F_n \in   I(W) \subset  A[x_1, \dots, x_n] $
and set $F:= ( F_1, \dots, F_n )$. Let $$ \cJ_F:=
\left(\frac{\partial F_i}{\partial x_j } \right)_{1 \le i , j \le
n} \ \in A[x_1, \dots, x_n]^{n \times n} $$ be the Jacobian matrix
of $F$  with respect to the dependent variables $x_1, \dots, x_n$,
and let $\Delta_F:= | \cJ_F | \in A[x_1, \dots, x_n]$ be its
Jacobian determinant.

Let $Z \subset \A^n$ be such that
 $\pi^{-1} (0) = \{0\} \times Z$.
 We assume that
 $Z$ is a 0-dimensional
 variety of cardinality $D$
and that $\cJ_F$ is non-singular at $\pi^{-1} (0) $, that is,
$\Delta_F(0,  \xi ) \ne 0$ for all $\xi \in Z$. Observe that this
means that
 $F(0,x)$ is a system of local equations
at $Z$. 


\smallskip

Under our assumptions, the elements in $W^e$ can also be
considered as power series:
For $\xi \in Z$, the fact that $\Delta_F(0, \xi) \ne 0$ implies that  there
exists a unique  \ $\gamma_\xi \in \C[[t_1, \dots, t_r]]^n$ \
such that:
$$
\gamma_\xi(0) = \xi \quad  \quad \hbox{and}  \quad \quad F(t_1,
\dots, t_r, \gamma_\xi) = 0.
$$
(See, for example,  \cite[Ch. 3, Section 4.5, Cor.
2]{Bourbaki61}.) It follows that   $f(t_1,\dots, t_r, \gamma_\xi)
= 0$ for all $f \in I(W)$ as $F$ is a system of local equations at
$\xi$, and so this also holds for all $f \in I(W^e)= I(W)^\e$.
Hence $ \gamma_\xi \in W^e$ and, as $\# Z = \# W^e = D$, we
conclude that the correspondence
$$
Z \to W^e \quad \quad , \quad \quad \xi \mapsto  \gamma_\xi
$$ is one-to-one.
In particular,  since $\No_{W^e}(h)$ is the
determinant of a matrix in ${\Q(t_1, \dots, t_r)}$,  \
$\No_{W^e} (h) \in \C[[t_1, \dots, t_r]] \cap
{\Q(t_1, \dots, t_r)} \subset \Q[[t_1,\dots,t_n]]$.

\medskip

The given  data ---the description of the fiber and its local equations--- suffices
 to determine $W^e$
uniquely, and in particular    allows to compute a rational
function $q$  which approximates the norm  $\No_{W^e}(h)$ to any
given precision $\kappa$ (we understand by this  that both Taylor
expansions coincide up to degree $\kappa$, that is $\No_{W^e}(h)
\equiv q  \ \ \mod{(t_1, \dots, t_r)^{\kappa +1 }}$). The rational
function $q$ can be obtained by a procedure based on an iterative
application of the Newton operator. This operator, defined as
$$
\Newton_F^{\, t}:= x^t - \cJ_F(x)^{-1} \, F(x)^t  \ \in
K(x)^{n\times 1}, $$ enables us to approximate the points in
$W^e$ from the points in the fiber $Z$. If we set
$\Newton_F^{(m)} \in K(x)^{1 \times n} $ for the $m$-times
iteration of $\Newton_F$, then, for every $\xi \in Z$, $$
\Newton_F^{(m)}(\xi) \equiv \gamma_\xi \ \ \mod {(t_1, \dots,
t_r)^{2^m }}$$ (see \cite[Section 2]{HeKrPuSaWa00}).

Procedure $\NumDenNewton$ (Subroutine \ref{NumDenNewton})
 computes polynomials $g_1^{(m)},\dots, g_n^{(m)}, f_0^{(m)}$ in $\Q[t_1,\dots,t_r]$
  such that
$$ \Newton_F^{(m)} = (g_1^{(m)}/  f_0^{(m)},\dots, g_n^{(m)}/ f_0^{(m)}).$$
Herein, $\Homog(f,d)$ is a procedure which computes the homogeneization of the polynomial
$f$ up to degree $d \ge \deg (f)$, $\JacobianMatrix(F,x)$ is a procedure which
 constructs the Jacobian matrix with respect to the variables $x$ associated
  to the system of polynomials $F$ and $\Adjoint(M)$ is a procedure which computes
   the adjoint of the matrix $M$.
For the correctness and complexity of the whole procedure, see
\cite[Lemma 30]{GiHaHeMoPaMo97}.

\newpage

\begin{algo}{Computation of numerators and denominators for the Newton operator}
{NumDenNewton}
{$\NumDenNewton(F, n, x, d, m)$}

\hspace*{4.1mm}{\# $F_1, \dots, F_n \in A[x_1,\dots,x_n]$ such
 that $J_F (x)\not\equiv 0${;}

\hspace*{4.1mm}\# $n$ is the number of dependent variables $x$,

\hspace*{4.1mm}\# $d$ is an upper bound for the degrees of the polynomials $F_1, \dots, F_n$,

\hspace*{4.1mm}\# $m$ is the number of iterations to be computed.

\smallskip

\hspace*{4.1mm}\# The procedure returns $g^{(m)}_1, \dots,
g^{(m)}_n, f^{(m)}_0$ such that $\cN_F(x)^{(m)} =
(g^{(m)}_1/f^{(m)}_0, \dots, g^{(m)}_n/f^{(m)}_0)$.

\begin{enumerate}

\item $\cJ_F := \JacobianMatrix(F, x)${;} \vspace{-1.5mm}

\item $\Delta_F:= |  (\cJ_F) |${;} \vspace{-1.5mm}

\item $A:= \Adjoint(\cJ_F)${;} \vspace{-1.5mm}

\item $\nu := nd +1${;} \vspace{-1.5mm}

\item {\bf for } $i$ { \bf from } 1 { \bf to } $n$ { \bf do } \vspace{-1.5mm}

\item \hspace*{4.1mm} $g_i^{(1)} := \Delta_F x_i - \sum_{j=1}^n A_{ij} f_j ${;} \vspace{-1.5mm}

\item \hspace*{4.1 mm} $G_i := \Homog(g_i^{(1)}, \nu)${;} \vspace{-1.5mm}

\item{\bf od;} \vspace{-1.5mm}

\item $f_0^{(1)} := \Delta_F${;} \vspace{-1.5mm}

\item $F_0:= \Homog(\Delta_F, \nu)${;} \vspace{-1.5mm}

\item  {\bf for } $k$  { \bf from } 2 { \bf to } $m$  { \bf do } \vspace{-1.5mm}

\item \hspace*{4.1mm}  {\bf for } $i$ { \bf from } 1 { \bf to } $n$ { \bf do } \vspace{-1.5mm}

\item \hspace*{8.1mm} $g_i^{(k)}:= G_i(f_0^{(k-1)}, g_1^{(k-1)}, \dots, g_n^{(k-1)})${;} \vspace{-1.5mm}

\item \hspace*{4.1 mm} {\bf od;}\vspace{-1.5mm}

\item \hspace*{4,1 mm} $f_0^{(k)} := F_0(f_0^{(k-1)}, g_1^{(k-1)}, \dots, g_n^{(k-1)})${;} \vspace{-1.5mm}

\item  {\bf od;} \vspace{-1.5mm}

\end{enumerate}

\hspace*{7.7mm} {\bf return}($g_1^{(m)}, \dots, g_n^{(m)}, f_0^{(m)}$);
}
\end{algo}

\bigskip

We summarize the procedure that approximates the norm of a given polynomial $h$ in
Procedure {$\Norm$} (Subroutine \ref{Norm}).
Herein,  $\CompanionMatrix$ is the procedure which
constructs the companion matrix of
a given univariate  polynomial. We keep our notations.

\begin{algo}{Approximation of the Norm  }
{Norm} {$ \Norm( h, \delta, n, x, p, v , F, d, \kappa)$ }

\hspace*{4.1mm}{\# $h \in A[x_1, \dots, x_n]$ is the polynomial whose norm we want to approximate,

\hspace*{4.1mm}\# $\delta$ is an upper bound for the
degree of $h$,

\hspace*{4.1mm}\# $n \in \N $ \ is the number of  dependent variables $x$,

\hspace*{4.1mm}\# $p \in \Q[t]$, \ $v  \in \Q[t]^n$   \ is a given
geometric resolution of $Z$,

\hspace*{4.1mm}\# $F= (F_1, \dots, F_n)  $ is a vector of polynomials in $I(W)$
\ such that $\Delta_F(0, \xi) \ne 0$ for all $\xi \in Z$,

\hspace*{4.1 mm}\# $d$ is an upper bound for
the degrees of the polynomials $F_1, \dots, F_n$,

\hspace*{4.1mm}\# $\kappa \in \N$ \ is the desired level of precision.

\smallskip

\hspace*{4.1mm}\# The procedure returns  $f, g \in A$ with $f(0) \ne 0$
 such that $g/f $ approximates the norm $\No_{W^e}(h)$

\hspace*{4.1mm}\# with precision  $\kappa$.

\begin{enumerate}

\item $m:=\lceil \log_2  (\kappa + 1) \rceil $;\vspace{-1.5mm}

\item $(g_1, \dots, g_n, f_0):= \NumDenNewton(F, n, x, d, m)${;}\vspace{-1.5mm}

\item $M_p := \CompanionMatrix(p)${;}\vspace{-1.5mm}

\item {\bf for } $i$ { \bf from } 1 { \bf to } $n$ { \bf do }\vspace{-1.5mm}

\item \hspace*{4,1mm} $M_i := g_i(v(M_p))${;}\vspace{-1.5mm}

\item {\bf od;}\vspace{-1.5mm}

\item $M_0 := f_0 (v(M_p))${;}\vspace{-1.5mm}

\item $H:= \Homog(h, \delta)${;}\vspace{-1.5mm}

\item $M:= H(M_0, M_1, \dots, M_n)${;}\vspace{-1.5mm}

\item $f:= | M_0 |^\delta${;}\vspace{-1.5mm}

\item $g:= | M  |${;}\vspace{-1.5mm}

\end{enumerate}

\hspace*{7.7mm} {\bf return}($f, g$);
}
\end{algo}

\begin{lem}  \label{newton}
Let notations be as before.
Assume that $h, F_1, \dots, F_n \in A[x_1, \dots, x_n]$ are
polynomials
encoded by slp's of length bounded by $L$ such that $\deg h \le \delta$ and
$\deg (F_i)\le d$, $1 \le i \le n$.

Then $\Norm$ (Subroutine \ref{Norm}) computes $f, g \in A$ with $f(0)
\ne 0$  such that $g/f $ approximates $\No_{W^e}(h) $ with precision
$\kappa$,  within complexity
$\cO((\log_2 \kappa ) \,  n^7  \delta^2 d^2\, D^4 L)$.
\end{lem}

\begin{proof}{Proof.--}

For the correctness of the algorithm we refer to
\cite[Section 2]{HeKrPuSaWa00} and the arguments given there.
 Now, we estimate its complexity:

\smallskip

First, the complexity of Subroutine \ref{NumDenNewton} applied to our situation is
of order
$\cO ((\log_2 \kappa) \, n^7 d^2 L)$  (see \cite[Lemma 30]{GiHaHeMoPaMo97}).

Then, the algorithm computes the matrices $v_j(M_p)$ $(1 \le j \le n)$ with complexity
of order
$\cO(nD^3)$ (note that, as the companion matrix is {\em very} sparse, the multiplication
 by $M_p$ can be done with
complexity $\cO(D^2)$). Now, the matrices $M_i : = g_i(v(M_p))$
$(1 \le i \le n)$ and $M_0 := f_0 (v(M_p))$ are obtained within
complexity $\cO ((\log_2 \kappa) \,  n^7 d^2 D^3 L)$. As $h$ is
encoded by a slp of length $L$, its homogeneous components up to
degree $\delta$ are encoded by   slp's of length $\cO ( \delta^2
L)$. Therefore, the complexity of the computation of $M$ is of
order $\cO ( \delta^2 L D^3 + (\log_2 \kappa) \, n^7d^2D^3L)$.

Finally, $f$ and $g$ can be computed  within complexity $\cO(D^4
+(\log_2 \delta) \, D^3 + (\log_2 \kappa) \, n^7 d^2 D^3 L)$ and
$\cO (D^4 + \delta^2 D^3 + (\log_2 \kappa) \, n^7 d^2 D^3 L)$
respectively.
\end{proof}


\typeout{A product formula}

\subsection{A product formula} \label{A product formula}


In what follows, we establish a product formula for the Chow form
of an affine variety. This formula is an analogue of the
classical Poisson formula for the resultant (\cite[Ch. 3, Thm.
3.4]{CoLiOs98}). It describes, under certain conditions,  the
Chow form in a recursive manner.

\medskip

Let $V \subset \A^n$ be an equidimensional {\em affine}  variety
of dimension $r$ and degree $D$ which satisfies Assumption
\ref{assumption}. Let $U_0, \dots, U_{r}$ be $r+1$ groups of
$n+1$ variables each and let $L_i : = U_{i0} + U_{i1}\, x_1+\cdots
+ U_{in}\, x_n$, $0 \le i \le r$, be the affine linear forms
associated to these groups of variables.  Set $K:= \Q(U_0, \dots,
U_{r-1})$ and  let  $I(V)^\e$ denote the extension of the ideal
of $V$ to the polynomial ring $K[x_1, \dots, x_n]$ (or to  any
other ring extension of $\Q [x_1,\dots,x_n]$ which will be clear
 from the context).
We also set \ $V^0 := V(I(V)^\e ) \cap V(L_0, \dots, L_{r-1})  \
\subset \A^n(\overline{K}) $, which is a 0-dimensional variety of
degree $D$.

\smallskip

For $0\le i \le  r$, we set  $ V_i :=V \,\cap
V(x_{i+1},\dots,x_r)  \subset \A^n $, which
 is an
equidimensional variety of dimension $i$ and degree $D$
as $V$ satisfies
Assumption \ref{assumption}.
Observe that these varieties satisfy Assumption \ref{assumption} as well.

Let
$K_i:= \Q(U_{0}, \dots, U_{i-1}) \hookrightarrow K$ and set
$$
V_i^0 :=  V(I(V_i)^\e) \cap V(L_{0},\dots, L_{i-1}) \ \subset \A^n(\overline
{K_i}).
$$
Observe  that $V_i^0 $ is also a 0-dimensional variety of degree
$D$.

Under these notations we have that $V_0^0 = V_0$, $K_r = K$ and
$V_r^0 = V^0$.

\begin{prop}{(Product formula)} \label{forma de chow}
Let $V \subset \A^n$
 be an
equidimensional variety of dimension $r$ which satisfies
 Assumption
\ref{assumption}.
Let notations be as in the previous paragraph.
Then
$$
\Ch_V(U_0, \dots, U_r)= \frac{\prod\limits_{i=0}^{r} \Ch_{V_i^0}(
{U_i})}{\prod\limits_{i=1}^{r} \Ch_{V_i^0}({e_{i}})} \ \ \in
{\Q(U_0, \dots, U_{r-1})}[U_r].
$$
\end{prop}

The proof of this fact
is based on the following lemma:

\begin{lem} \label{poisson}
Let $V  \subset \A^n$ be an
equidimensional variety of dimension $r$.
Let $\cF_V \in
\Q[U_0,\dots, U_r]$ and $\cF_{V^0} \in K[U_r]$ be
Chow forms of $V$ and $V^0$ respectively.
Then there exists $\lambda \in K^* $
such that
$$
\cF_V =
\lambda \, \cF_{V^0}.
$$
\end{lem}

\begin{proof}{Proof.--}

As before, we denote by $I(V)^\e$ the extension of the ideal
$I(V)$ to a ring extension of $ \Q[x_0, \dots, x_n]$ which will
be clear from the context. Let   $U_i^\lin$,  $0\le i \le r$,
denote the group  of $n$ variables $U_i \setminus \{ U_{i \,
0}\}$. We consider the map $$ \Q[U_0, \dots, U_r][x_1, \dots,
x_n]/ (I(V)^\e +(L_0, \dots, L_r)) \to \Q[U_0^\lin , \dots,
U_r^\lin  ][x_1, \dots, x_n] / I(V)^\e $$ defined by $ U_{i \, 0}
\mapsto -( U_{i \, 1} \, x_1 + \dots  + U_{i\, n} \, x_n),\ U_{i
 j} \mapsto  U_{i  j}\,$ and $\, x_j \mapsto  x_j$ for $ 0\le i \le
r$, $1 \le j \le n$.

As it is a ring isomorphism,
$I(V)^\e + (L_0, \dots, L_r)$ is a radical ideal.
Following  notations in
 Subsection \ref{The Chow form of a quasi-projective
  variety}, it follows that
this ideal  coincides   with the defining ideal  of the incidence variety
$\Phi_V$ and, therefore,
$$
(\cF_V) = (I(V)^e+ (L_0, \dots, L_r)) \, \cap \, \Q[U_0, \dots, U_r].
$$
Similarly
$$
(\cF_{V_0}) = (I(V^0)^e+ (L_r)) \, \cap \, K[U_r].
$$

We have that $I(V)^e + (L_0, \dots, L_r) \subset I(V^0)^e + (L_r) $ and so
$(\cF_V) \subset (\cF_{V_0})$, that is, there exists $\lambda \in
K[U_r] \setminus \{ 0\} $
such that
$\cF_V = \lambda \, \cF_{V^0}$.
As $\deg_{U_r} \cF_V = \deg V= \deg \cF_{V^0} $, $\lambda$ is an element in $
K^*$.
\end{proof}

\begin{proof}{Proof  of Proposition \ref{forma de chow}.--}
Let $1\le i\le r$. From Lemma \ref{poisson}, there exists
$\lambda_i\in K_i^*  $ such that
\begin{equation} \label{1}
 \Ch_{V_i}(U_0,\dots,U_i) =
\lambda_i \ \Ch_{V_i^0}(U_i).
\end{equation}
Hence
\ $
\Ch_{V_i}(U_0, \dots,U_{i-1}, e_i) =
\lambda_i \, \Ch_{V_i^0}(e_{i})
.$

Now, it is easy to see that
$ \Ch_{V_{i-1}}(U_{0}, \dots,U_{i-1}) $ divides
$Ch_{V_i}(U_0, \dots,U_{i-1}, e_i)
$. From Assumption \ref{assumption}, it follows
that $\deg V_{i-1} = \deg V_i = D$ and, therefore, both polynomials have the same
degree.
Moreover, the normalization imposed to both Chow forms
implies that they coincide. So
\begin{equation} \label{2}
\Ch_{V_{i-1}}(U_0, \dots, U_{i-1} )
 = \lambda_i \,  \Ch_{V_i^0}(e_{i})
.
\end{equation}
{}From Identities (\ref{1}) and (\ref{2}) we deduce that
\begin{equation} \label{3}
\frac{\Ch_{V_i}(U_0,\dots,U_i) }{\Ch_{V_{i-1}}(U_{0},\dots,U_{i-1})}
= \frac{\Ch_{V_i^0}(U_i)}{\Ch_{V_i^0}(e_{i})}.
\end{equation}

Multiplying these identities  for $i=1, \dots, r$ we obtain
$$
\frac{\Ch_V(U_0,\dots,U_r)}{\Ch_{V_0}(U_0)}  = \prod_{i=1}^{r}
\frac{\Ch_{V_i}(U_0,\dots,U_i)
}{\Ch_{V_{i-1}}(U_{0},\dots,U_{i-1})} =
\prod_{i=1}^{r}\frac{\Ch_{V_i^0}(U_i)}{\Ch_{V_i^0}(e_{i})}
$$
which gives the formula stated in Proposition \ref{forma de chow}.
\end{proof}

Observe that $ \Ch_{V_i^0}(U_i) = \prod_{\gamma \in V_i^0} L_{i}
(\gamma) = \No_{V_i^0} (L_i)$ and  \ $
 {\Ch_{V_i^0}(e_{i})} =  \prod_{\gamma \in V_i^0} x_{i} (\gamma)
= \No_{V_i^0} (x_{i})$, which implies that the Chow form of $V$
can also be presented as the quotient of two products of norms of
polynomials:

\begin{cor}\label{forma de chow - bis}
Let notations and assumptions be as before. Then
$$
\Ch_V(U_0,\dots, U_r) = \frac{\prod\limits_{i=0}^{r} \No_{V_i^0}(
{L_i})}{\prod\limits_{i=1}^{r}
 \No_{V_i^0}({x_{i}})} \ \ \in {\Q(U_0, \dots, U_{r-1})}[U_r].
$$
\end{cor}

This formula enables us to compute $\Ch_V$ as a quotient of power
series. To do so, we are going to prove a technical lemma first.

\begin{lem}
\label{ideal radical} Let $V \subset \A^n$ be an equidimensional
variety of dimension $r$ which satisfies Assumption
\ref{assumption}. Assume that $V$ is Cohen-Macaulay at every
point of $Z:= V \cap V(x_1, \dots, x_r)$. Then, the ideal $I(V)
+(x_1, \dots, x_r) \subset \C[x_1, \dots, x_n]$ is radical.

\end{lem}

\begin{proof}{Proof.--}

Let $\overline{V} $ denote the projective closure of $V \subset \A^n \hookrightarrow \P^n$.
Let $Z := {V} \cap V(x_1, \dots, x_r)$. The fact that $\# Z = \deg V $ implies that
$Z= \overline{V} \cap V(x_1, \dots, x_r)$.

Take $\xi \in Z$ and let
$Q_\xi $
be the primary component of the ideal
$ I(V)+(x_1, \dots, x_r) \subset \C [x_1, \dots, x_n]$
which corresponds to
$\xi$.

We consider the length $\ell(\overline{V}, V(x_1, \dots, x_r);\,
\xi)$ which under our assumptions can be defined as
\begin{equation}\label{length}
\ell(\overline{V}, V(x_1, \dots, x_r); \xi)= \dim_{\C} \C [x_1,
\dots, x_n] / Q_\xi .
\end{equation}

By a suitable version of B{\'e}zout theorem (see \cite[Prop.
3.30]{Vogel84})
$$
\sum_{\xi \in Z} \ell(\overline{V}, V(x_1, \dots, x_r); \xi) \le
\deg V.
$$

On the other hand, as $\ell(\overline V, V(x_1, \dots, x_r);\,
\xi)$ is a positive integer for each $\xi \in Z$ and $\# Z = \deg
V$, it follows that
$$\sum_{\xi \in Z} \ell(\overline{V}, V(x_1, \dots, x_r); \xi) \ge
\deg V.$$

Then $ \ell(\overline{V}, V(x_1, \dots, x_r); \, \xi) = 1 $ for
all $\xi \in Z$, and so (\ref{length}) implies that $Q_\xi =
(x_1- \xi_1, \dots, x_n - \xi_n)$  which is a prime ideal.

As $I:= I(V)+ (x_1, \dots, x_r)$ is 0-dimensional it has no
embedded components. Hence $I   = \cap_\xi  \, Q_\xi  $ is a
radical ideal.
\end{proof}

\smallskip
The following corollary shows that the coordinates of all the
points in $V_i^0 $
 belong to the
subring $\C[[U_{0}-e_{1}, \dots, U_{i-1}-e_i]] \cap
\overline{K_i}$, and that there is a one-to-one correspondence between
the points of $Z := {V} \cap V(x_1, \dots, x_r)$ and the points of
$ V_i^0 $.

\begin{cor}\label{series}
Let notations and assumptions be as in Lemma \ref{ideal radical}
and before. Let $0 \le i \le r$  and  $\xi  \in Z$. Then there
exists a unique $\gamma^{(i)}_\xi  \in \C[[U_{0}-e_{1}, \dots,
U_{i-1} -e_i]]^n$ such that $ \gamma^{(i)}_\xi \in V_i^0$ and $
\gamma^{(i)}_\xi (e_{1}, \dots, e_i) = \xi $.
\end{cor}

\begin{proof}{Proof.--}

Suppose $I(V)$ is generated by the polynomials $h_1, \dots, h_t$.
Since we are in the conditions of the previous lemma, the  Jacobian criterion
\cite[Thm. 18.15]{Eisenbud95}
implies that the
Jacobian matrix associated to the generators $h_1, \dots, h_t, x_1, \dots, x_r$ of the ideal
   $I(V) + (x_1, \dots, x_r)$
has maximal rank $n$ at $\xi$. In other words, there are $n$
polynomials $g_1, \dots, g_n$ among $h_1, \dots, h_t, x_1, \dots,
x_r$ such that the associated Jacobian determinant is non-zero.
Now, as the rank of the Jacobian matrix of $h_1,\dots,h_t$ at
$\xi$
 is bounded by
 the codimension $n-r$ of $V$ at $\xi$, we can assume
 w.l.o.g. that $g_1:=x_1, \dots, g_r:=x_r$.

Let $$ \Delta:= \left| \, \left(\frac{ \partial g_i }{ \partial
x_j}\right)_{r+1 \le i,  j \le n}
    \right|  
$$
be the Jacobian determinant of $g_{r+1}, \dots, g_n$ with respect to the
variables $x_{r+1}, \dots, x_n$. Then $\Delta(\xi ) \ne 0$
since   $\Delta$ coincides with  the Jacobian determinant of the system
$g_1, \dots, g_n$.

On the other hand, let $\Delta_i \in K_i[x_1, \dots, x_n]$ denote the
Jacobian determinant of the system
\newline  $ L_0, \dots, L_{i-1}, x_{i+1}, \dots, x_r,
g_{r+1}, \dots,
g_n $.
An easy verification shows that
   $\Delta_i (e_1, \dots, e_{i-1})(\xi) = \Delta(\xi) \ne 0$.
The statement follows from
 the arguments in Subsection \ref{Newton's algorithm}.
\end{proof}

\bigskip
Now, set
$$
\Psi:= \prod_{i=0}^r \No_{V_i^0}(L_i ) \ \in K[U_r], \quad \quad  \Phi:= \prod_{i=1}^{r} \No_{V_i^0}(x_{i} ) \ \in K^*
,
$$
so that, by Corollary \ref{forma de chow - bis}, \ $\Ch_V:= \Psi/
\Phi$.

\smallskip
{}From Corollary  \ref{series},   $\Psi \in \Q[[U_{0}-e_{1},
\dots, U_{r-1}-e_{r}]][U_r]$ and  $\Phi \in \Q[[U_{0}-e_{1},
\dots, U_{r-1}-e_{r}]]$.

\smallskip
The following lemma gives
the order of the denominator $\Phi$ at $E:= (e_1, \dots, e_r) \in \A^{r \, (n+1)}$
together with its  graded component of lowest degree:

\begin{lem} \label{orden de Gamma}
Let notations be as in the previous paragraph and let $D:= \deg
V$. Then $\ord_E (\Phi) = r  D$ and its graded component of degree
$r D$ is
$$
\Phi_{r  D} =
\pm \, \prod_{i=1}^r \Ch_{V_0}(U_{i-1}).
$$

\end{lem}

\begin{proof}{Proof.--}
Clearly, $\ord_E (\Phi) =\sum_{i=1}^r \ord_E
\Big(\No_{V_i^0}(x_i)\Big)$.

\smallskip
Let $1\le i \le n$.  Recall that $\No_{V_i^0}(x_i)=
\Ch_{V_i^0}(e_i)$. From Identity (\ref{3}) in the proof of
Proposition \ref{forma de chow} we have:
 $$
{\Ch_{V_i^0}(e_{i})} \, \Ch_{V_i}(U_0, \dots, U_i)
= \Ch_{V_i^0}(U_{i}) \ \Ch_{V_{i-1}}
(U_{0}, \dots, U_{i-1})
.$$

As $\Ch_{V_i^0}(e_{0}) =1$, then
\ $ {\Ch_{V_i^0}(e_{i})} \ \Ch_{V_i}(U_0, \dots, U_{i-1}, e_0)
=  \Ch_{V_{i-1}}
(U_{0}, \dots, U_{i-1})$.

We also have that
$$\Ch_{V_i}(e_{1},  \dots, e_{i}, e_0)=
\pm \Ch_{V_i}(e_{0}, e_{1},  \dots, e_i) = \pm 1.
$$
This shows that $ \Ch_{V_i}(U_0, \dots, U_{i-1},e_0) $ is
invertible in $ \Q[[U_0 - e_1, \dots, U_{i-1} - e_i]]$ and,
therefore, if $m := \ord_E \Big({\Ch_{V_i^0}(e_{i})}\Big)$,

$$
{\Ch_{V_i^0}(e_{i})}\ \equiv \  \pm \Ch_{V_{i-1}}
(U_{0}, \dots, U_{i-1}) \
 \  \ \ \mod{(U_{0}-e_{1}, \dots, U_{i-1}
-e_i)^{m+1}}.
$$
By Lemma \ref{poisson}, there exists $\lambda_{i-1} \in \Q(U_0, \dots, U_{i-2}) \setminus \{ 0 \}$ such that
$$
\Ch_{V_{i-1}}(U_{0}, \dots, U_{i-1}) = \lambda_{i-1} \,
\Ch_{V_{i-1}^0} (U_{i-1})
.$$

As $\Ch_{V_{i-1}^0} (U_{i-1})$ is a homogeneous polynomial of degree D in the group of variables
$U_{i-1}$ and does not depend on
$U_{i-1\, i}$, it is also homogeneous as a polynomial expanded in $U_{i-1} - e_i$.
Then, the order of $\Ch_{V_{i-1}}$ at
  $e_{i}$
with respect to the group of variables
$U_{i-1}$ equals
  $D$.
On the other hand, we have that
\ $ \Ch_{V_{i-1}}(e_1, \dots, e_{i-1}, U_{i-1})=
\pm \, \Ch_{V_0}(U_{i-1}) \ne 0$.  This implies that the series $\Ch _{V_{i-1}}$ in $  \Q[[U_0 - e_1, \dots, U_{i-1} - e_i]]$ has a term of degree $D$ depending only on the group of variables $U_{i-1} - e_i$.
We conclude that $m = D$ and
$$
(\No_{V_i^0}(x_i))_D=({\Ch_{V_i^0}(e_{i})})_D = \pm ( \Ch_{V_{i-1}}(U_{0}, \dots, U_{i-1}))_D = \pm \, \Ch_{V_0}(U_{i-1}).
$$

Therefore, $\ord_E (\Phi) = \sum_{i=1}^r \ord_E
\Big({\No_{V_i^0}(x_i)}\Big)= r D$ and the graded part of lowest
degree of $\Phi$ is $ \Phi_{r D} = \prod_{i=1}^r
(\No_{V_i^0}(x_i))_D = \pm \, \prod_{i=1}^r \Ch_{V_0}(U_{i-1})$.
\end{proof}


\typeout{The algorithm}

\subsection{The algorithm}

\label{The algorithm}


Here, we are going to put the previous results together in order
to obtain the algorithm underlying Main Lemma \ref{fibra} and to
estimate its complexity.

\medskip

Let notations be as in Main Lemma \ref{fibra}. As we have
already  noted, the imposed conditions  imply that both $V$ and
$Z$ have no component in the hyperplane $\{ x_0 =0\}$. Hence $V$
equals the projective closure of its affine part $V_{x_0}:= V
\setminus \{x_0=0\}$ and so their both Chow forms coincide. Hence
we concentrate w.l.o.g. on the affine case. We  use affine
coordinates and  keep the notation of the previous subsection.

\smallskip
{}From Corollary \ref{forma de chow - bis}, we have that $$ \Ch_V
= \frac{\prod\limits_{i=0}^{r} \No_{V_i^0}( {L_i})}
{\prod\limits_{i=1}^{r} \No_{V_i^0}({x_{i}})}. $$

Now, we approximate the norms appearing in this formula.

\smallskip

Set
$$
\cV_i:= V(I(V_i)^\e) \cap V(L_0, \dots, L_{i-1})) \ \subset \A^{i
  (n+1)} \times \A^n .
$$
The map  \ $ \pi_i : \cV_i \to \A^{i
  (n+1)} $ defined by $(U, x) \mapsto U$ is {dominant} of degree $D: = \deg V$.
We set $Z:= V_0 = V \cap V(x_1, \dots, x_r) \subset \A^n$ and let $E_i:= (e_1, \dots,
e_i) \in \A^{i   (n+1)} $.
Then
$$
\cZ_i:= \pi_i^{-1} (E_i) = \{ E_i\} \times Z
$$
and so this fiber is a 0-dimensional variety of cardinality $D$.
Furthermore,
 it is easy to check that
$$
L_0, \dots, L_{i-1}, \, x_{i+1}, \dots, x_r, \, f_{r+1}, \dots,
f_n \ \in  \Q[U_0, \dots, U_{i-1}][x_1, \dots, x_n]
$$
is a system of local equations of $\cV_i$ at $\cZ_i$.

\smallskip
Since by definition, $\No_{\pi_i}(x_i)$ and $\No_{\pi_i}(L_i)$
coincide with $\No_{V_i^0} (x_i) \in \Q[[U_0 -e_1,
\dots, U_{i-1} - e_i]]$ and $\No_{V_i^0} (L_i) \in \Q[[U_0 -e_1,
\dots, U_{i} - e_{i+1}]]$ respectively, we can compute
 any approximation of the latter applying Procedure $\Norm$
 (Subroutine \ref{Norm})
 modulo a change of variables
$(U_0, \dots, U_{i}) \mapsto (\widetilde U_0 +e_1,
\dots, \widetilde U_{i} + e_{i+1})$ (in order to center the series at
$0$).

We multiply the computed approximations for $ 0\le i \le r$ to
obtain rational functions $\psi$ and $\varphi$ which approximate
the power series
$$
\Psi:= \prod_{i=0}^r \No_{V_i^0}(L_i ) \ \in \Q [[ U_0 - e_1,
\dots, U_{r-1} -e_r]] [U_r], \quad \ \Phi:= \prod_{i=1}^{r}
\No_{V_i^0}(x_{i} ) \ \in \Q [[ U_0 - e_1, \dots, U_{r-1} -e_r]]
^*
$$
respectively.

{}From these approximations, we compute the graded parts  of
$\Phi$ and $\Psi$ of  degrees between $r  D$ and  $(2 \, r +1) \,
D$  centered at $( E,0) \in \A^{(r+1)  (n+1)}$ (where $E:=  (e_1,
\dots, e_r) \in \A^{r  (n+1)}$)
 by
applying
Procedure $\GradedParts$ (see Subsection \ref{Effective division procedures}).
\smallskip

By Lemma \ref{orden de Gamma}, we have that $\ord_{(E,0)} (\Phi) =
\ord_E (\Phi) = r D$. We also have $\deg \Ch_V =(r+1)\, D$. We use
this information  together with the obtained graded parts in
order to
 apply Procedure $\PowerSeries$ (Subroutine \ref{PowerSeriesDivision}).
This yields   a
polynomial $Q \in \Q[U_0, \dots, U_r]$ such that
$$
Q = \Phi_{r  D}^{(r+1) \, D + 1} \, \Ch_V.
$$
Again, from Lemma \ref{orden de Gamma},
the denominator $\Phi_{r \, D}^{(r+1) \, D + 1} $
does not vanish at
 $E^0:= (e_0, \dots, e_0) \in \A^{r  (n+1)}$.
We apply  Procedure $\PolynomialDivision$ (Subroutine
\ref{PolynomialDivision}) to the polynomials $ Q$ and $ \Phi_{r
D}^{(r+1) \, D + 1} $ and the point $E^0$.

\smallskip

We summarize this procedure in Procedure $\ChowForm$ (Subroutine \ref{algo-fibra})
 which computes
the Chow form of an affine equidimensional variety $V$ satisfying Assumption
 \ref{assumption}.

\begin{algo}
{Chow form from a fiber and local equations} {algo-fibra} {$
\ChowForm(n,x, r,D , p, v, f, d)$ }

\hspace*{4.1mm}{\# $n$ is the number of variables $x: = (x_1, \dots, x_n)$,

\hspace*{4.1mm}\# $r, D$ are the dimension and the degree of $V$
respectively,

\hspace*{4.1mm}\# $p \in \Q[t]$, \ $v  \in \Q[t]^n$ \ is a given
geometric resolution of  the fiber $Z$,

\hspace*{4.1mm}\# $ f = (f_{r+1}, \dots, f_n) \in \Q[x_1, \dots,
x_n]^{n-r} $  is a system of
local equations of $V$ at $Z$ of degrees

\hspace*{4.1mm}\#  bounded by $d$.

\smallskip

\hspace*{4.1mm}\# The procedure returns the normalized Chow form $\Ch_V$.

\begin{enumerate}

\item  {\bf for } $i$ { \bf from } 1 { \bf to } $r$ { \bf do}
\label{do}\vspace{-1.5mm}

\item  \hspace*{3mm} $( \varphi_i^{(1)}, \varphi_i^{(2)}):= \Norm(x_i, 1,n, x, p, v,
L_0, \dots, L_{i-1}, x_{i+1}, \dots, x_r, f_{r+1}, \dots, f_n, d,
(2\, r +1)\, D) $;\vspace{-1.5mm}

\item {\bf od}; \label{od1}\vspace{-1.5mm}

\item  {\bf for } $i$ { \bf from } 0 { \bf to } $r$ { \bf do}\vspace{-1.5mm}

\item  \hspace*{3mm} $(\psi_i^{(1)}, \psi_i^{(2)}):= \Norm(L_i, 1,n, x, p, v,
L_0, \dots, L_{i-1}, x_{i+1}, \dots, x_r, f_{r+1}, \dots, f_n, d, (2\, r +1)\, D)
 $;\vspace{-1.5mm}

\item {\bf od}; \label{od}\vspace{-1.5mm}

\item $\varphi^{(1)} := \prod_{i=1}^r \varphi_i^{(1)}$, \quad \quad $\varphi^{(2)}
 := \prod_{i=1}^r \varphi_i^{(2)}$; \label{multiplicacion}\vspace{-1.5mm}

\item $\psi^{(1)} := \prod_{i=0}^r \psi_i^{(1)}$, \quad \quad $\psi^{(2)}
:= \prod_{i=0}^r \psi_i^{(2)}$;
\label{multiplicacion2}\vspace{-1.5mm}

\item $(\Phi_{0} , \dots, \Phi_{(2r+1)\, D}):=
\GradedParts(\varphi^{(1)}, \varphi^{(2)}, (e_1,\dots, e_r), (2\,
r+1) \, D)$; \label{gradedparts}\vspace{-1.5mm}

\item $(\Psi_{0} , \dots, \Psi_{(2r+1)\, D}):=
\GradedParts(\psi^{(1)}, \psi^{(2)}, \, (e_1,\dots, e_r,0), (2\,
r+1) \, D)$; \label{gradedparts2}\vspace{-1.5mm}

\item $Q:= \PowerSeries((r+1)\, (n+1), r\, D, (r+1)\, D, \,
\Phi_{r\, D} , \dots, \Phi_{(2r+1)\, D}, \, \Psi_{r\, D} , \dots,
\Psi_{(2r+1)\, D})$; \label{powerseries}\vspace{-1.5mm}

\item $\Ch_V:= \PolynomialDivision(Q, \Phi_{r \, D}^{(r+1) \, D+1} ,
 (r+1) D, (e_0,\dots, e_0))$; \label{polydivision}\vspace{-1.5mm}

\end{enumerate}

\hspace*{7.7mm} {\bf return}($\Ch_V$);
}
\end{algo}

\begin{proof}{Proof of Main Lemma \ref{fibra} .--}
As we have already observed,  we may suppose w.l.o.g that $V$ is an affine variety and
that the  polynomials $f_{r+1}, \dots, f_n$ are in $\Q[x_1, \dots, x_n]$.

We apply  Procedure $\ChowForm$ (Subroutine \ref{algo-fibra}) to $V$  in order
 to compute its normalized Chow form. The correctness of this
procedure follows from our previous analysis. The announced
complexity is a consequence of the complexity of the subroutines
we call during this procedure:
\begin{itemize}

\item By Lemma \ref{newton}, the complexity of  lines \ref{do} to
\ref{od} is of order  $\cO(r \, \log_2 (r \,  D) \, n^7 \, d^2 \,
D^4 L)$. The products in lines \ref{multiplicacion} and
\ref{multiplicacion2} do not change this estimate.

\item The computation of the graded parts in lines
  \ref{gradedparts} and \ref{gradedparts2} has complexity  $\cO(r^3 \,
  \log_2 (r \,  D) \, n^7 \, d^2 \, D^6 L) $.

\item Finally, the subroutines $\PowerSeries$ and $\PolynomialDivision$ in lines
  \ref{powerseries} and \ref{polydivision} add complexity
$\cO(r^8 \, \log_2 (r \,  D) \, n^7 \, d^2 \, D^{11} L)$.
\end{itemize}

We conclude that the overall complexity is \ $\cO(r^8 \, \log_2
(r \,  D) \, n^7 \, d^2 \, D^{11} L) $.
\end{proof}

\begin{rem} \label{dominante}
The Chow form $\Ch_V$ is  the numerator of
$\No_{V^0}(L_r)$ by Lemma \ref{poisson}.
Unfortunately this norm is  a rational
function, due to the fact that the map
$\pi_r$ is not finite but just dominant.

The product formula is the tool which allows  to overcome this
difficulty, as it gives an expression for $\Ch_V$  without any
extraneous
denominator.
\end{rem}

\smallskip

We directly derive
the following estimate for the length of a slp representation
of the Chow form of an equidimensional variety:

\begin{cor} \label{L(Ch)}
Let $V \subset \P^n$ be an equidimensional variety of dimension
$r$ and degree $D$. Let $ f_{r+1}, \dots, f_{n}\in I(V)$ be a
system of local  equations at a dense open subset of $V$, encoded
by slp's of length bounded by $L$. Then, if $d:= \max\{ \deg(f_i)
: r+1 \le i\le n\}$, we have
$$
L(\cF_V) \le
  \cO (r^8 \,
  \log_2 (r \,  D) \, n^7 \, d^2 \, D^{11} L).
$$

\end{cor}

\begin{proof}{Proof.--}
Let $\ell_1, \dots, \ell_r \in \Q[x_0, \dots, x_n]$ be
linear forms
such that
\ $ Z:= V \cap V(\ell_1, \dots, \ell_r)$ \ is a 0-dimensional variety of
cardinality $D$.
We can choose these linear forms so that
$Z$ lies in the dense open subset where  $f_{r+1}, \dots, f_n$
is a system of
local equations.

Furthermore, let
 \ $\ell_0, \, \ell_{r+1}, \dots,  \ell_n$ \ be linear forms which
complete the previous ones to a change of variables such that $Z
\cap \{ \ell_0\ = 0 \} = \emptyset$.

Then $V$ satisfies Assumption \ref{assumption} with respect to these
variables, and  the statement
follows directly from Main Lemma \ref{fibra}.
\end{proof}


%% file: ComputacionChow.tex

\typeout{The computation of the Chow form}

\section{The computation of the Chow form}

\label{The computation of the Chow form}

\vspace{2mm}

We devote this section  to  the description and complexity
analysis of the algorithm underlying Theorem 1. The first
subsections gather some results which   lead to the proof of the
theorem.


\typeout{Geometric resolutions}

\subsection{Geometric resolutions}

\label{Geometric resolutions}

\vspace{2mm}

Geometric resolutions where first introduced in the works of Kronecker
and K{\"o}nig in the last years of the   XIXth century.
Nowadays they are  widely  used in computer
algebra, especially in the 0-dimensional case, but there are also
important applications
in the general case.
We refer to \cite{GiHe00} for a complete
historical account.

In what follows we show  that
we can compute any ---sufficiently generic---  geometric resolution of an equidimensional variety from
a Chow form in polynomial time.
This computation and the procedure described in Section \ref{The representation of the Chow form} imply that, from the point of view of complexity,
Chow forms and
geometric resolutions  are
equivalent representations of an equidimensional variety.

\medskip

Let $V \subset \A^n$ be an equidimensional affine variety of
dimension $r$ and degree $D$. For $0\le i\le r$, let $L_i$ denote as usual the
generic affine  forms. Let $c_i
\in \Q^{n+1}$. We set $$ \ell_i := L_i(c_i) = c_{i   0} + c_{i
1} \, x_1+ \cdots + c_{i  n} \, x_n \ \in \Q[x_1, \dots, x_n] .$$
 We assume that the
projection map $$ \pi_{(\ell_1,\dots,\ell_r)} : V \to \A^{r} \quad \quad ,
\quad \quad x \mapsto (\ell_1(x), \dots ,  \ell_{r}(x)) $$ is {\em
finite}, that is the affine linear forms $\ell_1, \dots, \ell_r$ are in
Noether position with respect to $V$. Let  $y_1, \dots, y_r$ be
new variables. Set $$ K:= \Q ( y_1, \dots, y_r) \quad , \quad L=
\Q ( \ell_1, \dots, \ell_r) \otimes_{\Q[\ell_1, \dots, \ell_r]}
\Q[V] $$ and consider the morphism $$ K \rightarrow L \quad \quad , \quad \quad
y_i \mapsto \ell_i.$$ Then $ K \hookrightarrow L$ is a finite
extension of degree $[L:K] \le D$. We assume furthermore that
$\ell_0 $ is a primitive element of this extension, that is \ $ L
= K [\ell_0 ] $.

Then the  {\em geometric resolution of $V$ associated to
 $\ell: = (\ell_0, \dots, \ell_r)$ } is
the pair
$$
p:= p_{V,
\ell}  \in K[t] \quad \quad , \quad \quad  w:= w_{V , \ell}
\in K[t]^n
$$
where  $p$ is the monic minimal polynomial of $\ell_0$ with
respect to the extension $K \hookrightarrow L$, and  \ $w = (w_1,
\dots, w_n)$ verifies  $ \deg w_i < [L:K]$ and $ p'(\ell_0)
\, x_i = w_i  (\ell_0) \, \in L $ for $1\le i\le n$, where $p':=
\partial p /\partial t$. These polynomials are uniquely
determined and because of the Noether position assumption, we
have  that $p, w_i$ lie in fact in  $ \Q[y_1, \dots, y_r][t]$, see
e.g. \cite[Section 3.2]{GiLeSa01}.

A geometric resolution gives a parametrization of
a dense  open set of $V$ in terms of  the points of a hypersurface in $\A^{r+1}$:
there is a map
\begin{eqnarray*}
V(p(t,y_1, \dots, y_r)) \setminus V(p'(t,y_1, \dots, y_r))  & \to
& V \setminus V(p'( \ell_0(x), \ell_1(x), \dots, \ell_r(x)))\\
(t,y_1, \dots, y_r) &\mapsto& \frac{w}{p'}(t, y_1, \dots, y_r).
\end{eqnarray*}

\smallskip

Note that, in case the considered variety  is 0-dimensional, this
definition of geometric resolution essentially coincides with the
one given in Section \ref{The representation of the Chow form}:
the passage from one to the other can be made by considering the
resultant with respect to the variable $t$ between $p$ and $p'$.

\bigskip

The following construction shows that  the  geometric resolution associated
to the generic affine linear forms $L_0, \dots, L_r$ can be expressed in terms of the
characteristic polynomial of the variety,
and hence in terms of the Chow form:

Let $U_0, \dots, U_r$ be $r+1$ groups of $n+1$ variables which correspond to the coordinate functions of $\A^{(r+1)(n+1)}$ and
let
 $T:= (T_0,\dots,T_r)$  be
a group of $r+1$ variables which  correspond to the coordinate
functions of $\A^{r+1}$.
We recall that
a {\em characteristic polynomial} \ $\cP_V \in \Q[U_0, \dots, U_r][T_0, \dots,
 T_r]$  \
 of $V$  is defined as  any
defining equation of the Zariski closure of the image of the map
$$
\varphi_V:  \A^{(r+1)(n+1) }\times V \to \A^{(r+1)(n+1)}\times
\A^{r+1}, \quad (u_0, \dots, u_r ; \, \xi) \mapsto ( u_0, \dots,
u_r; \ L_0(u_0, \xi),  \dots, \ L_r(u_r, \xi))
$$
which is a hypersurface.
This is a   multihomogeneous polynomial  of degree
$D$ in each group of variables $U_i\cup\{ T_i\}$. Its degree in
the group of variables $T$ is also bounded by $D$.

\smallskip

A characteristic  polynomial of $V$  can be derived from a Chow form $\cF_V$.
For
 $1\le i\le  r$ we set
$\zeta_i:=(U_{i 0}-T_i, U_{i 1},\dots,U_{i n})$.
Then
\begin{equation} \label{P_V}
\cP_V=  (-1)^D \, \cF_{V}( \zeta_0,\dots,\zeta_r )
\end{equation}
is a characteristic polynomial of $V$.

Set $\cP_{ V}:= a_D\, T_0^D+ \cdots +a_0 $ for the expansion of $
\cP_V $ with respect to $T_0$.
Then $a_D$ lies  in $
\Q[U_1, \dots, U_r] \setminus\{ 0\} $,  and
in fact
it  coincides  with the coefficient of
$U_{00}^D$ in $\cF_V$, that is
$$
a_D(U_1, \dots, U_r) =  \cF_V(e_0, U_1, \dots, U_r).
$$

In case $V$ satisfies Assumption~\ref{assumption}, we define {\em the}
characteristic polynomial of $V$ as
$$
(-1)^D \, \Ch_{V}( \zeta_0,\dots,\zeta_r )$$ where $Ch_{V}$ is
the normalized Chow form of $V$. We refer to \cite[Section
2.3.1]{KrPaSo99} for  further  details as well as for the proof
of the stated facts.

\begin{lem} \label{resolucion generica}

Let $V \subset \A^n$ be an equidimensional variety of dimension $r$
and degree $D$. Let $U_0, \dots, U_r$
be $r+1$ groups of $n+1$ variables and let $L_0, \dots, L_r$ be the generic affine
 forms associated to
$U_0, \dots, U_r$.  Set $E:= \Q(U_0, \dots, U_r)$ and let $V^\e $
denote the Zariski
closure of
$  {V}$ in  $\A^n(\overline{E})$.  Let $T_0, \dots, T_r$ be new indeterminates.

Then the
geometric resolution of $V^\e$ associated to $L_0, \dots, L_r$
is  given by
$$
P:= \frac{\cP_V}{a_D}  \in E[T_1, \dots, T_r][T_0] \quad \quad ,
\quad \quad
W:=  - \frac{1}{a_D} \left( \frac{\partial \cP_V }{ \partial U_{0 \, 1}} , \dots,
    \frac{\partial \cP_V }{ \partial U_{0 \, n}}  \right) \in E[T_1, \dots, T_r][T_0]^n
$$
where $\cP_V$ is a characteristic polynomial of $V$ and $a_D$
is the leading coefficient of $\cP_V$ with respect to $T_0$.
\end{lem}

\begin{proof}{Proof.--}

Using the fact that
the extended ideal $I(V)^\e \subset E[x_1, \dots, x_n]$
is radical,
it is easy to check that $I(V^\e) = I(V)^\e$.
Consider then  the morphism
$$
\cA:= E[T_1, \dots, T_r] \to \cB:= E[x_1, \dots, x_n] / I(V)^\e
\quad \quad , \quad \quad
T_i \mapsto L_i(U_i, x).
$$
Our first aim is to prove that  this is an integral  inclusion,
or in other words, that  the projection map $\pi_{(L_1,\dots,L_r)} : V^\e \to
\A^r(\overline{E})$
is finite.

By definition
\begin{equation} \label{P(L)}
 P(U_0, \dots, U_r)(L_0(U_0, x), \dots, L_r(U_r, x)) \equiv 0
\ \ \ \ \mod{I(V)
\otimes_{\Q[x]} \Q(U)[x]}.
\end{equation}

Specializing $U_0$  by  the $(i+1)$-th element of the canonical
basis $e_i$  in this identity,  we deduce that \ $  P(e_i,U_1,
\dots, U_r) (T_0, T_1, \dots, T_r)  \in  \cA[T_0]$ is a monic
equation for $x_i$ for $i=1, \dots, n$. Therefore $\cA
\hookrightarrow \cB$ is an integral extension.

Set
$
\cK:=
E(T_1, \dots, T_r) $ and $
\cL:= \cK \otimes_\cA \cB$.

It is immediate   that
 $P:= \cP_V/a_D $ is a monic
  polynomial equation for $L_0$ with respect to the extension
$\cK \hookrightarrow \cL$. As  $\cA \hookrightarrow  \cB$ is an integral extension,
from the definition of $\cP_V$ we deduce that $P$ is the minimal monic polynomial of $L_0$.
This implies that  $[\cL:\cK] = D$ and that
$L_0$ is a primitive element of this extension.

\smallskip

Write
$$ Q(U,x): = \cP_V(U_0, \dots, U_r)(L_0(U_0, x), \dots,
L_r(U_r, x)) = \sum_{\beta} b_{\beta} \, U_0^{\beta}$$ with
$b_\beta \in \Q[U_1, \dots, U_r][x_1, \dots, x_n]$.

As $  b_\beta \in I(V)^\e \subset \Q[U_1, \dots, U_r][x_1, \dots,
x_n] $ for all $\beta$,
 $\displaystyle{\frac{\partial Q(U,x)}{\partial
U_{0 i}}} \in I(V)^\e \subset \Q[U_0, \dots, U_r][x_1, \dots, x_n]$ \
for  $ i=1, \dots, n$.  Therefore, $\displaystyle{\frac{\partial Q(U,x)}{\partial
U_{0 i}}} = 0 $ in $\cL$ for  $ 1\le i\le n$.
Then the  chain rule implies that the identity
$$
\frac{\partial P_V}{\partial T_0}(U,L(U, x )) \, x_i = -
\frac{\partial P_V}{\partial U_{0  i}}(U,L(U,x))
$$ holds in $\cL$ and the lemma follows.
\end{proof}

Now we show how a particular geometric resolution can be obtained by
 direct specialization of the generic one.

\smallskip
Using the same notation  as in the beginning of this subsection,
we will assume
 that
$ V \cap V(\ell_1,
\dots ,  \ell_{r})$
is a 0-dimensional variety of cardinality $D$.
This condition is satisfied provided that   $\ell_1, \dots, \ell_r$
are generic enough \cite[Prop. 4.5]{KrPaSo99}.
After a linear change of variables, we may assume w.l.o.g.  that $
  \ell_i= x_i$ for $i=1, \dots, r$, so that
the stated  condition is Assumption~\ref{assumption}.

Thus, for the rest of this section we fix the following notations:
$$Z:= V \cap V(x_1, \dots, x_r)\quad, \quad  K := \Q(x_1, \dots, x_r)\quad , \quad
L : = K
\otimes_{\Q[x_1, \dots, x_r]} \Q[V].$$

We also assume
that
$\ell_0= L_0(c_0,x)  \in \Q[x_1, \dots, x_n]$
separates the points of
$Z$.
This is also a generic condition:
if we set $\rho:= \discr_{T_0} \,  \cP_Z \in \Q[U_0] \setminus \{ 0\}$,
this  condition   is satisfied provided that
$\rho (c_0) \ne 0$.

\smallskip
These two conditions  ensure the existence of the
associated geometric resolution of $V$:

\begin{lem} \label{especializacion}

Let $V \subset \A^n$ be an equidimensional variety of dimension
$r$ and degree $D$ which satisfies Assumption~\ref{assumption}.
Let
$\ell_0: = L_0(c_0,x) \in \Q[x_1,
\dots, x_n]$ be an  affine linear form which separates the points of $Z$.

Then the projection map $\pi: V \rightarrow \A^r \ , \ \pi(x) =
(x_1, \dots, x_r)$ is finite and $\ell_0 $ is a primitive element
of the extension $K \hookrightarrow L$. The geometric resolution
of $V$ associated to  $\ell:= (\ell_0, x_1, \dots, x_r)$ is given by
\begin{eqnarray*}
p & := &  \cP_V (c_0, e_1, \dots, e_r )(  t, x_1, \dots, x_r)
\ \ \in \Q[x_1, \dots, x_r][t] , \\[2mm]
w& :=&  - \left(
 \frac{\partial \cP_V }{ \partial U_{0 1}} , \dots,
 \frac{\partial \cP_V }{ \partial U_{0  n}}  \right)
(c_0, e_1, \dots, e_r)(  t, x_1, \dots, x_r ) \ \ \in \Q[x_1,
\dots, x_r][t]^n,
\end{eqnarray*}
where $\cP_V$ is the normalized characteristic polynomial of $V$.

\end{lem}

\begin{proof}{Proof.--}

The fact that $\pi$
 is finite follows from \cite[Lemma 2.14]{KrPaSo99}.
On the other hand, the normalization imposed on
 $\cP_V$ implies that
$$
p_0(t):= p(t, 0, \dots, 0) =
\cP_V (c_0, e_1, \dots, e_r )( t, 0, \dots,
0)  \ \in \Q[t]
$$
is a monic ---and thus non-zero--- polynomial of degree $D$ which
vanishes on $ \ell_0 (\xi) $ for all $\xi \in Z$.
The hypothesis that $\ell_0$ separates the points of $Z$ implies
that
$p_0$ is the minimal polynomial of $\ell_0$ over $Z$; in particular it is
 a {\em squarefree}  polynomial of degree $D$ and so, as $p$ is monic,
$$
0 \ne \discr \, p_0 = (  \discr_{t} \, p )(0, \dots,0).
$$

In particular, $ \discr_{t} \, p \ne 0$ and thus $p$ is also a  squarefree polynomial
which annihilates
$\ell_0$ over $V$.

Now, as the map $\pi$ is finite, the minimal
polynomial $m_{\ell_0}  \in K[t]$ of $\ell_0$  lies in $\Q[x_1,
\dots, x_r][t]$.
Hence $ m_{\ell_0}(0, \dots, 0, t) $
vanishes on $ \ell_0 (\xi) $ for all $\xi \in Z$.
This implies that $\deg_t m_{\ell_0} =  D $. As  $p$ is a monic polynomial of degree
 $D$ in $t$, then
  $p = m_{\ell_0} $. So
   $\ell_0$ is a primitive element of the  extension
$K \hookrightarrow L$, and $p$ is its minimal polynomial.

\smallskip

Using the same notation  of Lemma \ref{resolucion generica} we have
$$
 \frac{\partial P_V}{\partial T_0}(U,L(U, x )) \, x_i = -
\frac{\partial P_V}{\partial U_{0 \, i}}(U,L(U,x)) \ \ \ \in \cL.
$$
As this identity only involves  polynomials in $\Q[U_0, \dots, U_r][x_1, \dots,x_n]$,
it can be directly evaluated to
obtain the parametrization $w$.
\end{proof}

In particular this shows
that the {\em total } degree of the polynomials in the geometric
resolution is bounded by $\deg p \le D$ and  $\deg w_i \le D$
(See also \cite[Prop. 3]{GiLeSa01}).

\smallskip

Lemma \ref{especializacion} can  be directly applied to compute a geometric
resolution of an equidimensional variety $V$ which satisfies Assumption
\ref{assumption} from a given Chow form of $V$:

\begin{cor} \label{GeomResFromChow}
Let notations and assumptions  be as in Lemma~\ref{especializacion}.
Suppose that there is  given a Chow form $\cF_V$  of $V$, encoded by a
 slp of length
$L$. Then, there is an algorithm which computes a geometric
resolution of  $V$ associated to $\ell$ within complexity $\cO( n
\, L) $. All polynomials arising in this geometric resolution are
encoded by
  slp's of length $\cO(L)$. \hfill $\square$
\end{cor}

Lemma \ref{especializacion} also  yields, from $\Ch_V$,  a
geometric resolution of the fiber $Z$ associated to an  affine
linear form $\ell_0$, as $\Ch_Z(U_0) = \Ch_V(U_0, e_1, \dots,
e_r)$. This is summarized in Procedure $\GeomRes$ (Subroutine
\ref{algo-geomres}).

\begin{algo}
{Computing a geometric resolution of a fiber} {algo-geomres} {$
\GeomRes (n,r , D, \Ch_V, \xi, c)$ }

\hspace*{4.1mm}{\# $n$ is the number of variables,

\hspace*{4.1mm}\# $r, D$ are the dimension and an upper bound for
the degree of $V$ respectively,

\hspace*{4.1mm}\# $\Ch_V$ is the normalized  Chow form of $V$,

\hspace*{4.1mm}\# $\xi:= (\xi_1,\dots, \xi_r)\in \A^r$ such that
$\# Z_\xi  = \deg V$, where $Z_\xi:=V \cap V(x_1 - \xi_1,\dots, x_r -
\xi_r)$.

\hspace*{4.1mm}\# $c_0\in \Q^{n+1}$ s.t. $\ell_0:=L_0(c_0,x)$
is the considered  affine linear form.

\smallskip

\hspace*{4.1mm}\# The procedure returns $(D_0, p, v)$, where $D_0$
is the degree of $V$ and $(p, v)\in  \Q[t]^{n+1}$ is the

\hspace*{4.1mm}\#   geometric resolution of $Z_\xi$ associated to $\ell_0$ in case $\ell_0$
separates the points in $Z_\xi$.

\hspace*{4.1mm}\#
Otherwise, it returns error.

\begin{enumerate}

\item $P(U_0,t):=\Ch_V ((U_{0 0} -t, U_{0 1},\dots, U_{0 n}),e_1 - \xi_1
e_0, \dots, e_r - \xi_r e_0)$;\vspace{-1,5mm}

\item $(p_0,\dots, p_D):= \Expand(P(c_0, t), t, 0, D)$;
\vspace{-1,5mm}

\item $D_0:= D$;\vspace{-1,5mm}

\item {\bf while} $p_{D_0} =0$ and $D_0\ge 0$ {\bf do} \vspace{-1,5mm}

\item \hspace*{4.1mm} $D_0:= D_0 - 1$; \vspace{-1,5mm}

\item {\bf od};\vspace{-1,5mm}


\item  $p:= (-1)^{D_0} P(c_0,t)$;\vspace{-1,5mm}

\item $(\rho, q_1, q_2) := {\rm Res}(p, p', D_0, D_0-1)$\label{resultante1};\vspace{-1,5mm}

\item {\bf if } $\rho = 0 $ {\bf then}\vspace{-1,5mm}

\item \hspace*{4.1mm} {\bf return }(``error'');\vspace{-1,5mm}

\item {\bf else}\vspace{-1,5mm}

\item \hspace*{4.1mm} $(w_1,\dots, w_n):= ((-1)^{D_0+1}
\partial P/\partial U_{01}(c_0, t), \dots, (-1)^{D_0+1}
\partial P/\partial U_{0n}(c_0, t))$;\vspace{-1,5mm}

\item \hspace*{4.1mm} $(v_1,\dots, v_n):=
 ( {\rm Mod}(\frac1\rho \ q_2 \,w_1,\, p\,, 2D_0-1, D_0), \dots, {\rm Mod}(\frac1\rho \ q_2 \,w_n,\, p\,, 2D_0-1, D_0 ))$ \label{mod};\vspace{-1,5mm}

\end{enumerate}

\hspace*{7mm} {\bf return}$(D_0, p, v_1, \dots, v_n)$;}

\end{algo}

In Procedure $\GeomRes$ (Subroutine \ref{algo-geomres}), as we do
in all zero-dimensional situations, we use the definition of
geometric resolution stated in Section \ref{The representation of
the Chow form} to avoid divisions by $p'$. In line
\ref{resultante1} of this subroutine, ${\rm Res} (f, g, d_1,d_2) $
is a procedure that, using basic linear algebra, computes $(\rho,
q_1, q_2)$ where $\rho$ is the resultant between the univariate
polynomials $f$ and $g$ of degrees $d_1$ and $d_2$ respectively,
and $q_1$ and $q_2$ are polynomials of degrees bounded by $d_2 -1$
and $d_1 -1 $ respectively satisfying $\rho = q_1 f + q_2 g$. In
line \ref{mod}, ${\rm Mod} (f, g, d_1,d_2)$ is a procedure that
computes the remainder of the division of the polynomial $f$ of
degree bounded by  $d_1$  by  the polynomial $g$ of degree bounded
by $d_2$.

\smallskip
\begin{prop} \label{GeomResFiber}
Let $V\subset \A^n$ be an equidimensional variety of dimension $r$
and degree bounded by  $D$. Let $(\xi_1, \dots, \xi_r)\in \A^r$
be such that $Z_\xi:=V \cap V(x_1 - \xi _1, \dots, x_r- \xi_r)$
is a 0-dimensional variety of cardinality $\deg V$. Let be given
$\Ch_V$ encoded by a slp of length $L$ and the coefficients of an
affine linear form $\ell_0\in \Q[x_1,\dots,x_n]$ which separates
the points in $Z_\xi$. Then, Procedure $\GeomRes$  (Subroutine
\ref{algo-geomres}) computes a geometric resolution of  $Z_\xi$
(in the sense of Section \ref{The representation of the Chow
form}) within complexity $\cO(nD^2L+D^4)$. \hfill $\square$
\end{prop}

\medskip

On the other hand, next result shows the converse of Corollary
\ref{GeomResFromChow}: To derive a Chow form from a given
geometric resolution is quite standard in the zero-dimensional
case,   but  it was  by no means clear up to now how to generalize
that for varieties of arbitrary dimension.  Here we show how to do
that  within polynomial complexity. This is done by deriving, from
the given geometric resolution of $V$, a geometric resolution of
the fiber $Z$ and  a system  of local equations for $V$ at $Z$,
which enables  us to apply Procedure $\ChowForm$ (Subroutine
\ref{algo-fibra}).

\begin{prop}\label{3.5}
Let $V \subset \A^n$ be an equidimensional variety of dimension
$r$ and degree $D$ which satisfies Assumption~\ref{assumption}.
Let
$\ell_0 \in \Q[x_1, \dots, x_n]$ be a
linear form which separates the points of
$Z$.

Suppose that there is given a  geometric resolution
 $(p, w)$ of $V$  associated to $\ell:= (\ell_0, x_1, \dots,
 x_r)$,
encoded by slp's of length $L$.
Then there is a bounded probability algorithm which computes (a slp
for)
$\Ch_V$  within complexity
$\cO(n^{16}D^{19} (D+L)) $.
\end{prop}

\begin{proof}{Proof.--}

First we derive a geometric resolution of $Z$ associated to
$\ell_0$:

We know that  $ \Ch_Z(U_0) = \Ch_V(U_0,
e_1, \dots, e_r)$. Thus,
$$
\cP_Z (U_0)(t) = \cP_V (U_0, e_1, \dots, e_r)(t, 0, \dots, 0) \
\in \Q[U_0][t].
$$
The geometric resolution $(p,w)$ of $V$ associated to $\ell$ is given by
Lemma~\ref{especializacion}. Applying the
same lemma  to $Z$, we deduce
that the geometric resolution $(p_0,w_0)$ of
$Z$ associated to $\ell_0$   is $p_0(t):=p(t, 0, \dots, 0) \in \Q[t]$ and
$ w_0(t):= w(t, 0,
\dots, 0) \in \Q[t]^n$.

\smallskip
Now,   let us derive a system of local equations of $V$ at $Z$:

Let $c_i \in \Q^{n+1}$, $r+1\le i \le n $,
be such that the  affine linear forms
 $
\ell_i := c_{i  0} + c_{i 1}x_1 + \cdots + c_{i n}x_n \  \in \Q[x_1, \dots, x_n]
$
are linearly independent and such that each of them
 separates the points of $Z$.

 For $r+1\le i \le n$ define

$$H_i:= |\,p'(M_p)\, t -  (p'\, \ell_i(w/p')) (M_p)\,|
$$
where
$ M_p \in \Q[x_1, \dots, x_r]^{D \times
  D}$ denotes the companion matrix of $p$.
  Since $p' \, \ell_i (w/p')$ belongs to $ \Q[x_1, \dots, x_r][t]$,
  $H_i  \in \Q[x_1, \dots, x_r][t]$.

  Observe that $x_i = (w_i/p')(\ell_0(x))$ in
  $L$ implies that in $L$
$$p'(\ell_0)\ell_i = c_{i 0} p'(\ell_0) + c_{i 1}w_1(\ell_0) +
\cdots + c_{i n} w_n (\ell_0) = ( p'\, \ell_i ( w/p'))(\ell_0).$$
Thus, as $M_p$ is the matrix of multiplication by
$\ell_0$ with respect to $K \hookrightarrow  L$, we conclude that
$H_i= |\,p'(M_p)\,|\,m_{\ell_i}$ where $m_{\ell_i}$ is the minimal
polynomial of $\ell_i$
over $K$.

The assumptions that
$\ell_i$ separates the points of $Z$ and
 the projection $ \pi: V \rightarrow \A^r  ,
x \mapsto (x_1, \dots, x_r)$ is finite, imply that $m_{\ell_i}$
belongs to   $\Q[x_1, \dots,
x_r][t] $.

 Therefore, for $r+1\le i \le n$, we can define
 $$f_i:= m_{\ell_i}(\ell_i) = \frac {1}{| p'(M_p) | }H_{i}(x_1, \dots, x_r)(\ell_i)
$$
These are
squarefree
 polynomials in separated  variables which vanish over $V$, and
 so, it is easy to verify from the Jacobian criterion that
$ f_{r+1}, \dots, f_{n}$ is a system  of reduced local equations of $V$ at
$Z$.

Observe that as
there exist $a,b \in \Q[x_1,\dots,x_r][t] $ such that
$\discr (p) = a(t)\, p (t)+ b(t) \, p'(t)$, $\discr(p)\,\Id = b(M_p)\,p'(M_p)
$. On the other hand, as  $\deg Z = \deg V$, $ \discr (p) (0,\dots,0)= \discr (p_0) \ne 0$.
Therefore,    $|\, b(M_p)\,p'(M_p)\,| (0,\dots,0)=(\discr (p_0))^D \ne
0$. We conclude that $|\,p'(M_p)\,|(0,\dots,0)\ne 0$ and hence, we can
use the point $(0,\dots,0)$  to perform Procedure $\PolynomialDivision$
(Subroutine \ref{PolynomialDivision}) in order to obtain division free slp's for
$f_{r+1},\dots,f_n$.

\smallskip
Finally, we apply
procedure $\ChowForm $ (Subroutine \ref{algo-fibra}) to $Z$  and $\{ f_{r+1},\dots,f_n\}$.

\smallskip

Let us decide now the random choices in order to
insure that the algorithm has an error probability
bounded by $1/4$:

We need  $c_{r+1},\dots, c_n \in \Q^{n+1}$ satisfying the stated
conditions of independence and separability. These conditions are
satisfied provided that $$ \rho (c_{r+1}) \cdots \rho (c_{n}) \, \left|
(c_{i j} )_{\substack{r+1 \le i \le n \\ 1 \le j \le n-r}}\right| \ne 0,
$$ where   $\rho:= \discr_{t} \,  \cP_Z \in \Q[U_0] \setminus \{
0\}$. As $\cP_Z$ is an homogeneous polynomial of degree $D$ and
$\deg_t \cP_Z = D$, $\deg \rho \le D\,(2D-1)$. Thus the degree of
the polynomial giving bad choices is bounded by $(n-r)\,D\,(2D-1)
+ (n-r)$. We choose $\ell:= 8\,n\,D^2$ in order to apply Schwartz
lemma.
\smallskip

Now we compute the complexity of the algorithm:

\smallskip

The dense  representation of the geometric resolution of $Z$
associated to $\ell_0$
is computed within complexity
$\cO(n\, D^2 \, L)$ (using Procedure $\Expand$).

The construction of the random choice for the  affine linear forms  $\ell_{r+1}, \dots, \ell_n$
is not relevant here.

The computation of each polynomial $H_i$ requires $\cO(D^4)$ operations
for the computation of the
determinant plus the computation of each coefficient of the matrix,
that is $\cO(D^3 L)$ more operations, Hence, computing $H_i$ requires
 $\cO( D^3\,(D+L))$ operations.

By  Lemma~\ref{vermeidung}, taking into account that the total
degree of each $f_i$ is bounded by $D$ (since it is the minimal
polynomial of the  affine linear form $\ell_i$), and that the
lengths of $H_i$ and $| p'(M_p) |$ are of order $\cO(D^3 (D +
L))$,
 the complexity of the final division
for computing each
$f_i$ is $\cO(D^2(D+ D^3 (D+L)) )= \cO(D^5 (D + L))$.

Finally, Lemma~\ref{fibra} gives the final complexity
$\cO(r^8\, \log_2(r\,D) n^7 D^{13} D^5 (D + L))= \cO(n^{16}D^{19} (D+L))$.
\end{proof}


\typeout{Intersection of a variety with a hypersurface}

\subsection{Intersection of a variety with a hypersurface}

\label{Intersection of a variety with a hypersurface}

\vspace{2mm}

Let  $V\subseteq \A^n$ be  an equidimensional variety defined over $\Q$, and let $f
\in \Q[x_1, \dots, x_n] $ be a non-zero divisor modulo $I(V)$.
In this section, we compute, from the Chow form of $V$ and the equation $f$,
 the Chow form of the
set-theoretic intersection $V \cap V(f) \subset \A^n$. In order to do this,
we use {\em generalized Chow forms}, which we define
now.
We refer to \cite{Philippon86} and \cite[Section 2.1.1]{KrPaSo99}
for a more extensive treatement of these generalized Chow forms.

\medskip
We assume that $\dim V = r$ and that $\deg f \le d$.
As before, for $i=0, \dots, r$, we introduce a group $U_i=(U_{i 0},\dots,U_{i n})$
of $n+1$ variables;  we
introduce also a  group    $U(d)_{ r}$  of \ $d+n \choose {n}$ \
variables. We set
$$
L_i:=U_{i0}+U_{i1}x_1+ \cdots + U_{in}x_n \quad \quad , \quad
\quad F_r:=\sum_{|\alpha|\le d} U(d)_{r \alpha} \, x^\alpha
$$
for the generic affine linear forms in $n$ variables associated to $U_i$
and the generic polynomial of degree $d$ in $n$ variables
associated to $U(d)_r$.

Set $N:= r\, (n+1) +{{d+n}\choose{n}} $ and let $W\subset \A^N
\times V$ be the incidence variety of $L_0, \dots, L_{r-1}, F_r$
with respect to  $V$, that is
$$ \displaylines{\qquad
W:= \{ (u_0, \dots, u_{r-1}, u(d)_r; \xi) \in \A^N \times \A^n ; \
\hfill\cr \hfill \xi\in V,  \ L_0(u_0,\xi) =0, \dots,
L_{r-1}(u_{r-1},\xi)=0,\, F_r(u(d)_r,\xi) =0 \}.\qquad}$$

Let  $\pi : \A^N \times \A^n \to \A^N$ denote the canonical
projection   onto the first coordinates. Then $\overline{\pi(W)}$
is a hypersurface in $\A^N$. A  {\em generalized Chow form} or
{\em $d$-Chow form } of $V$ is
any squarefree polynomial $\cF_{d, V}\in \Q[U_0, \dots, U_{r-1},
U(d)_r]$ defining $\overline{\pi (W) } \subseteq \A^N $.

A $d$-Chow form \ $\cF_{d,V} $
\ happens to be a
  multihomogeneous polynomial
  of
degree $d\, \deg V$ in each group of variables $U_i$, and
 of degree $\deg V$ in the group $U(d)_r$.

\smallskip

If the variety $V$ satisfies Assumption \ref{assumption}, we
define the {\em normalized} $d$-Chow form of $V$
as the
unique $d$-Chow form $\Ch_{d,V}\in \Q[U_0,\dots, U_{r-1},U(d)_r]$ of $V$ satisfying $\Ch _{d,V} (e_0, \dots,
e_{r-1}, e(d)) = 1$, where $e(d)$ is the vector of coefficients of
the polynomial $x_r^d$.

\bigskip

Let $\overline{V}$ and $ \overline{V(f) }$ denote the closure in
$\P^n$
 of $V$ and
$V(f)$ respectively. Set \ $\overline V \cap \overline{V( f)} =
\bigcup_{C} C$ \ for the
 irreducible
decomposition of $\overline V \cap \overline{V( f)} \subset \P^n$
and, for each irreducible component $C$, let $\cF_C \in \Q[U_0,
\dots, U_{r-1}]$ denote a Chow form of $C$. Then \cite[Prop.
2.4]{Philippon86} states that
\begin{equation}\label{Chowdeve}
\Ch_{d, V}(U_0, \dots, U_{r-1} , f)  = \lambda \,\prod_{C}
\cF_{C}^{m_C}
\end{equation}
for some $\lambda \in \Q^*$ and some positive integers $m_C \in \N$.
(Here we wrote $\Ch_{d, V}(U_0, \dots, U_{r-1} , f)$  for the
specialization of the  group $U(d)_r$  into
the coefficients of the polynomial $f$.)

On the other hand, as $V \cap V(f) = \bigcup_{C\not\subset \{ x_0
= 0\}} C$, the polynomial
$$
\prod_{C\not\subset \{x_0=0\}} \cF_{C}
$$
is a Chow form of $V \cap V(f)$.

Hence, in order to compute $\Ch_{V \cap V(f)}$, the goal is to
compute first
$\Ch_{d,V}(f):= \Ch_{d, V}(U_0, \dots, U_{r-1} , f)$ and then
clean the multiplicities and the Chow forms of components
contained in the hyperplane $\{x_0=0\}$.

\bigskip

The following result enables us to compute a $d$-Chow form
from the standard one. We recall some of the notation of
Subsection~\ref{A product formula}: for an
equidimensional variety $V \subset \A^n$ of dimension $r$ and degree
$D$
satisfying Assumption \ref{assumption}, we set $K:= \Q(U_0, \dots,
U_{r-1})$ and  $I(V)^\e$ for  the extension of the ideal of
$V$ to  $K[x_1, \dots, x_n]$. Also recall that
$$
V^0 := V(I(V)^\e ) \cap V(L_0, \dots, L_{r-1})  \
\subset \A^n(\overline{K}),
$$
 is a 0-dimensional variety of degree $D$, and that $\No_{V^{0}}$ refers to the Norm
  as defined
in Section \ref{Newton's algorithm}.

\begin{lem} \label{Ch-d}
Under Assumption \ref{assumption}, we have
$$
\Ch_{d, V} = \Ch_{V}(U_0, \dots, U_{r-1} , e_0)^d \ \No_{V^0}(F_r).
$$

\end{lem}

\begin{proof}{Proof.--}

Let $\Ch_{d,  V^0} \in K[U_r]$ be the $d$-Chow form of $V^0$.
First one shows ---exactly as in Lemma~\ref{A product formula}---
that there exists \
 $\lambda_d \in K^*$ such that \ $ \Ch_{d,
 V}
= \lambda_d \, \Ch_{d, V^0} $. Set  $e(d)_0 $ for the vector of
coefficients of the polynomial $x_0^d$. Evaluating this identity
at $U(d)_0 \mapsto e(d)_0$ we obtain
$$
 \Ch_{d, V}(U_0, \dots, U_{r-1}, e(d)_0) = \lambda_d
\, \Ch_{d,V^0} (e(d)_0) = \lambda_d.
$$

Consider the morphism $ \varrho_{d}: \Q[U_0, \dots, U_{r-1},
U(d)_r]  \to  \Q[U_0,\dots, U_{r-1}, U_r] $ defined by $
\varrho_d(L_i) = L_i$ for $0\le i \le r-1$ and $ \varrho_d(F_r)
= L_r^d$. Then $\varrho_d(\Ch_{d, V}) = \Ch_V^d $ (see
\cite[Lemma 2.1]{KrPaSo99}), which implies that
$$
\Ch_{d, V}(U_0, \dots, U_{r-1}, e(d)_0) = \varrho_d( \Ch_{d, V} )
(U_0, \dots, U_{r-1} , e_0) = \Ch_{ V}(U_0, \dots, U_{r-1},
e_0)^d.
$$
Therefore $\lambda_d =  \Ch_{ V} (U_0, \dots, U_{r-1} , e_0)^d  $.
The statement follows immediately from this identity and the
observation that $ \Ch_{d, V^0} = \No_{V^0}(F_r)$.
\end{proof}

To clean the components of $\overline V \cap \overline{V(f)}$
lying in the  hyperplane $\{ x_0 = 0 \}$ we use the
following criterion:

\begin{lem}
Let $W \subset \P^n$ be an irreducible variety
of dimension $r-1$.
Then  $W \subset  \{x_0 =0\}$ if and only if
 $\cF_{W}$   does not
depend on the variable $U_{00}$.
\end{lem}

\begin{proof}{Proof.--}

In case $\cF_W$ does not depend on $U_{00}$ we have that
$$
\cF_W (e_0, U_1, \dots, U_{r-1}) = 0,
$$
which is equivalent to the fact that $W$ is contained in the
hyperplane $\{x_0 =0 \}$.

\smallskip

On the other hand, assume that $W \subset \{ x_0 =0\} \cong
\P^{n-1}$. Then $\cF_W$ coincides with the Chow form of $W$
considered as a subvariety of this linear space, see e.g. the
proof of \cite[Lemma 2.6]{KrPaSo99}. Hence $\cF_W$ does not depend
on $U_{00}$, and as a matter of fact, it does not depend on any
of the variables $U_{i\, 0} $ for $0\le i  \le  r-1$.
\end{proof}

Let again  $\cF_C \in \Q[U_0,
\dots, U_{r-1}]$ denote a Chow form of  an irreducible component $C$ of
$\overline V \cap \overline{V( f)}  \subset
\P^n$. Recalling Identity \ref{Chowdeve}, set
$$
\cF_1 := \prod_{C \subset \{x_0 =0\} } \cF_{C}^{m_C}
 \quad \quad
\mbox{and} \quad \quad
\cF_2 := \prod_{C \not\subset \{x_0 =0\} } \cF_{C}^{m_C}.
$$

Then $ \Ch_{d, V}(f) = \lambda\, \cF_1 \cF_2$ for $\lambda \in
\Q^*$ and the squarefree part $(\cF_2)_\red$ of $\cF_2$ is a Chow
form of $V \cap V(f)$.

By the previous lemma, $\cF_1$ does not depend on
$U_{00}$, while all the factors of $\cF_2$ do.

Therefore
$$
\frac{\partial \Ch_{d, V}(f)}{\partial U_{0 0}}
 =  \lambda\, \cF_1 \frac{\partial \cF_2}{\partial U_{0 0}}
$$
and so
\begin{equation}
\label{chauinfinito}
 \cF_{V\cap V(f) } := \frac{\Ch_{d, V}(f)}{{\rm gcd}(
\Ch_{d, V}(f), {\partial \Ch_{d, V}(f)}/ {\partial U_{00}})}
\end{equation}
is a Chow form of $V \cap V(f)$.

\bigskip

\begin{lem}
\label{intersection} Let $V\subset \A^n$ be an equidimensional variety of degree $D$
which satisfies Assumption \ref{assumption} and let $f\in \Q[x_1,\dots,x_n]$
of degree bounded by $d$ be a non-zero divisor modulo $I(V)$.
Assume that $\Ch_V$ and
$f$ are encoded by  slp's of length bounded by $L$.

Then  there  is a bounded probability algorithm (Procedure
$\Inter $ (Subroutine~\ref{algo-inter}) below) which computes the
Chow form $\Ch_{V \cap V(f)}$ of the intersection variety $V \cap
V(f)$ within (worst-case) complexity \ $\cO((ndD)^{12}L)$ .
\end{lem}

\begin{proof}{Proof.--}
Our first goal is to compute  $\Ch_{d, V}(f)\in \Q[U_0,\dots,U_{r-1}]$
 by means of  Lemma \ref{Ch-d}.
To obtain $\No_{V^0}(f)$ we derive first
 a geometric resolution of
$V^0$ from its characteristic
polynomial   and  Lemma~\ref{resolucion generica}. It is easy to
check that the polynomial
$$
p(t):= (-1)^D \Ch_V(U_0, \dots, U_{r-1}, (U_{r0}-t, U_{r1}, \dots,
U_{rn}))$$ is a characteristic polynomial of $V^0$, with leading
coefficient  $a := \Ch_V(U_0, \dots, U_{r-1}, e_0)$.

Then, the geometric resolution of $V^0$ associated to $L_r$ is given by
$$\frac{1}{a}\,  p(t) \in K[U_r][t] \qquad {\rm and} \qquad
\frac{1}{a} \, w(t)  \ \in K[U_r][t]^n \quad {\rm where} \quad w:=-\left(
 \frac{\partial p}{\partial U_{r 1}} , \dots,
\frac{\partial p}{\partial U_{r n}} \right)   .$$


For $\gamma\in V^0$, if we denote by $f^h$ the
homogeneization up to degree $d$ of $f$ with respect to a new variable
$x_0$ and $p'$ the derivative of $p$ with respect to $t$,
we have
$$p'(L_r(\gamma))^d\,f(\gamma)=
f^h(p'(L_r(\gamma)),w_1(L_r(\gamma)),\dots,w_n(L_r(\gamma))).$$
Thus, if $M$ denotes the companion matrix of $\frac{1}{a}\,  p(t)
$, we get
$$|p'(M)|^d\, \No_{V^0}(f) = | f^h(p'(M),w_1(M),\dots,w_n(M))|.$$
In order to avoid  divisions  (since $M\in K[U_r]$),
we replace $M$
by $M_p:=a\,M $ and $p'$, $w_1,\dots,w_n$ by their
homogeneizations $(p')^h$, $w^h_1,\dots,w^h_n$ up to degree $D$
such that $M_0:=a^D p'(M) = (p')^h(a\,\Id,M_p)$ and for $1\le i
\le n $, $M_i:=a^D w_i(M) = w^h_i(a\,\Id, M_p)$. Therefore,
multiplying both sides by $a^{dD^2}=|a^D\,\Id|^d$, we obtain
$$|M_0|^d \,\No_{V^0}(f) = | f^h(M_0,M_1,\dots,M_n)|.$$

Finally, from  Lemma \ref{Ch-d}, we conclude that \begin{equation}\label{ChdVdeF}
\Ch_{d,V}(f)= a^d
\No_{V^0}(f) = \frac{a^d\, | f^h(M_0,M_1,\dots,M_n)|}{|M_0|^d}
\ \in \Q[U_0,\dots,U_{r-1}].
\end{equation}

We  compute this quotient applying Procedure $\PolynomialDivision$ (Subroutine
 \ref{PolynomialDivision}).

Next we apply Identity \ref{chauinfinito} to
compute a Chow form $\cF:= \cF _{V\cap V(f)}$ from $\Ch_{d,V}(f)$:
we first compute the polynomial
\begin{equation}\label{Ge}
G:= {\rm gcd} (\Ch_{d,V}(f), \partial \Ch_{d,
V}(f)/\partial U_{00})\end{equation}
applying Procedure $\GCD $ (Subroutine
\ref{GreatestCommonDivisor})
and then perform the division $\cF = \Ch_{d,V}(f) / G$ applying again
Procedure $\PolynomialDivision $.

Finally, as Assumption \ref{assumption} holds,
$\cF(e_0,\dots,e_r)\ne 0$ and we obtain  the normalized Chow form
$\Ch_{V\cap V(f)}= \cF /\cF(e_0,\dots,e_r)$.

\medskip

Now let us check the size of the sets of points we have to take to
insure that the algorithm has an error probability bounded by
$1/4$:

First, in order to compute $\Ch_{d,V}(f)$
 we need $u\in \Q^{(r+1)(n+1)}$ such that $| M_0 | (u) \ne
0$. But let us observe that in fact $| M_0 |(e_1,\dots,e_r,U_r)\ne 0 \in \Q[U_r]$
(so it is enough to choose randomly $u_r \in \Q^{(n+1)}$ such that
$| M_0 |(e_1,\dots,e_r,u_r)\ne 0$). This is due to the fact
 that $a(e_1,\dots,e_{r})=\Ch_V(e_1,\dots,e_{r},e_0)=\pm 1$, thus Assumption
 \ref{assumption}
implies
that $\Ch_V(e_1,\dots,e_r,U_r)=\pm \Ch_Z(U_r)$. Hence,
 $p_Z(U_r,t):=p(e_1,\dots,e_r)(U_r,t)$
 is a characteristic polynomial of $Z$, whose
 discriminant does not vanish, and then, the
 polynomial
$| M_0|(e_1,\dots,e_r,U_r)\ne 0 \in \Q[U_r]$.

Now, as $\deg | M_0|(e_1,\dots,e_r,U_r) \le D^2$, if we take $u_r:=
\Random (n+1, 12\,D^2)$ we infer that with probability at least $
1 - {1/12}$, $(e_1,\dots,e_r,u_r)$ is a good base point to apply
Procedure $\PolynomialDivision $ and obtain $\Ch_{d,V}(f)$.

\smallskip
Next   we compute $G$ applying Procedure $\GCD$
$\lceil{6\,(1+\log 12)}\rceil=26$ times (see Remark
\ref{remarkrepeticion} so that its error probabiliy is at most
${1/12}$).

\smallskip
Finally as $G$ is a polynomial of degree bounded by $r\,d\, D$ in
$r\,(n+1)$ variables, choosing $u:= \Random(r\,(n+1),12\,r\,d\,
D)$ we also guarantee that the probability that $u$ is a good base point to perfom
the last division is at least $ 1 -  {1/ 12}$.

\smallskip
Thus, the error probability of the whole algorithm is at most
${1}/{4}.$

\medskip
Now let us compute the (worst-case) complexity of this algorithm:

The whole complexity of computing the numerator and denominator in
Identity \ref{ChdVdeF} is of order $\cO(n (d^2 D^3 L +
D^4))=\cO(nd^2D^4L)$. By Lemma \ref{vermeidung} the complexity of
computing $\Ch_{d,V}(f)$ is of order $\cO((n d D)^{2} (ndD +
nd^2D^4)L)=\cO(n^3d^4D^6L)$.

Then, we apply Lemma \ref{mcd} and Proposition \ref{Probabilidad de Ms}
 to compute a slp of length $\cO(n^7d^8D^{10}L)$ for $G$ of Identity  \ref{Ge}
within complexity $ \cO(n^9d^8D^{10} L)$.

Finally, when we perform the last division, the overall complexity
of computing $\Ch_{V\cap V(f)}$ is of order $ \cO((ndD)^{12} L)$.
\end{proof}

We summarize the algorithm in Procedure $\Inter$ (Subroutine
\ref{algo-inter}).

\begin{algo}
{Intersection with a hypersurface} {algo-inter} {$ \Inter(n,r, D
,f, d,  \, \Ch_V)$ }

\hspace*{4.1mm}{\# $n$ is the number of variables,

\hspace*{4.1mm}\# $r, D$ are the dimension and the degree of $V$
respectively,

\hspace*{4.1mm}\# $ f \in \Q[x_1, \dots, x_n] $ is a non-zero
divisor modulo $I(V)$ of degree bounded by $d$,

\hspace*{4.1mm}\# $\Ch_V$ is the normalized  Chow form of $V$,

\smallskip

\hspace*{4.1mm}\# The procedure returns the  normalized Chow form
$\Ch_{V \cap V(f)}$ of the intersection variety $V \cap V(f)$.

\begin{enumerate}

\item $p:= (-1)^D \Ch_V(U_0, \dots, U_{r-1}, (U_{r0}-t, U_{r1},
\dots, U_{rn}))$;\vspace{-1,5mm}

\item $a := \Ch_V(U_0, \dots, U_{r-1}, e_0)$;\vspace{-1,5mm}


\item  $\displaystyle{
w :=  - \left( \frac{\partial p}{\partial U_{r 1}} , \dots,
\frac{\partial p}{\partial U_{r n}} \right)}$ ;\vspace{-2mm}


\item $M_p:= a\,  {\rm CompanionMatrix}(p/a)$;
\label{principio1}\vspace{-1,5mm}

\item {\bf for} $i$ {\bf from} 1 {\bf to} $n$ {\bf do}\vspace{-1,5mm}

\item \hspace*{4.1mm} $w_i^h := \Homog(w_i, D)$;\vspace{-1,5mm}

\item \hspace*{4.1mm} $ M_i := w_i^h (a, M_p)$;\vspace{-1,5mm}

\item {\bf od};\vspace{-1,5mm}

\item $(p')^h := \Homog(\frac{\partial p}{\partial t}, D)$;\vspace{-1,5mm}

\item $M_0 := (p')^h (a, M_p)$;\vspace{-1,5mm}

\item $f^h := {\rm Homog}(f, d)$;\vspace{-1,5mm}

\item $M_f:= f^h (M_0, M_1, \dots, M_n)$;\vspace{-1,5mm}

\item $H_1:= |\,M_f \, |$; \vspace{-1,5mm}

\item $H_2:= |\, M_0 \, |$; \label{final1}\vspace{-1,5mm}

\item $u_r:= {\rm Random}(n+1, 12\,D^2)$;\vspace{-1,5mm}

\item {\bf if} $ H_2(e_1, \dots, e_{r},u_r) = 0$ {\bf then}\vspace{-1,5mm}

\item \quad {\bf return}(``error'');\vspace{-1,5mm}

\item {\bf else}\vspace{-1,5mm}

\item \hspace*{4.1mm}  $\Ch_{d, V}(f):= \PolynomialDivision(a^d H_1, H_2^{d}, r d D,
(e_1,\dots, e_{r},u_r))$; \label{principio2}\vspace{-1,5mm}

\item \hspace*{4.1mm} $G := \GCD(\Ch_{d, V}(f), \partial \Ch_{d,
V}(f)/\partial U_{00}, (U_0,\dots, U_{r-1}), r d D; 12)$;
\label{principio3}\vspace{-1,5mm}

\item  \hspace*{4.1mm} $u:= {\rm Random}((r+1)(n+1), 12 r d D)$;
\label{elegiru}\vspace{-1,5mm}

\item  \hspace*{4.1mm} {\bf if} $ G(u) = 0$ {\bf then}\vspace{-1,5mm}

\item  \hspace*{8.2mm} {\bf return}(``error'');\vspace{-1,5mm}

\item  \hspace*{4.1mm} {\bf else}\vspace{-1,5mm}

\item  \hspace*{8.2mm} $\cF := {\rm PolynomialDivision}(G, \Ch_{d, V}(f),
r D d, u)$;\vspace{-1,5mm}

\item  \hspace*{8.2mm} $\Ch_{V \cap V(f)}:= \cF /\cF(e_0, \dots,
e_r)$;\label{final3}\vspace{-1,5mm}

\end{enumerate}

\hspace*{7.7mm} {\bf return}($\Ch_{V \cap V(f)}$); \vspace{-1,5mm}}
\end{algo}


\typeout{Separation of varieties}

\subsection{Separation of varieties}

\label{Separation of varieties}

\vspace{2mm}

Let $V \subset \A^n$ be an equidimensional variety of dimension
$r$. Let $g \in \Q[x_1, \dots, x_n] \setminus\{0\}$, and set $Y$
for the union of the irreducible components of $V$  contained in $V(g)$ and
$W$ for the union of the other ones. Hence $Y$ and $W$ are
equidimensional varieties of dimension $r$ such that
$V = Y \cup W$, $ Y \subset V(g)$ and $ g $ is not
a zero divisor modulo $ I(W)$.

The following procedure (Subroutine~\ref{algo-sepchow}) computes the
Chow forms of $Y$ and $W$ from a Chow form of $V$ and the
polynomial $g$.

\medskip

For the sake of simplicity we assume that $V$
---and therefore $Y $ and $W$---  satisfies Assumption~\ref{assumption}.

\begin{lem} \label{sepchow} Let $V\subset \A^n$ be an equidimensional variety of degree
bounded by $D$ which satisfies Assumption \ref{assumption}. Let $g\in \Q[x_1,\dots,x_n]\setminus \{0\}$ of degree bounded by $d$
and
$Y$ and $W$ defined as above.
Assume that $\Ch_V$ and
$g$ are encoded by  slp's of length bounded by $L$.

Then  there  is a bounded probability algorithm (Procedure $\Sep$
(Subroutine~\ref{algo-sepchow}) below) which computes the Chow
forms $\Ch_Y$ and $\Ch_W$ within (worst-case) complexity \
$\cO((ndD)^8 L)$.
\end{lem}

\begin{proof}{Proof.--}
Let $\cP_V \in \Q[U_0, \dots, U_r][T_0, \dots, T_r]$ be the
normalized characteristic polynomial of $V$, as defined in
Subsection \ref{Geometric resolutions} and set $P': = \partial
\cP_V/\partial T_0$.

We consider the following map, already introduced in  Subsection
\ref{Geometric resolutions}:
$$
\varphi_V  :   \A^{(r+1)(n+1) }\times V \to V(\cP_V) \quad  ,
\quad (u_0, \dots, u_r ; \, \xi) \mapsto ( u_0, \dots, u_r; \
L_0(u_0, \xi),  \dots, \ L_r(u_r, \xi)).
$$
By Lemma~\ref{resolucion generica} $\varphi_V$  is a birrational map
which in fact is an isomorphism when restricted to   $$ U:=
(\A^{(r+1)(n+1) }\times V) \setminus V(P'(L_0, \dots, L_r)) \to
\cU:= V(\cP_V) \setminus V(P'), $$  with inverse
$$
\psi_V : (u_0, \dots, u_r; \, t_0, \dots t_r) \mapsto \left( u_0,
\dots, u_r; \, - \frac{1}{P'}\, \frac {\partial \cP_V}{\partial
U_{01}}, \dots, - \frac{1}{P'}\, \frac{\partial \cP_V}{\partial
U_{0n}} \right).
$$
Define  $$G:= (P')^d \, \psi_V^*(g)=
g^h \left(P',- \frac{\partial \cP_V}{\partial U_{01}} ,
\dots, - \frac{\partial \cP_V}{\partial U_{0n}}\right),$$
where $ g^h := {\rm Homog}(g, d)$. Thus
$\varphi_V$ induces an isomorphism between $V(g) \cap U$ and $
V(G) \cap \cU$. Hence $V(\cP_Y)$ equals the union of the
components in $V(\cP_V) $ which are contained in the hypersurface
$V(G) \subset \A^{(r+1)\, (n+1)+(r+1)}$, and $V(\cP_W)$ is the
union of the other ones. As $\cP_V$ is a squarefree polynomial we
conclude that
$$
\cP_Y:= \gcd (G, \cP_V) \quad \quad , \quad \quad \cP_W=
\frac{\cP_V}{  \gcd (G, \cP_V)},
$$
and therefore, from Identity (\ref{P_V}) of Section  \ref{Geometric
resolutions},
we obtain that
 $$ \cF_Y = \cP_Y(U)(0) \quad \mbox{and} \quad \cF_W =
\cP_W(U)(0)  $$ are Chow forms of $Y$ and $W$ respectively.

\medskip

Note that as $\cP_Y \,|\, \cP_V$,  $\cP_Y(e_0,\dots,e_r,0,\dots,0) \ne 0$,
thus $e:=(e_0,\dots,e_r,0,\dots,0)$ is a good base point
 to apply Procedure $\PolynomialDivision$. Thus the only
 probability step of this algorithm is the computation of the
 Greatest Common Divisor between $\cP_V$ and $G$.

\medskip

Now we estimate the (worst-case) complexity of the algorithm:

The characteristic polynomial  $\cP_V$ can be computed from
$\Ch_V$ with complexity $\cO(L)$ using  Identity~(\ref{P_V}) in
Section \ref{Geometric resolutions}. Its partial derivatives
with respect to $T_0$ and $U_{01},\dots,U_{0n}$ can  be  computed
within complexity  $\cO(n\, L)$.

The polynomial $G$ is obtained within
complexity $\cO(d^2 (d+nL))$.

As $ \deg \cP_V = (r+1) \, D$ and $\deg G \le d\,((r+1)D -1)$,
both bounded by $(r+1)\,d\,D$, the gcd computation of $\cP_Y$
requires $(n  d D)^{6} (d+L) $ additional arithmetic operations.

>From  Lemma~\ref{vermeidung}, the polynomial division for
$\cP_W$  is then  performed within complexity $\cO((n d D)^{8}
L)$.

The final specialization $T\mapsto 0$ does   not change this
estimate.

Therefore, the (worst-case) complexity of the algorithm is of
order $\cO((n  d D)^{8} L)$.
\end{proof}

We summarize the algorithm in Procedure $\Sep$ (Subroutine
\ref{algo-sepchow}.

\begin{algo}
{Separation of varieties} {algo-sepchow} {$ \Sep(n,r,  D ,g, d,
\, \Ch_V)$ }

\hspace*{4.1mm}{\# $n$ is the number of variables,

\hspace*{4.1mm}\# $r, D$ are the dimension and an upper bound for
the degree of $V$ respectively,

\hspace*{4.1mm}\# $ g \in \Q[x_1, \dots, x_n] \setminus \{ 0\}$,

\hspace*{4.1mm}\# $d$ is a bound for the degree of $g$,

\hspace*{4.1mm}\# $\Ch_V$ is the normalized  Chow form of $V$.

\smallskip

\hspace*{4.1mm}\# The procedure returns the  normalized   Chow
forms $\Ch_Y$ and $\Ch_W$.

\begin{enumerate}

\item $\cP_V:=  \Ch_V((U_{00}-T_0,U_{01}\dots, U_{0r}), \dots,
(U_{r0}-T_r, U_{r1}, \dots, U_{rn}))$;\vspace{-1,5mm}

\item $ g^h := {\rm Homog}(g, d)$;\vspace{-1,5mm}

\item $\displaystyle{ G:=g^h \left(P',- \frac{\partial \cP_V}{\partial U_{01}} ,
\dots, - \frac{\partial \cP_V}{\partial U_{0n}}\right)}$;\vspace{-1,5mm}

\item  $\cP_Y:= \GCD ( G, \cP_V, (U_0, \dots, U_r, T_0, \dots,
T_r),  (r+1)\,d\, D
)$; \vspace{-1,5mm}

\item  $\cP_W:= \PolynomialDivision (
\cP_Y, \cP_V,  (r+1) \, D,  (e_0, \dots, e_r; 0, \dots,
0))$; \vspace{-1,5mm}

\item $ \Ch_Y:= \cP_Y(U)(0) / \cP_Y(e_0, \dots, e_r; \, 0 \dots, 0)$;
 \vspace{-1,5mm}

\item $ \Ch_W:= \cP_W(U)(0) / \cP_W(e_0, \dots, e_r ; 0 , \dots, 0)$;
 \vspace{-1,5mm}

\end{enumerate}

\hspace*{7.7mm} {\bf return}($\Ch_Y$,  $\Ch_W$); }
\end{algo}

\begin{rem}
In case the variety $V$ does not satisfy Assumption \ref{assumption},
 this procedure can be modified within the same bounds of complexity so that,
 from a Chow form of $V$, we obtain Chow forms of $W$ and $Y$.
The only problem that may appear in the previous lemma is that $\cP _Y(e)$ may be
 zero and we will not be able to accomplish the polynomial division.
To solve this, we can modify  Subroutine \ref{algo-sepchow} in the following way:
  we choose a random point so that we can apply the polynomial
  division subroutine  with error probability bounded
  by $\frac{1}{8}$ and we change the error probability of the GCD computation also
   by
  $\frac{1}{8}$ (by repeating it several times) in order that the  error probability of the
whole procedure is still bounded by $\frac 1 4$.

\end{rem}


\typeout{Equations in general position}

\subsection{Equations in general position}

\label{Equations in general position}

The algorithm we construct in Subsection \ref{Proof of Theorem 1}
works under some genericity hypotheses on the input polynomial
system. This is one of the main reasons
---but not the only one--- for the introduction of
non-determinism in our algorithm: there are no known efficient
deterministic procedures to obtain these hypotheses from a given
system. In order to achieve them we replace the system  and the
variables by {\em random} linear combinations.
Effective versions of
Bertini's and Noether normalization theorems enable us to
  estimate
the probability of success of this preprocessing.

\smallskip

The complexity of our algorithm is controlled by the geometric
degree of the input system, that is the maximum degree of the
varieties successively cut out by the equations obtained by this
preprocessing.

To define  this parameter, which is a suitable generalization
of the geometric degree of a 0-dimensional system
introduced  in \cite{GiHeMoMoPa98},
we first give the following definition:

\begin{defn}
Let  $g\in \Q[x_0,\dots,x_n]$ be an homogeneous polynomial,
$ I_g \subset \Q[x_0 \dots, x_n]_g$ be an homogeneous ideal,  and
$V \subset \P^n$ be a projective variety.
We say that $I_g$ is {\em radical of dimension $r$  outside $V_g:=V - V(g)$}  if
every primary component $\cQ$ of $I_g$ such
that $V(\cQ)_g \not\subset V_g$ is prime of dimension $ r$.

An analogous definition holds for an ideal in $\Q[x_1,\dots,x_n]_g$ and an affine variety
$V_g\subset \A^n_g$.
\end{defn}


Let $f_1, \dots, f_s, g  \in \Q[x_0, \dots, x_n] $ be homogeneous
polynomials of degree bounded by $d$, and set
$V_g =  V(f_1, \dots
, f_s)_g \subset  \P^n_g$.
We assume that $V_g \ne \P^n_g$, that is
$f_j \ne 0$ for some $j$. We also assume w.l.o.g $\deg f_j = d$
for every $j$: if this were not the case, we replace the input
system by
$$
x_i ^{d- \deg f_j} \, f_j \quad \quad  , \quad \quad  0\le i \le  n, \ 1\le j \le s.
$$

For $a_1, \dots, a_{n+1} \in \Q^s$  we set
$$
Q_i(a_i):= a_{i \, 1} \, f_1 + \cdots + a_{i \, s} \, f_s
$$
for the associated linear combination of $f_1, \dots, f_s$, which
---by the assumption that $\deg f_j = d$---
is also a system of homogeneous polynomials of degree $d$.

Let $\Delta$  be the set of $(n+1)  \times s$-matrices $A =(a_1,
\dots,  a_{n+1})^t \in \Q^{(n+1) \times s}$ such that the ideals
$ I_i(A) := (Q_1(a_1) , \dots, Q_i(a_i)) \ \subset \Q[x_0, \dots,
x_n]$, $1\le i \le n+1$, satisfy:
\begin{itemize}
\item $V(I_{n+1}(A) )_g = V_g$ in  $\P^n_g$.
\item For $1\le i \le n$,  if $V(I_i(A))_g \ne V_g$, then $ I_i(A)_g$
is a radical ideal of dimension $n-i$ outside $V_g$.
\end{itemize}

These are the first genericity hypotheses the polynomials should
verify
in order that our algorithm works.

\smallskip

For every $A \in \Delta$ we set \ $\delta(A) := \max \{ \deg
V(I_i (A) )_g \, ; \, 1 \le i \le n+1 \} $.

\begin{defn}
Keeping these   notations, the {\em geometric degree\/}
of the system
$$
f_1= 0, \ \dots, \ f_s =0 , \quad g \ne 0
$$
is defined as
$$
\delta:=\delta (f_1, \dots, f_s; \, g) := \max \{ \delta(A) \, ; \, A \in \Delta\}.
$$
\end{defn}

Note that  B{\'e}zout inequality implies \ $\delta \le d^n$.

\begin{rem}

For a system of polynomials $F_1, \dots, F_s, G  \in \Q[x_1,
\dots, x_n] $ (non-necessarily homogeneous) of degree bounded by
$d$,  the  affine analogue $\delta_{\rm aff}$ of  the
geometric degree is defined in exactly the same  manner,  but without
preparing the polynomials to make their degrees coincide.

In fact, if for $1\le i\le s$, $d_i := \deg F_i$, $d:= \max_i d_i$
 and $F_i^\h\, , \, G^\h \in \Q[x_0,
\dots, x_n]$ are the homogenizations of $F_i$ and $G$
respectively, then

$$
\delta_{\rm aff} (F_1, \dots, F_s; \, G) = \delta (x_0^{d -d_1} \, F_1^\h,
\dots,
x_0^{d - d_s} \, F_s^\h\, ; \ x_0 \, G^\h).
$$

\end{rem}

\bigskip

Let
$$
V_g = V_0 \cup \cdots \cup V_{n-1}
$$
be the equidimensional decomposition of $V_g$ in $\P^n_g$, where $V_i$ is either
empty or of dimension $i$, and let $A=(a_1,
\dots, a_{n+1})^t \in \Delta$.

For $i =1, \dots, n+1$, as $I_i(A) \subseteq (f_1,\dots,f_s)$,
$V_g\subseteq V(I_i(A))_g$ always holds.
Moreover, if $V(I_i(A))_g=V_g$ for some
$i$, then $V(I_j(A))_g=V_g$ for all $j\ge i$. Also, observe that
the ideal $I_i(A)$ is generated by $i$ polynomials, so every irreducible component of $V(I_i(A))$ has dimension at least $n-i$.
Thus, we infer that for $r:=n-i, \ 0\le r\le n-1$, we have
$$
V(I_{n-r}(A) )_g  = V_r' \cup V_r \cup
\cdots\cup V_{n-1}
$$
where $V_r'$ is an equidimensional variety of
dimension $r$. (We set $V_r'= \emptyset$ for every $r$ such that
$V(I_{n-r}(A))_g=V_g$ since in these cases $V_g=V_r \cup \dots \cup V_n$.)

 From now on, $Q_i(a_i)$ will be
denoted simply by $Q_i$.

The condition that $A \in \Delta$ implies that, in case $V_r' \ne
\emptyset$, $Q_{n-r+1}$ is not a zero divisor modulo $I(V_r')_g$.
In this case,  we have
\begin{eqnarray*} V'_{r-1}\cup
V_{r-1}\cup V_r\cup \dots  \cup V_{n-1} &=& V(Q_1,\dots,Q_{n-r+1})_g \\
&=& (V'_r\cup V_r \cup \dots \cup V_{n-1})\cap V(Q_{n-r+1})\\ &=&
(V'_r \cap V(Q_{n-r+1})) \cup V_r \cup \dots \cup V_{n-1},
\end{eqnarray*}

 as for all $i$, $V_i\subset V(Q_{n-r+1})$. Hence, since
 $\dim (V'_r \cap V(Q_{n-r+1}))=r-1$, we deduce that
\begin{equation} \label{muchosV}
V_r' \cap  V(Q_{n-r+1}) = V_{r-1}' \cup V_{r-1} \cup
\widetilde{V}_{r-1}
\end{equation}
where  $$ \widetilde{V}_{r-1}=\bigcup \{\, C ; \, C \mbox{
component of } V'_r \cap V(Q_{n-r+1}) \, \cap \,( V_r  \cup
\cdots \cup V_{n-1})  \mbox{ of dimension } r-1 \, \} $$ is an
equidimensional  subvariety of
 $ V_r  \cup \cdots \cup V_{n-1}$ of dimension $r-1$.
We set $\widetilde{V}_{n-1}:= \emptyset$ and $V_{-1}': =
\emptyset$.

\bigskip

Now for $b_{0}, \dots, b_n \in \Q^{n+1}$ we consider the linear change
of variables
$$
y_k(b_k) := b_{k   0} \, x_0  + \cdots + b_{k  n} \, x_n
\qquad, \qquad   0\le k\le  n.
$$
We say that $(b_0, \dots, b_n)$ is {\em admissible} if, under this
linear change of variables, for $0\le r\le n-1$,
\begin{itemize}
\item the varieties
$
V_{r}' \cup V_{r} \cup \widetilde{V}_{r}
$
satisfy Assumption~\ref{assumption},
\item the polynomials $Q_1,\dots, Q_{n-r}\in I(V'_r)_g$ are a system of local equations
of $V_r'$ at $Z_r:=V_r' \cap V(y_1,\dots, y_r)$.
\end{itemize}

\bigskip

We  construct the polynomials $Q_1, \dots, Q_{n+1}$ and the
variables $y_0, \dots, y_n$ by choosing the coefficient vectors $a_i$,
$1 \le i \le n+1 $, and $b_k$, $ 0 \le k \le n $, at random in a
given  set. In what follows we estimate the error probability of
this procedure:

\begin{lem}
\label{prob} Let notation be as in the previous paragraphs. Let
$N$ be a positive integer and let
$$
a_i \in [0,   8  N  (d+1)^{2 n} )^s \quad  1\le i\le  n+1,
\quad ,  \quad b_k \in [0, 2   N   n^2  d^{2 n})^{n+1} \quad
 0\le k \le  n
$$
be chosen at random. Then the error probability of $A:=(a_1, \dots,
a_{n+1})^t $ being in $\Delta$ and $(b_0, \dots, b_n) $ being
admissible is bounded by   $1/N$.

\end{lem}

\begin{proof}{Proof.--}

The set of matrices $\Delta$ contains in fact a non-empty open
set of $ \Q^{{(n+1)} \times s} $: by the effective Bertini's
theorem in \cite[Lemmas 1 and 2]{Lecerf00} or a local version of
\cite[Prop. 4.3 and Cor. 4.4]{KrPaSo99} there is a non zero
polynomial $F $ with $\deg F \le 4\, (d+1)^{2 n}$  such that
$F(a_1, \dots, a_{n+1}) \ne 0$ implies that $A=(a_1, \dots,
a_{n+1}) ^t \in \Delta$.

\smallskip

Assume now that $A\in \Delta$.
By the effective Noether theorem version of
\cite[Prop. 4.5]{KrPaSo99} there is a non zero polynomial $G\in k[U_0,\dots,U_n]$ with
$$
 \deg G \le 2 \, \sum_{r=0}^{n-1} r\, \deg(V_{r}' \cup V_{r} \cup
 \widetilde{V}_{r})^2
$$
such that $G(b_0, \dots, b_n) \ne 0$ implies that under the linear change of variables given by $(b_0, \dots, b_n)$, the varieties $V_r\cup V_r' \cup \widetilde V_r$ satisfy Assumption \ref{assumption}. Since, from Identitiy \ref{muchosV},
$$
\deg (V_{r}' \cup V_{r} \cup \widetilde{V}_{r}) \le d\, \deg
V_{r+1}  \le d^n ,
$$
$\deg G \le n (n-1) d^{2 n} $.

Now we will define a polynomial $H\in k[U_1,\dots, U_{n-1}]$ such that $H(b_1,\dots, b_{n-1})\ne 0$ implies that the second condition for admissibility holds.

Fix $r$, $0\le r \le n-1$. We know that $(Q_1,\dots,Q_{n-r})_g$ is a
radical ideal of dimension $r$ outside $V_g$ whose associated variety coincides
with $V'_r$ outside $V_g$. Thus, localizing
at any $\xi \in V'_r$, $\xi \notin V_g$,
we get $((Q_1,\dots,Q_{n-r})_g)_\xi = I(V'_r)_\xi$, that is, $Q_1, \dots, Q_{n-r}$ is a system of local equations of $V_r'$ at $\xi$.

Therefore, it suffices to take new variables $y_0, \dots, y_n$ such that $V_r'\cap V_g \cap V(y_1,\dots, y_r) = \emptyset.$

{}From the definition of $V_0'$, it is clear that $V_0' \cap V_g = \emptyset$.

For $1\le r \le n-1$,
as $V_g$ is definable by polynomials of degrees bounded by $d$ and no irreducible component of $V_r'$ is contained in $V_g$, there exists
a polynomial $g_r\in k[x_1,\dots, x_n]$ with $\deg(g_r)\le d$ such that $V_g\subset V(g_r)$ and $V_r' \cap V(g_r)$
is equidimensional of dimension $r-1$. Let $\cF_r\in k[U_1,\dots, U_{r}]$ be a Chow form of $V_r' \cap V(g_r)$.

Set $H:= \prod_{r=1}^{n-1}\cF_r\ \in k[U_1,\dots, U_{n-1}]$. The condition $H(b_1,\dots, b_{n-1})\ne 0$ implies that, for every $1\le r \le n-1$, $\cF_r(b_1,\dots, b_r) \ne 0$ and
so,
$$V_r' \cap V_g \cap V(y_1(b_1), \dots, y_r(b_r))\subset V_r' \cap V(g_r) \cap V(y_1(b_1), \dots, y_r(b_r)) = \emptyset.$$

Observe that $H$ is a non zero polynomial with
$$\deg H = \sum_{r=1}^{n-1} \deg \cF_r = \sum_{r=1}^{n-1} \deg V_r'\cap V(g_r) \le
 \sum_{r=1}^{n-1} d \deg V_r'\le \sum_{r=1}^{n-1} d^{n-r+1} \le (n-1) d^n.$$
Therefore, there exists a non zero polynomial condition of degree bounded by $n(n-1) d^{2n} +(n-1)d^n \le n^2 d^{2n}$ which ensures that the matrix $(b_0,\dots, b_n)$  is admissible.

The conclusion follows  as usual from the Zippel-Schwartz test.
\end{proof}


\typeout{Proof of Theorem 1}

\subsection{Proof of Theorem 1}

\label{Proof of Theorem 1}

\vspace{2mm}

Let $f_1, \dots, f_s, g  \in \Q[x_0, \dots, x_n] $ be homogeneous
polynomials of degree bounded by $d$, and set  $V_g:= V(f_1,
\dots , f_s)_g \subseteq \P^n_g$. Set $\delta$ for the geometric
degree of the system $f_1=0,\dots,f_s=0,g\ne 0$.

\medskip

The algorithm is iterative and consists of two main steps, besides
the preparation of the input (equations and   variables).

The preparation of the input enables us to work with an affine
variety $W$ instead of the input quasi-projective variety $V_g$ and local systems of equations of certain auxiliary varieties appearing in the process.

The first main step computes recursively the
Chow forms of a
non-minimal equidimensional decomposition of $W$. Here the crucial
point which controls the explosion of the complexity
is that the size of the input of an iteration does not depend on the size of
the output of the previous step: the input of each recursive step has the same
controlled size.

The second main step clears extra components and computes
the Chow forms of  the equidimensional components of the minimal decomposition of
$W $ from which the Chow forms of the equidimensional components of $V_g$ are
obtained straightforwardly.

This is a bounded error probability algorithm
whose expected complexity is of order $s (n d \, \delta)^{\cO(1)}L$.
Its worst-case complexity is $s(n d^n)^{\cO(1)}L$.

\medskip
For the rest of this proof, we set $N:=d^{56n}$.

\medskip
{\em Input Preparation}

\smallskip
Set $V_g=V_0 \cup \dots \cup V_n$ for the  minimal equidimensional
decompositon of $V_g$, where each $V_r$ is either empty or of pure
dimension $r$.

\smallskip
First, applying Procedure Deg
described at the end of  Section
\ref{Complexity of basic computations} to $f_1,\dots, f_s$,
we compute with error probability bounded by
$ 1 / 6N$ the exact degree of the polynomials  $f_1,\dots, f_s$ within complexity
$\cO( sd^2L + n \log(sdN))$.
This   also states
whether these polynomials are the zero polynomial and, therefore,
whether $V_g = \P^n_g$.  In that case  $\cF_{V_n}= |( U_0,\dots,U_n) |$ and for $i< n$,
$\cF_{V_i}=1$.

\smallskip
Thus, with error probability bounded by $1/6N$ we
can assume we know the exact degree of the polynomials
$f_1,\dots,f_s$, and
 that $V_n=\emptyset$ and $\dim V_g \le n-1$.

\smallskip

We consider the polynomials
$$\tilde f_{ij}:=x_i^{d-\deg f_j} f_j \qquad , \qquad 0\le i \le n\ , \ 1\le j \le s,$$
hence we have now $t\le (n+1)s $ polynomials $\tilde f_{ij}$ of
degree $d$, that we rename $\tilde f_1,\dots,\tilde f_t$.

We apply Lemma \ref{prob} to choose randomly a matrix
$A=(a_1,\dots,a_{n+1})^t \in \Q^{(n+1)\times t}$
and a matrix $B=(b_0,\dots,b_n) \in \Q^{(n+1)\times(n+1)}$ such
that the error probability that $A\in \Delta$ and $B$ is admissible is
bounded by $1/6N$.

We can assume thus that
 the linear
combinations $(Q_1,\dots,Q_{n+1})=A(\tilde f_1,\dots,\tilde f_t)$
 satisfy

 $$V_g =
V(Q_1,\dots,Q_{n+1})_g$$
 and, for $0\le r\le  n-1$,

\begin{itemize}
\item    $(Q_1,\dots,Q_{n-r})_g$
is either empty outside $V_g$ or a radical ideal of dimension $r$ outside $V_g$.
 \item
$
V(Q_1, \dots, Q_{n-r} )_g  = V_r' \cup V_r \cup \cdots \cup V_{n-1},
$
where $V_r'$ is  either  empty   or  an equidimensional
variety of dimension $r$ with no irreducible component included
in $V_r \cup \dots \cup V_{n-1} $.
\item
$
V_r' \cap V(Q_{n-r+1})_g = V_{r-1}' \cup V_{r-1} \cup \widetilde
V_{r-1}$,
where \ $ {\widetilde V_{r-1}}$ is either empty or an
equidimensional variety of dimension $r-1$ included in $  V_r
\cup \cdots \cup V_{n-1}$. We set $V_n':= \P^n_g$  to extend this property to
$r=n$.
\end{itemize}

We can assume  moreover  that  the change of coordinates
$y=B\,x$ verifies
\begin{itemize}
\item $
B V_{r}' \cup B V_{r} \cup B \widetilde{V}_{r}
$
satisfies Assumption~\ref{assumption}
\item $
Q_1(B^{-1}y),\dots, Q_{n-r}(B^{-1}y)$ is a system of local equations
of $BV'_r$ at $BV'_r \cap V(y_1,\dots,y_r)$.
\end{itemize}

The complexity of constructing the random matrices $A$ and $B$
 and the inverse of the matrix
$B$
is of order
$\cO( s n^4(\log N + \log d))$.


Now,
Assumption \ref{assumption} implies that the  varieties
have no irreducible component at  infinity.
Hence we restrict to  the affine space: we set $y_0=1$ and
denote by $q_1,\dots,q_{n+1}, h$ the set of polynomials in the new
variables obtained from $Q_1,\dots,Q_{n+1}, g$, that is:
$$(q_1,\dots,q_{n+1})= A \, F (B^{-1}(1,y_1,\dots,y_n) ) \quad ,
\quad h= g(B^{-1}(1,y_1,\dots,y_n)),$$

where $F:=(\tilde f_{1},\dots, \tilde f_{t})$.

\smallskip
We define
 $$W: = \overline{V
(q_1, \dots, q_{n+1} )_h } = \overline{B \,V  \cap \A^n_h} \subset
\A^n.$$
 Let
  $W=W_0 \cup \dots \cup W_{n-1}$
be the minimal equidimensional decomposition of $W$, where for $0\le r\le n-1$,
 $W_r$ is
either empty or of dimension $r$,
and let  $W_r'$ and
$\widetilde{W}_r$ defined by the same construction  as $V_r'$
and $\widetilde V_r$ before, that is
\begin{itemize}
\item $V(q_1,\dots, q_{n-r})_h=W_r' \cup W_r \cup \cdots \cup W_{n-1}$
\item $W_r' \cap V(q_{n-r+1})_h = W_{r-1}'\cup W_{r-1} \cup \widetilde W_{r-1}$
\end{itemize}
As the identity
$$
W_r= \overline{B \, V_r  \cap  \A^n_h}
$$
holds,
from a  Chow form of $W_r$ we obtain a Chow form of the corresponding  $V_r$
by means of the change of variables:
$$\cF_{V_r} (U_0, \dots, U_r) = \Ch_{W_r} ( U_0 B^{-1}, \dots, U_{r} B^{-1}).$$

\medskip
We observe that $W_r'= \overline{B \, V_r'  \cap  \A^n_h}$, and then
$q_1,\dots,q_{n-r}$ is a system of
local equations of $W'_r$ at $W'_r\cap V(y_1,\dots,y_r)$.

\medskip

The error probability of this preparation step is bounded by
$1/3N$. Once the matrices $A$ and $B$ are fixed, we have that the
complexity of computing the polynomials $q_1,\dots,q_{n+1}, h$ and
their length are all of order   $\cO( s n^2 d L)$.

\medskip
{\em First Main Step}

\smallskip
{}From $r = n-1$ to $0$, the algorithm computes the Chow form of
$W_r \cup \widetilde W_r$ and a geometric resolution of the fiber
$Z_r: = W_r'\cap V( y_1, \dots, y_r)$ (which also gives the degree
$D_r$ of $W_r'$). The former will be the input of the second main
step while the latter is the input of the next step in this
recursion. Each step of this recursion is a bounded probability
algorithm whose error probability is bounded by $1/3nN$ provided
that the input of the iteration step was correct.

\smallskip
We begin
with the fiber $Z_n   = V(y_1,\dots, y_n)= (0,\dots,0)$ and its geometric
resolution $(t, (t,\dots,t))$ associated to  $\ell=x_1$. We also set $D_n:=1$.

\smallskip

Now, we are going to describe a step of the recursion. From a
geometric resolution of $Z_{r+1}$ we compute a Chow form for
$W_r \cup \widetilde W_r$ and a geometric resolution of $Z_r$, which
is the input of the next recursive step. Set $D_{r+1}$ for the given estimate of
$\deg W'_{r+1}$.

\smallskip
\begin{itemize}
\item Computation of $\Ch_{W'_{r+1}}$:

{}From the geometric resolution $(p_{r+1}, (v_1,\dots, v_n))$
associated to the  affine linear form $\ell_{r+1}$  of $Z_{r+1}$,
and the system of local equations $q_1,\dots,q_{n-r-1}$ of
$W'_{r+1}$ at $Z_{r+1}$, we compute  the Chow form of $W_{r+1}'$
applying  Procedure $\ChowForm$ (Subroutine \ref{algo-fibra}).
This step of the algorithm
  is deterministic and computes $\Ch_{W_{r+1}'}$ provided
that the  polynomials and variables satisfy the genericity
conditions and that the geometric resolution of $Z_{r+1}$ is
accurate.  Observe that by  Main Lemma \ref{fibra}   applied to the local system
of equations $q_1,\dots,q_{n-r-1}$ of degree $d$ and length $\cO(s n^2 d L)$,
the complexity and the
length of the output are both of order
$$L(\Ch_{W'_{r+1}})= \cO\Big((r+1)^8\log_2((r+1)D_{r+1})n^7 d^2 D_{r+1}^{11} (s n^2 d L)\Big)=
 \cO\Big(s n^6 (n d D_{r+1})^{12}L\Big).$$

\item Computation of $\Ch_{W_{r+1}' \cap V(q_{n-r})}$:

Now we apply sufficient  times Procedure $\Inter$
(Subroutine \ref{algo-inter}) to
compute
the Chow form of $W_{r+1}' \cap V(q_{n-r})$ with error probability bounded by
$
1 /18 n N$: by Lemma \ref{intersection}, the length of the output Chow form and
the complexity of one iteration are
both of order
$$ L( \Ch_{W_{r+1}' \cap V(q_{n-r})})= \cO\Big((n d D_{r+1})^{12} L(\Ch_{W_{r+1}'})\Big)
=\cO\Big(  s n^6(n d D_{r+1})^{24}L \Big),$$

while, from Corollary \ref{remarkrepeticion} for the choice
$s=\lceil6(\log (18 n N )+1)\rceil$,
 the   complexity   of this step is of order
$$\begin{matrix}
\cO\Big( ((r+1)(n+1)+1) \,  \log (18 n N) (L( \Ch_{W_{r+1}'}) +
L( \Ch_{W_{r+1}' \cap V(q_{n-r})})) + \log^2 (18 n N)\Big) \ = \\
= \  \cO\Big( \log^2\!( N)\,  s\, n^9(n d D_{r+1})^{12}
L\Big).\end{matrix} $$

\item Computation of $\Ch_{W_{r}\cup W_r \cup \widetilde W_r}$:

Observe that
each
irreducible component of $W_{r+1}' \cap V(q_{n-r})$ is either an
irreducible component of $W_r'\cup W_r \cup \widetilde W_r$ or an
irreducible variety included in $V(h)$.
Therefore, we apply sufficient times Procedure $\Sep$ (Subroutine
\ref{algo-sepchow}) to compute the Chow form of $W_r'\cup W_r \cup \widetilde
W_r$ with error probability bounded by
$
1 /18 n N$: by Lemma \ref{sepchow}, the length of the output Chow form and
the complexity of one iteration are
both of order
$$ L( \Ch_{W_{r }' \cup W_r \cup \widetilde
W_r })=\cO\Big((n d(d  D_{r+1}))^{8} L(\Ch_{W_{r+1}'\cap
V(q_{n-r})}) \Big)=
 \cO\Big( s n^6 d^{8}(n d D_{r+1})^{32} L\Big),$$

while the  complexity   of this step is of order
$$\cO\Big(   \log^2\!(N)  s n^9d^{8}(n d D_{r+1})^{32} L\Big). $$

\item Computation of $\Ch_{W'_{r}}$ and $\Ch_{W_r \cup \widetilde W_r}$:

Next, since $W_r \cup \widetilde W_r \subset V(q_{n-r+1})$ and no component of
$W'_r$ does, we use  $q_{n-r+1}$ to  separate $\Ch_{W_r'}$ from
$\Ch_{W_r \cup \widetilde W_r}$.
We apply sufficient times Procedure $\Sep$ (Subroutine
\ref{algo-sepchow}) to compute the Chow forms of $W_r'$ and $W_r \cup \widetilde
W_r$ with error probability bounded by
$
1 /18 n N$:  the length of the output Chow forms and
the complexity of one iteration are
both of order
$$ L( \Ch_{W_{r }'},\Ch_{ W_r \cup \widetilde
W_r })=\cO\Big((n d (d D_{r+1}))^{8} L(\Ch_{W_{r}'\cup W_r \cup
\widetilde W_r})\Big) =
   \cO\Big( s n^6 d^{16}(n d D_{r+1})^{40} L\Big),$$
while the  complexity   of this step is of order
$$\cO\Big(   \log^2\!(N)\, s n^9 d^{16}(n d D_{r+1})^{40} L \Big) .$$

\item Computation of a geometric resolution of $Z_r:=W'_r \cap
V(y_1,\dots,y_r)$:

We apply here Procedure $\GeomRes$ (Subroutine \ref{algo-geomres}). It requires a
random choice of the coefficients  of  a separating  linear
form $\ell_r$. We do that in order to insure that the error probability is $1/6nN$.
 The condition that a linear form separates the points
of the fiber $Z_r$ is given by a polynomial of degree bounded by
$\binom{\deg Z_r }{2}\le \frac{d^{2(n-r)}}{2}$  as $\deg Z_r\le d^{n-r}$.
So we choose the set of coefficients of $\ell_r$ in
$[0,3nNd^{2(n-r)})^{n+1}$. The complexity of constructing these
coefficients is thus of order $\cO( (n+1)( \log (n N) + (n-r)\log
d))=\cO(n^2(\log N + \log d))$ and the complexity of computing afterwards
 the geometric resolution of
$Z_r$ (that is, all its constant coefficients) adds, as $D_r\le d
D_{r+1}$,
$$\cO\Big(n(d D_{r+1})^2 L(\Ch_{W'_r}) +  d^4D_{r+1}^4\Big)=
\cO\Big(s n^5 d^{16}(n d D_{r+1})^{42} L \Big) $$ operations.
\end{itemize}

Summarizing, from the geometric resolution of $Z_{r+1}$ and the
polynomials $q_1,\dots, q_{n-r}$,   the algorithm produces,
within complexity $\cO\Big(   \log^2\!(N)  s n^7 d^{16}(n d
D_{r+1})^{42} L \Big)$,
all the coefficients of the geometric resolution of $Z_r$
and a slp of length $ \cO\Big( s n^6 d^{16}(n d D_{r+1})^{40} L\Big)$ for the Chow
form of $W_r\cup \widetilde W_r$.
The
error probability that the computed objects are not the correct ones, provided that the
input was right, is bounded by $1/3nN$.

Therefore, provided that the input preparation was correct, this algorithm is expected
to compute $\Ch_{W_r \cup \widetilde W_r}$, for $0\le r\le n-1$,  with error probability bounded by $1/3N$, within complexity of order
$$\cO\Big( \log^2\!(N) \, s n^{7}d^{16}(\sum_{k=r+1}^{n-1}(n d D_{k})^{42}) L\Big),$$
and, by the iterative character of the algorithm,  to compute all $\Ch_{W_r \cup
\widetilde W_r}$, $0\le r\le n-1$, within the same complexity than
that of computing $\Ch_{W_0 \cup \widetilde W_0}$.

\medskip

{\em Second Main Step}

\smallskip
For $0\le r\le n-1$,
in order to extract from the Chow form  $\Ch_{W_r \cup \widetilde
W_r}$ the factor $\Ch_{W_r}$, we define a hypersurface $V(G_r)$
such that, probabilistically, $\widetilde W_r$ is exactly the
union of all the irreducible components of $W_r\cup \widetilde
W_r$ contained in $V(G_r)$, and then we apply Procedure $\Sep$ (Subroutine
\ref{algo-sepchow}) to compute $\Ch_{W_r}$.

\smallskip
Fix $k$, $1\le k \le n-1$.
We  define a polynomial $H_k\in \Q[y_1,\dots,y_n]$ such that, with
error probability bounded by $1/6(n-1)N$, the following conditions hold:
\begin{enumerate}
\item $W_k\cup \widetilde W_k \subseteq
V(H_{k})$,
\item no irreducible component of $W_r$ is contained in $V(H_{k})$ for $r=0,\dots, k-1$.
\end{enumerate}

Let $\cP$ be the characteristic polynomial of $W_k\cup \widetilde
W_k$. For any affine linear form $\ell_0=L_0(c_0,x)$, we have that
$H_k:=\cP(c_0,e_1,\dots, e_k)(\ell_0,
y_1,\dots,y_k)$ vanishes on $W_k\cup \widetilde W_k$. We determine
now randomly $\ell_0$ such that Condition 2 holds with error probability
bounded by $1/6(n-1)N$. This is a standard argument that can be found for instance
in \cite[Section
2.3.5]{GiHe91}):

For any irreducible component $C$ of $W_0\cup \dots \cup W_{k-1}$
there exists $\xi_C:=(\xi_{1}^{C},\dots,\xi_n^{C}) \in C - (W_k \cup \widetilde W_k)$.
Now, if a linear form $\ell_0$ satisfies that for any
$\xi \in (W_k \cup \widetilde W_k)\cap
V(y_1-\xi_1^{C},\dots, y_k-\xi_k^{C})$ (which is a zero-dimensional
variety of degree bounded by $d\,\delta$), $\ell_0(\xi) \ne
\ell_0(\xi_C)$,
$\cP(c_0,e_1,\dots,e_k)(\ell_0(\xi_C),\xi_1^{C},\dots,\xi_k^{C}) \ne 0$
holds. Hence $C$ is not included in $V(H_k)$.

The condition to be satisfied is thus given by
$$\prod_{C,\xi}(\ell_0(\xi) -
\ell_0(\xi_C))\ne 0,$$ where $C$ runs over the irreducible components of
$W_0\cup \dots \cup W_{k-1}$ and $\xi \in (W_k \cup \widetilde W_k)\cap
V(y_1-\xi_1^{C},\dots, y_k-\xi_k^{C})$. The polynomial has degree bounded by
$d\,\delta^2\le d^{2n+1}$ since $\deg W_0\cup \dots \cup W_{k-1}
\le \delta$.
Choosing $c_0:=(0,c_{01},\dots,c_{0n}) \in [0,6(n-1)N
d^{2n+1})^n$, the probability that $H_k$ does not satisfy
Condition 2 is bounded by $1/6(n-1)N$. Therefore  the
probability that, for $1\le k\le n-1$, at least one $H_k$ does not
satisfy Condition 2 is bounded by $1/6N$.

\smallskip
Now, for $r=0,\dots, n-2$ we define $G_r := \prod_{k=r+1}^{n-1}
H_k$. Clearly, as $\widetilde W_r \subset W_{r+1}\cup \dots \cup
W_{n-1}$, $G_r$ vanishes on
$\widetilde W_r$ by Condition 1. On the other hand, as,  by
Condition 2, no irreducible component
of $W_r$ is contained in $V(H_{k})$ for $r+1\le k\le n-1$,
$G_r$ splits $W_r$ and $\widetilde W_r$.

\medskip
For $1\le k\le n-1$, $\deg H_k \le  d D_{k+1}$ and
$L(H_k)=\cO(L(\Ch_{W_k\cup\widetilde W_k}))$ since we derive
$\cP$ from the corresponding  Chow form by Identity \ref{P_V}.
Hence $L(H_k) = \cO\Big(  s n^6 d^{16} (n d D_{k+1})^{40} L\Big)$.
Thus, for $0\le r\le n-2$, $\deg G_r \le d \sum_{k\ge r+1}
D_{k+1}$ and $$L(G_r) = \cO\Big(  s n^{6} d^{16 }(\sum_{k\ge r+1}
(n d  D_{k+1})^{40}) L\Big).$$ The computation of  all $H_k$,
$1\le k\le n-1$, involves the computation of the random
coefficients of each linear form $\ell_0$, that is $\cO(n^2(\log N+ n \log d))$
operations for each one of them, plus the complexity
of computing and specializing each characteristic polynomial.
Thus the total complexity of computing all $H_k$ is of order
$\cO\Big( n^2 \log N+s n^{6}d^{16}(\sum_{k\ge 2} (n d
D_{k})^{40})L\Big)$. We conclude that the complexity of computing
all $G_r$, $0\le r\le n-2$, is also of the same order.

This algorithm is expected to compute the right polynomials
$G_0,\dots,G_{n-2}$, provided that the Input Preparation and the First Main Step were
correct, with error probability bounded by $1/6N$.

\medskip

Now we apply sufficient times Procedure $\Sep$ (Subroutine
\ref{algo-sepchow}) to $\Ch_{W_r\cup \widetilde W_r}$ and $G_r$ in order
 to compute $\Ch_{W_r}$ with error probability
bounded by $1/6nN$:  the length of the output Chow forms and the
complexity of one iteration are both of order
$$ L( \Ch_{W_{r }})=\cO \Big( (n (d \sum_{k=r+1}^{n-1}D_{k+1}) (d D_{r+1}))^{8}
L(\Ch_{W_{r} \cup \widetilde W_r}, G_r)\Big)
= \cO (s n^7 d^{16} (n d \overline D )^{56}L),$$
where $\overline D = \max\{D_k: 1\le k \le n-1\}$,
 while the  total complexity   of computing all $\Ch_{W_r}$
with error probability bounded by $1/6N$, provided that the
polynomials $G_0,\dots,G_{n-2}$ were correct, is of order
$$\cO\Big( \log^2\!(N) s n^{11} d^{16}(n d \overline D)^{56} L + s \log^2(s)\, n^2 \log(d) L \Big) .$$

Thus, the total error probability of the second main step is
bounded by $1/3N$.

\medskip

Finally, the Chow form $\cF_{V_r}$ is obtained by changing variables back.
 This computation does not  change the order of complexity involved.

\smallskip

The total error probability
of the whole algorithm is bounded by $1/N$.
Moreover, in case each of the random choices was right, $D_k\le \delta$ for every $k$,
and therefore the
Chow forms $\cF_{V_r}$ of the equidimensional components $V_r$ of
$V_g$ are encoded by slp's of length
$$L(\cF_{V_r})=\cO \Big(  s n^7 d^{16}(n d \delta)^{56} L\Big),$$
and computed within complexity
$$\cO\Big(   \log^2\!(N) s n^{11} d^{16}(n d \delta)^{56} L + \Big) .$$

Since in any case, $D_k \le d^{n-k}\le d^{n-1}$ for every $1\le k\le n-1$,
the worst-case complexity of the computation is of order
$$\cO\Big(   \log^2\!(N) s n^{67} d^{16} d^{56 n} L \Big) .$$

\smallskip

Therefore the expected complexity of the algorithm is
 $$\cO\Big((1-\frac{1}{N}) (\log^2\!(N) s n^{11} d^{16}(n d \delta)^{56} L) +
 \frac{1}{N}  (\log^2\!(N) s n^{67} d^{16}d^{56 n} L) \Big).$$

Fixing $N:=d^{56 n}$, we conclude that the expected complexity
of our bounded probability algorithm is of order
$$\cO\Big( \log^2 (d^{56 n}) s n^{11} d^{16}(n d \delta)^{56} L +
\log^2(d^{56 n}) s n^{67}d^{16} L \Big) = s(nd\delta)^{\cO(1)} L, $$
while the error probability is bounded by $1/N$.

\bigskip

We summarize in Procedure $\Equidim$
(Subroutine \ref{algo-equidim}) the algorithm underlying the Proof of Theorem
1.

\begin{algo}
{Equidimensional decomposition} {algo-equidim} {$ \Equidim(n,d,
\, f_1, \dots, f_s, \, g, x )$ }

\hspace*{4.1mm}{\# $f_1,\dots, f_s, g$ are homogeneous polynomials
in $\Q[x_0,\dots, x_n]$ and $x:=(x_0,\dots, x_n)$;

\hspace*{4.1mm}\# $d$ is an upper bound for the degrees of $ f_1 ,
\dots, f_s, \, g$.

\smallskip

\hspace*{4.1mm}\# The procedure returns the  Chow forms  \ $
\cF_{V_0}, \dots, \cF_{V_n}$ \ of the equidimensional components

\hspace*{4.1mm}\# of $V_g := V(f_1, \dots, f_s)_g \subset \P^n_g$.
\vspace{-2mm}
\begin{enumerate}

\item $N:= {d^{56 n}}$;\vspace{-2mm}

\item $(d_1,\dots, d_s) := {\rm Deg}(f_1,\dots, f_s, x, d;
6 s N)$;\vspace{-2mm}

\item {\bf if} $(d_1,\dots, d_s) := (-1, \dots, -1)$ {\bf
then}\label{ifpn1}\vspace{-2mm}

\item \hspace*{4.1mm} $(\cF_{V_0}, \dots, \cF_{V_{n-1}}, \cF_{V_n}):= (1, \dots, 1, |(U_0, \dots , U_n)|)$;
\label{ifpn2}\vspace{-2mm}

\item {\bf else}\vspace{-2mm}

\item \hspace*{4.1mm} $F:= (x_0^{d-d_1} f_1, \dots, x_n^{d-d_1}
f_1, \dots, x_0^{d-d_s}f_s, \dots, x_n^{d-d_s} f_s)$;\vspace{-2mm}

\item \hspace*{4.1mm} $A:={\rm RandomMatrix}(n+1, s(n+1),48 N (d+1)^{2n})$;\vspace{-2mm}

\item \hspace*{4.1mm} $B:= {\rm RandomMatrix}(n+1,n+1,\, 12 N  n^2 d^{2n} )$;\vspace{-2mm}

\item \hspace*{4.1mm} $(y_0,\dots, y_n) := B\, (x_0,\dots,
x_n)$;\vspace{-2mm}

\item \hspace*{4.1mm} $(q_1, \dots, q_{n+1}) := A \,
F(B^{-1}(1, y_1, \dots, y_n))$;\vspace{-2mm}

\item \hspace*{4.1mm} $ h := g (B^{-1} (1,y_1, \dots, y_n))
$;\label{polyh}\vspace{-2mm}

\item \hspace*{4.1mm} $\cF_{V_n}:= 1$;\vspace{-2mm}

\item \hspace*{4.1mm} $(c^{(n)}, D_n, p_n, v^{(n)} ) := ( e_1, 1, t, (t ,\dots, t)
)$;\vspace{-2mm}

\item \hspace*{4.1mm} {\bf for }  $i $ { \bf from } 1  { \bf to } $n$ { \bf do
}\vspace{-2mm}

\item \hspace*{8.2mm} $r:= n-i$;\vspace{-2mm}

\item \hspace*{8.2mm}
$\Ch_{W_{r+1}'} := \ChowForm(n, r+1, D_{r+1}, c^{(r+1)}, p_{r+1},
v^{(r+1)}, q_1, \dots, q_{n-r-1}, d) $;
\vspace{-1.5mm}

\item \hspace*{8.2mm}  $ \cF := \Inter(n, r+1, D_{r+1}, q_{n-r}, d,
\Ch_{W_{r+1}'};
18 n N)$;\vspace{-1.5mm}

\item \hspace*{8.2mm} $\Ch_{W_{r}' \cup W_{r} \cup \widetilde{W}_{r}}: =
( \Sep  (n, r+1, d\, D_{r+1}, h, d, \, \cF;
18 n N))_2$;\vspace{-1.5mm}

\item \hspace*{8.2mm} $(\Ch_{W_r'}, \Ch_{W_{r} \cup \widetilde{W}_{r}}):=
 \Sep  (n, r, d\, D_{r+1}, d, \, \Ch_{W_{r}' \cup W_{r} \cup \widetilde{W}_{r}}, q_{n-r+1};
18 n N)$;\vspace{-1.5mm}

\item \hspace*{8.2mm} $c^{(r)}:= {\rm Random}(n+1,
3n N d^{2(n-r)})$;\vspace{-1.5mm}

\item \hspace*{8.2mm} $(D_{r}, p_{r}, v^{(r)}):= \GeomRes(n, r, d D_{r+1}, \Ch_{W_r'},
(0,\dots, 0), c^{(r)})$; \vspace{-2mm}

\item \hspace*{4.1mm} {\bf od};\vspace{-2mm}

\item \hspace*{4.1mm} {\bf for }  $k$ { \bf from } 0 { \bf to } $n-1$ { \bf do}
\vspace{-2mm}

\item \hspace*{8.2mm} $\cP_k := \Ch_{W_k \cup \widetilde W_k} ((U_{00}-T_0, U_{01}, \dots, U_{0n}), \dots,
(U_{k0}-T_r, U_{k1},\dots, U_{kn}))$;\vspace{-1.5mm}

\item \hspace*{8.2mm} $u^{(k)} := {\rm Random}(n+1,
6(n-1) N d^{2n+1}$;\vspace{-1.5mm}

\item \hspace*{8.2mm} $H_k := \cP_k(u^{(k)}, e_1, \dots, e_k)(u^{(k)}_0 + u^{(k)}_1
y_1 + \cdots + u^{(k)}_n y_n, y_1,\dots, y_k)$; \vspace{-1.5mm}

\item \hspace*{4.1mm} {\bf od};\vspace{-2mm}

\item \hspace*{4.1mm} {\bf for }  $r$ { \bf from } 0 { \bf to } $n-2$ { \bf
do}\vspace{-2mm}

\item \hspace*{8.2mm} $G_r := \prod_{k=r+1}^{n-1}
H_k$;\vspace{-1.5mm}

\item \hspace*{8.2mm}  $\Ch_{W_r}:= \Sep ( n, r, d\, D_{r+1}, G_r, d (D_{r+2} +
\cdots + D_{n}), \, \Ch_{ W_{r} \cup \widetilde{W}_{r}};
{6nN})$;\vspace{-1.5mm}

\item \hspace*{8.2mm} $\cF_{V_r} := \Ch_{W_r}(U_0 B^{-1}, \dots, U_r
B^{-1})$;\vspace{-1.5mm}

\item \hspace*{4.1mm} {\bf od};\vspace{-2.2mm}

\item {\bf fi};\vspace{-3mm}

\end{enumerate}

\hspace*{7.7mm} {\bf return}($\cF_{V_0}, \dots,
\cF_{V_n}$);\vspace{-1.5mm} }
\end{algo}

\clearpage


%% file: AplicacionesChow.tex

\typeout{Applications}

\section{Applications}

\label{Applications}

\vspace{2mm}

We  present some algorithmical applications of our results,
concerning the computation of resultants and  the resolution of
generic over-determined systems.


\typeout{Computation of resultants}

\subsection{Computation of resultants}

\label{Computation of resultants}

\vspace{2mm}

\subsubsection{The classical $d$-resultant}

\label{The classical $d$-resultant}

As a first application of our results, we  compute a slp for the
classical resultant of $n+1$ generic homogeneous polynomials of
degree $d$ in $n+1$ variables. The algorithm follows directly from
Lemma \ref{fibra} and is therefore deterministic. For the
definition and basic properties of the classical resultant we
refer for instance to \cite[Chapter 3] {CoLiOs98}).

\begin{cor}\label{resdensa} There is a deterministic algorithm
which computes  (a slp for) the classical resultant $\Res _{n,d}$
of $n+1$ generic homogeneous polynomials of degree $d$ in $n+1$
variables  within complexity $ (n d^n )^{\cO(1)}$.
\end{cor}

 \begin{proof}{Proof.--} It is a well-known fact that the resultant $\Res_{n,d}$ is the
 Chow form of the Veronese  variety $V(n,d)$ defined as the image
 of the morphism
 $$ \varphi _{n, d} : \P^n \to \P^N \quad , \quad \xi \mapsto
 (\xi^\alpha)_{\alpha \in \N_0^{n+1}, |\alpha| = d},$$

 where  $N : = \binom{n+d}{n} - 1$.  We recall that $V_{n,d}$ is an
 irreducible variety of dimension $n$ and degree $d^n$.

We compute here the resultant by defining a system of local
equations at an adequate fiber of $V(n,d)$ in order to apply Lemma
\ref{fibra}.

\smallskip
Let $\{ y_\alpha : \alpha  \in \N_0^{n+1}, |\alpha| = d\}$ be a
set of homogeneous coordinates of $\P ^N$ and consider the
projection
$$\pi : V(n,d) \to \P^n \quad , \quad (y_\alpha)_\alpha \mapsto (y_{de_0} : \dots :
y_{de_n}) $$ where $e_i$ is as usual the $(i+1)$-vector of the
canonical basis of $\Q^{n+1}$. This projection is finite
\cite[Ch. 1,  Thm 5.3.7]{Shafarevich72}. Moreover, $Z:=
\pi^{-1}((1:1:\dots: 1))$ verifies that $Z=\varphi_{n,d}(Z_0)$
with $Z_0:= \{ (1: \omega_1 :\dots :\omega_n) ; \, \omega_i^d =1
$ for $1\le i \le n \,\}$. Thus $\# Z = d^n = \deg V(n,d)$ and
the $n$-dimensional variety  $V(n,d)$ satisfies Assumption
\ref{assumption} for the fiber $Z $.

\smallskip
Let us define now a system of  local equations of $V_{n, d}$ at
$Z$:  For every $\alpha = (\alpha_0, \dots, \alpha_n) \in
(\N_0)^{n +1}$ such that $ |\alpha| = d$ and $\alpha \ne (d-1)e_0
+ e_i$ $( 0\le i \le n)$ we consider the polynomial
$$f_\alpha :=  y_{d e_0}^{d-1-\alpha_0} y_\alpha - y_{(d-1)e_0 +e_1}^{\alpha_1}
\dots y_{(d-1)e_0 +e_n}^{\alpha_n}. $$ These are $N-n$ non-zero
homogeneous polynomials of degree $d-\alpha_0$ which  vanish  at
$V_{n,d}$ since
$$f_\alpha((\xi^\beta)_\beta)=\xi_{0}^{d(d-1-\alpha_0)} \xi^\alpha -
\xi_0^{(d-1)(\alpha_1+\dots +\alpha_n)}\xi_1^{\alpha_1} \dots
\xi_n^{\alpha_n}=0.$$ From the Jacobian criterion one also checks
that, as $\frac{\partial f_\alpha}{\partial y_\alpha}=y_{d
e_0}^{d-1-\alpha_0}$ and $\frac{\partial f_\alpha}{\partial
y_\beta}=0$ for $\beta\ne \alpha$ and $\beta\ne (d-1)e_0+e_i$, the
Jacobian matrix of the system has maximal rank $N-n$ at any $\xi
\in Z$.

Observe that the  equations $f_\alpha$  can be encoded by   slp's
of length $\cO(d)$.

\smallskip
Next step in order to apply Lemma \ref{fibra} is to compute a
geometric resolution of the fiber $Z$. For that aim we compute its
characteristic polynomial (considering it as an affine variety in
$\{y_{d e_0}\ne 0\}$) and apply Lemma \ref{especializacion} for a
separating linear form.

\smallskip
Let $ L : = \sum _{|\alpha| = d} U_\alpha y_\alpha$ be a generic linear form in
$N+ 1$ variables, and let $P =  \sum _{|\alpha| = d} U_\alpha x^\alpha$
be the generic homogeneous polynomial of degree $d$ in $n+1$ variables associated to $L$.

 The
characteristic polynomial of $Z$ is
 \begin{eqnarray*}
 \cP_Z(U,T)& = & \prod_{\xi\in Z}(T-L(U,\xi))\ =\ \prod_{(1:\omega)\in
Z_0}(T-L(U,\varphi_{n,d}(1,\omega)))\\
&= & \prod_{(1:\omega)\in Z_0}(T-P(U,(1,\omega)))\ =\
\prod_{\omega :\,(1:\omega)\in Z_0}(T-P^a(U, \omega)),
\end{eqnarray*}
 where $ P^a(U,\omega)=P(U,(1,\omega))$. Therefore, if
we set $A := \Q [x_1, \dots, x_n] / (x_1^d -1, \dots, x_n^d -1)$,
$\cP_Z$ is then computed  as the characteristic polynomial of the
linear map $A \to A$ defined by $g \mapsto P^a g$ within
complexity $d^{\cO(n)}$.

\smallskip
Finally, an easy computation shows that the linear form $ \ell =
y_{de_0} + d \,  y_{(d-1)e_0 + e_1} + \dots + d^n y_{(d-1)e_0 +
e_n}$ separates the points in $Z$. Thus $\ell$ yields a geometric
resolution of $Z$ and we  apply Lemma \ref{fibra} to compute
$\Res _{n,d}$ within the stated complexity.
\end{proof}


\typeout{Computation of sparse resultants}

\subsubsection{Sparse resultants}

\label{Sparse resultants}

\vspace{2mm}

Let $\cA = \{\alpha_0, \dots, \alpha_N\} \subset \Z^n$ be a finite set
of integer vectors.
We assume that $\Z^n$ is generated by the
differences of elements in $\cA$.

\smallskip

For $0\le i \le n$ let  $U_i$ be a group of $N+1$ variables
indexed by the elements of $\cA$,  and set
$$
F_i := \sum_{\alpha \in \cA} U_{i \alpha }  \, x^{\alpha} \ \in
\Q[U_i][x_1^{\pm 1}, \dots, x_n^{\pm 1}]
$$
for the generic Laurent polynomial
with support  in $\cA$.
Let $  W_\cA \subset (\P^{N})^{n+1} \times (\C^*)^n  $
be the incidence variety of $ F_0, \dots, F_n$ in
$ (\C^*)^n $, that is
$$
W_\cA= \{ (\nu_0, \dots, \nu_n; \, \xi) \in (\P^{N})^{n+1} \times (\C^*)^n :
  \ F_i (\nu_i, \xi)=0 , \  0\le i \le  n\},
$$
and let $ \ \pi : (\P^{N})^{n+1} \times (\C^*)^n \to
(\P^{N})^{n+1} \ $ be the canonical projection. Then $
\overline{\pi (W_\cA )} $ is an irreducible variety of
codimension 1. The {\it $\cA$-resultant} \ $\Res_\cA$ is defined
as the unique ---up to  a sign--- irreducible polynomial in
$\Z[U_0, \dots, U_n]$ which defines this hypersurface (see
\cite[Ch. 8, Prop.-Defn. 1.1]{GeKaZe94}).

This is a  multihomogeneous polynomial of degree
$\Vol(\cA)$ in each group of variables $U_i$,
where  $\Vol(\cA) $ denotes the  (normalized) volume
of the convex hull $\Conv(\cA)$, which is defined as  $n!$ times its volume
with respect to the Euclidean volume form of $\R^n$.

\smallskip

Consider the map
$$
(\C^*)^n \to \P^{N} \quad \quad , \quad \quad
\xi \mapsto (\xi^{\alpha_0}: \dots: \xi^{\alpha_N}).
$$
The Zariski closure of the image of this map is called the {\it
affine toric variety}  $X_\cA \subset \P^{N} $ associated to
$\cA$. This is an irreducible variety of dimension $n$  and degree
$\Vol(\cA)$. Its Chow form coincides ---up to a scalar factor---
with the sparse resultant \ $\Res_\cA  \in \Z[U_0, \dots, U_n]$
(see \cite[Ch. 8, Prop. 2.1]{GeKaZe94} and \cite[Ch. 7, Thm.
3.4]{CoLiOs98}).

\smallskip

For a broader background on toric varieties and sparse resultants
we refer to \cite{GeKaZe94} and \cite{CoLiOs98}.

\medskip

We  apply the algorithm underlying Theorem 1 to compute the
sparse resultant $\Res_\cA$ for  the case that  $\cA \subset
(\N_0)^n$ and the elements $ 0, e_1, \dots, e_n $ ---that is the
vertices of the standard simplex of $\R^n$--- lie in $\cA$. To do
so, we construct a set of equations which define $X_\cA$ in the
open chart $(\P^N)_{y_0}$, where $(y_0:\dots:y_N)$ is a system of
homogeneous coordinates of $\P^N$, and compute a Chow form of
this variety.

\begin{cor}
Let $\cA \subset \N_0^n $ be a finite set which contains $\{0 ,
e_1, \dots, e_n \}$. Then there is a bounded probability
algorithm which computes (a slp for) a scalar multiple of the
$\cA$-resultant $\Res_\cA$ within (expected) complexity $ (n +
\Vol(\cA))^{\cO(1)} $.
\end{cor}

\begin{proof}{Proof.--}
 W.l.o.g. we assume that in $\cA$,
$\alpha_0 = 0$ and $\alpha_i = e_{i}$ for $i=1, \dots, n$. Set $
d:= \mbox{\rm max}_{\alpha \in \cA} |\alpha| $.

For $n+1\le j\le  N$ we set
$$
f_j := y_0^{d-1} y_j - y_0^{d - |\alpha_j|} y_1^{\alpha_{j1} } \cdots y_n^{\alpha_{j n}}
\ \in \Q[y_0, \dots, y_N].
$$
Then, $ X_\cA \setminus \{y_0 =0\} = V:= V(f_{n+1}, \dots,
f_N)_{y_0} \subset (\P^N)_{y_0}$. Therefore the Chow form of
$X_\cA$ coincides with the one of $V$ and can be computed by
application of Procedure $\Equidim$ (Subroutine
\ref{algo-equidim}) to the polynomial system $f_{n+1}, \dots,
f_N; y_0$.

\smallskip

Each polynomial  $f_j$, $n+1\le j\le  N$, can be encoded by a slp
of length $\cO(d)$. Moreover, as for each $\alpha \in \cA$,
$|\alpha| = \Vol( \{0, e_1, \dots, e_n , \alpha\} ) \le \Vol
(\cA) $ since $ \{ 0, e_1, \dots, e_n , \alpha\}  \subset   \cA $,
then $ d \le \Vol (\cA). $ Therefore    $L(f_j) \le \cO(
\Vol(\cA)) $ for $n+1 \le j\le  N$.

\smallskip

Now, as the toric variety $X_\cA$ is non-degenerated (that is, it
is not contained in any hyperplane in $\P^N$),   \cite[Cor.
18.12]{Harris}  implies that
$$N+1 \le \dim X_\cA + \deg X_\cA = n+ \Vol(\cA).$$
This gives an estimation for the parameter $N$.

\smallskip

Finally, we have to estimate the geometric degree
$\delta(f_{n+1}, \dots, f_N; y_0)$. As we want to compute this
degree outside $\{ y_0 = 0\}$ it is enough to deal with linear
combinations of the dehomogeneized polynomials $\hat f_j$
obtained by specializing $y_0 = 1$ in the original $f_j$  for
$n+1\le j\le N$.

For $1\le i\le  N$, $n+1\le j\le N$ and $a_{i j}  \in \Q$  we set
$$
q_i  := a_{i  n+1} \, \hat f_{n+1} + \cdots + a_{i  N} \, \hat
f_N.
$$

For every $i$, the support $\Supp(q_i)$ ---that is the set of
exponents of its non-zero monomials--- is contained in $  (\cA
\times \{0\}) \cup \cS    \subset \Z^N$, where $\cS:= \{ e_{n+1},
\dots, e_N \} \subset \Z^{N}$ and then, by \cite[Prop.
2.12]{KrPaSo99},
$$
\deg V(q_{1}, \dots, q_i ) \le \Vol ( (\cA \times \{0\}) \cup
\cS) .$$ As we have that
$$
\Vol ( (\cA \times \{0\}) \cup  \cS)  =
 N! \; \vol_{\R^N} \Conv( (\cA \times \{0\}) \cup  \cS)
= n! \; \vol_{\R^n} \Conv( \cA )
= \Vol (\cA)
$$
(where $\vol _{\R^N}$ and $\vol_{\R^n}$ denote the standard
Euclidean volume forms) we infer that
 $$\delta:= \delta(f_{n+1}, \dots, f_N;y_0) \le
\Vol (\cA). $$

\smallskip


We conclude that  $\Res_{\cA}$ can be probabilistically computed
by means of subroutine $\Equidim$ within complexity
$ (N-n) ( N \, d\, \delta)^{\cO(1)} L(f_{n+1}, \dots, f_N) \le
(n + \Vol(\cA))^{\cO(1)}$.
\end{proof}

\begin{rem}

It would be interesting to improve this algorithm in order to
compute $\Res_\cA$ without any extraneous scalar factor. It would
suffice to compute this  factor as the coefficient of any
extremal monomial of $\cF_{X_\cA}$, as we know a priori that the
corresponding coefficient in $\Res_\cA$  equals $\pm 1$
(\cite[Ch. 8, Thm. 3.3]{GeKaZe94}, see also \cite[Cor.
3.1]{Sturmfels94}).

\end{rem}

\begin{exmpl}
We take the  following example  from   \cite[Exmpl. 4.13]{KrPaSo99}:
Set
$$
\cA(n,d):= \{ 0, e_1, \dots, e_n , e_1+\cdots + e_n, 2(e_1+\cdots + e_n), \dots, d \, (e_1,
+\cdots + e_n) \} \ \subset \Z^n.
$$
It is easy to check that $\Vol(\cA(n, d) ) = n \, d$, and so the
previous algorithm computes  a slp for (a scalar multiple of)
$\Res_{\cA(n,d)} $ within  \ $(n \,d)^{\cO(1)}$ \ arithmetic
operations.
\end{exmpl}


\typeout{Generic over-determined systems}

\subsection{Generic over-determined systems}

\label{Generic over-determined systems}

\vspace{2mm}

Our last application concerns the computation of the unique
solution of a generic over-determined system.

\medskip

Let $f_0, \dots, f_{n} \in \Q[x_0, \dots, x_n]$ be homogeneous
polynomials of degree  $d$.  The associated equation system  is
{\em generically} inconsistent, where generically means  if and
only if the vector of the coefficients of the polynomials does
not lie in the hypersurface $V(\Res_{n, d}) \subset
(\P^{N})^{n+1}$ defined by the  classical resultant $\Res_{n, d}$
of $n+1$ homogeneous $(n+1)$-variate polynomials of degree $d$,
and $N:= {d+n\choose n}-1 $.

\smallskip

Now assume that the  system {\em is} consistent. In this case
 the system is said to be {\em over-determined}, in the sense that its solution
set can be defined ---at least locally--- with less equations.

Under this condition the system has generically exactly  one
solution, which is a rational map of the coefficients of the
polynomials $f_0,\dots, f_n$ (see Corollary~\ref{resultante}
below). A natural problem is   to compute this rational
parametrization. In what follows we show that this
parametrization can be easily derived from the resultant, and
therefore  can be computed with our algorithm.

\medskip

In fact we treat the more general case of an over-determined
linear system on a variety. The following result seems to be
classical. However we could not find a proof in the existing
literature, so  we provide one here.

\begin{lem} \label{parametrizacion}
Let $V\subset \P^n $ be an equidimensional variety of dimension
$r$ definable over $\Q$. Let $\cF_V(U_0,\dots,U_r)$ be a Chow
form of $V$, and let $u := (u_0, \dots, u_r) \in V(\cF_V) \subset
(\P^n)^{r+1} $ be such that $ \partial \cF_V / \partial
U_{i_0j_0} (u) \ne 0 $ for some $0\le i_0 \le r, \ 0 \le j_0 \le
n$. For  $0\le i \le  r$, let $L_i(U_i,x): = U_{i0} x_0 + \dots +
U_{in} x_n$ denote the generic linear form associated to $U_i$.
Then $ V \cap V(L_0(u_0,x), \dots, L_r(u_r,x))$ consists of
exactly one element $ \xi(u) $, and
$$
\xi(u) =
\left(
\frac{\partial \cF_V }{ \partial U_{i_0 0}}
(u) : \cdots :  \frac{\partial \cF_V }{ \partial U_{i_0 n}}(u)
\right).
$$

\end{lem}

\begin{proof}{Proof.--}
As the formula stated by the Lemma is invariant by linear changes
of variables, we can assume w.l.o.g. that no irreducible
component of $V$ is contained in any hyperplane $\{ x_j =0\}$,
$0\le j\le  n$.

For $0\le i\le  r$ we set $ \ell_i(x):= L_i(u_i,x) = u_{i  0}
x_0 + \cdots + u_{i   n}   x_n \in \C[x_0, \dots, x_n]$ for the
linear form associated to $u_i \in \C^{n+1}$. Then $ V \cap
V(\ell_0, \dots, \ell_r) \ne \emptyset$ because of the assumption
$\cF_V(u) =0$. Let $\xi$ be a point in this variety. Suppose
$\xi_0 \ne 0$.

\smallskip

Set $ V^\aff \subset \A^n $ for the image of $V$ under the
rational map $\psi : \P^n \dashrightarrow \A^n$ defined by $(x_0:
\dots : x_n) \mapsto (x_1/x_0, \dots , x_n/x_0)$. Let $T:= \{
T_0, \dots, T_r\} $ be a group of $r+1$ additional variables, and
let $P := P_{V^\aff} \in \Q[U][T]$ be the characteristic
polynomial of  $V^\aff$, as defined in Subsection~\ref{Geometric
resolutions}. Then, for $0\le j\le  n$,
\begin{eqnarray*}\label{xijxi0}
0  & =&  \frac{\partial P (U,L) }{ \partial U_{i_0 j} }  (u, \xi) +
\frac{\xi_j}{\xi_0}
\, \frac{\partial P (U,L)}{ \partial T_{i_0} } (u, \xi) \\[2mm]
& =&  \frac{\partial \cF_V  }{ \partial U_{i_0 j} }(u) -
\frac{\xi_j}{\xi_0}
\,
\frac{\partial \cF_V
}{\partial U_{i_0 0} }(u).
\end{eqnarray*}
The first equality was shown in Lemma \ref{resolucion generica},
while the second follows directly from formula (\ref{P_V}) in
Subsection~\ref{Geometric resolutions}, and the fact that
$L_i(u_i,\xi) =0$ for $0\le i\le  r$.

{}From Identity (\ref{xijxi0}) and the assumption $\partial \cF_V
/ \partial U_{i_0 j_0} (u)\ne 0$,  we infer  that $\partial \cF_V
/
\partial U_{i_0 0} (u)\ne 0$ and $ \displaystyle
{\frac{\xi_j}{\xi_0} = \frac{ \partial \cF_V / \partial U_{i_0 j}
} {\partial \cF_V / \partial U_{i_0 0} }(u)} $. Therefore
$$ \xi = \left(
\frac{\partial \cF_V }{ \partial U_{i_0 0}} (u) : \cdots :
\frac{\partial \cF_V }{ \partial U_{i_0 n}}(u) \right).$$

This shows  in particular  that $\ell_0, \dots, \ell_r$ have
exactly one common root in $V \setminus \{ x_0 = 0\} $. Moreover,
as the formula for the coordinates of $\xi$ does not depend on the
chosen affine chart, we conclude that $\xi$ is the only common
root of $\ell_0, \dots, \ell_r$ in $V$.
\end{proof}

By the way, we point out a mistake in \cite[Ch. 3, Cor.
3.7]{GeKaZe94}. This Corollary would imply that the formula of
Lemma~\ref{parametrizacion} holds in case $\xi(u)$ is a simple
common root of $\ell_0, \dots, \ell_r$ in $V$. Denoting by
$\cO_{V, \xi}$ the local ring of $V$ at $\xi$, this is equivalent
to the fact that $\cO_{V, \xi} / (\ell_0, \dots, \ell_r) \cong
\C$.

\medskip

The following counterexample shows
that this is not true:  let
$$F(t,x,y):= x^2 \, (x+t) - t \, y^2 \in \Q[t, x, y] \quad , \quad
C:= V( F ) \ \subset \P^2 .
$$
$C$  is an elliptic curve with a node at $(1:0:0)$. The linear
forms $\ell_0:= L_0((0:1:0),(t:x:y))= x$ and $\ell_1:=
L_1((0:0:1),(t:x:y))= y$  have a single common root $(1:0:0)$ in
$C$, which is a simple root of $\ell_0, \ell_1$ in $C$.

On the other hand,  as  $C$ is  a hypersurface,  $\cF_C = F(M_0,
- M_1, M_2)$, where
 $M_j$ denotes the  maximal minor
obtained by deleting the $(j+1)$ column of the matrix
$$
\left (
\begin{array}{ccc}
U_{00} & U_{01} & U_{02}  \\[2mm]
U_{10} & U_{11} & U_{12}
\end{array}
\right ).
$$
A straightforward computation shows that $\partial \cF_C /
\partial U_{i  j} ((0:1:0), (0:0:1)) =0$ for every $i, j$.

\medskip

The proof given in \cite{GeKaZe94}
is based on the biduality theorem and on Cayley's trick,
and it holds in case $V$ is  {\em smooth}, and in case $u=(u_0,
\dots, u_r) $ does not lie in the singular locus of the hypersurface
$V(\cF_V)$.
This last condition is equivalent to ours.

\bigskip

Let $V\subset \P^n$ be an equidimensional variety of dimension
$r$, and set
 $ \Omega_V := V (\cF_V ) \subset (\P^n)^{r+1}$ for
 the set of (coefficients of) over-determined linear systems over $V$.
As $\cF_V$ is squarefree and  each of its irreducible factors depends on
every group of variables,
$$ \gcd (\cF_V ,
\partial \cF_V / \partial U_{i \, 0} , \dots, \partial \cF_V / \partial U_{i\, n})
= 1
$$
for $0\le i\le  r$. Then $\Theta_i := \Omega_V \setminus V
(\partial \cF_V / \partial U_{i \, 0} , \dots,
\partial \cF_V / \partial U_{i\, n}) $
is a dense open set of $\Omega_V$ and so
$$
\Psi_V: \Omega_V \dashrightarrow \P^n \quad \quad , \quad
\quad u \mapsto \xi(u):=
\left(\frac{\partial \cF_V }{ \partial U_{i 0}}
(u) : \cdots :  \frac{\partial \cF_V }{ \partial U_{i n}}(u)
\right) \quad \hbox{ if } u \in \Theta_i
$$
is a  rational map well-defined on
\ $ \Theta_0 \cup \cdots \cup  \Theta_r$.

\bigskip

Now let $V \subset (\P^n)_g$ be an arbitrary variety of dimension
$r$, and let $V = V_r \cup \cdots \cup  V_0$ be its
equidimensional decomposition. In what follows, for sake of
clarity,  we keep the same notations as previously for different
objects sharing analogous properties.

 Set (again) $\Omega_V\subset
(\P^n)^{r+1}$ for the set
$$\Omega_V= \{ (u_0, \dots, u_r) \in (\P^n)^{r+1} : \exists\, \xi \in V \ / \ L_0(u_0, \xi) = 0, \dots, L_r(u_r, \xi) = 0\}$$
of generic over-determined linear systems over $V$, which is a quasi-projective variety of codimension 1
in $ (\P^n)^{r+1}$.

For every  $0\le k \le r$, let $\Omega_{V_k}$ be the set of the coefficients of $r+1$
linear forms which have a common root in $V_k$. If $\cF_{V_k}$ is a Chow form
of $V_k$, we have that
$$\Omega_{V_k} \subset \bigcap_{0\le i_0<\dots<i_k\le r} V( \cF_{V_k}(U_{i_0},\dots, U_{i_k}))$$
and, therefore, $\Omega_{V_k}$ has codimension at least 2 for $0\le k \le r-1$.

Let $\Theta_i := \Omega_{V_r} \setminus V (\partial \cF_{V_r} / \partial U_{i \, 0} , \dots,
\partial \cF_{V_r} / \partial U_{i\, n}) $ for $i=0, \dots, r$.

Then
every over-determined linear system over $ V$ with
coefficients in the open set
$(\Theta_0\cup \dots \cup \Theta_r) \setminus (\Omega_{V_0} \cup \dots \cup \Omega_{V_{r-1}})$
of $\Omega_V$ has a unique solution in $V$ which, in fact, lies in $V_r$.
As before, this solution can be given by the rational map
$$
\Psi_V := \Psi_{V_r} :  \Omega_{V_r}\dashrightarrow (\P^n)_g \quad \quad , \quad
\quad u \mapsto \xi(u):=
\left(\frac{\partial \cF_{V_r} }{ \partial U_{i 0}}
(u) : \cdots :  \frac{\partial \cF_{V_r} }{ \partial U_{i n}}(u)
\right) \quad \hbox{ if } u \in \Theta_i.
$$

\smallskip

As an immediate  consequence of Theorem 1  and
Lemma~\ref{parametrizacion} we obtain:

\begin{cor} \label{Psi_W}

Let
$ f_1, \dots, f_s , \, g \, \in {\Q} [x_0, \dots, x_n]$
be homogeneous
polynomials  of degree bounded
by $d$
encoded by straight-line programs of length bounded by $L$.

Set $V :=  V(f_1, \dots, f_s)\setminus V(g) \subset \P^n $ for
the quasi-projective variety $\{ \, f_1=0,\dots,f_s=0,g\ne 0\,\}$
and let $ V = V_{0} \cup \cdots \cup V_{n} $ be its minimal
equidimensional decomposition. Let $\delta:= \delta( f_1,\dots,
f_s; \, g)$ be the geometric degree of the input polynomial
system.

Then there is a bounded probability algorithm which computes
(slp's for) the  coordinates of the rational map $\Psi_V$ defined
above within (expected) complexity \ $
 s (n\, d\, \delta )^{\cO(1)} L
$.
\end{cor}

The previous result  can be applied directly to compute the
solution of a generic over-determined system of $n+1$ homogeneous
polynomials in $n+1$ variables of degree $d$ by means of
$\Res_{n,d}$:

\begin{cor} \label{resultante}
Let
$u = (u_0, \dots, u_n) \in
(\P^{N})^{n+1}$
where
$N:= {d+n\choose n}-1$,
 and for $0\le i\le n$, set
$$
f_i := \sum_{|\alpha|=d } u_{ i  \alpha} \, x^\alpha .
$$

Assume that $\Res_{n, d} (u) =0$ and that ${\partial \Res_{n, d}}/ {
\partial U_{i_0 \, \beta} } (u) \ne 0$ for some $0\le i_0 \le n,
\beta = (\beta_0, \dots, \beta_n)  \in (\N_0)^{n+1}$ with
$|\beta|=d$.

Then $V(f_0, \dots, f_n) $ consists of exactly one element $
\xi(u) \in \P^n$, and
$$
\xi(u)= \left( \frac{\partial \Res_{n, d} }{\partial U_{i_0,
(d-1) \, e_j +e_0}}(u):  \cdots : \frac{\partial \Res_{n, d}
}{\partial U_{i_0,  (d-1) \, e_j +e_n }}(u) \right)
$$
for any  $ 0 \le j \le n $ such that $\beta_j \ne 0$.

\end{cor}

\begin{proof}{Proof.-}
{}From Lemma \ref{parametrizacion} applied to the Veronese variety
$V(n,d) \subset \P^{N} $ (see Section \ref{The classical
$d$-resultant}) we have that $V(f_0,\dots, f_n)$ has only one
point  $\xi(u)$ and  that
$$
(\xi(u)^\alpha)_{|\alpha | = d} = \left (\frac{\partial \Res_{n,
d} }{\partial U_{i_0   \alpha }}(u) \right)_{|\alpha | = d}.
$$
Let $\beta = (\beta_0, \dots, \beta_n) $ be such that $|\beta| =
d$ and $ \partial \Res_{n, d} / \partial U_{i_0 \beta  }(u) \ne 0
$, and let $ 0 \le j \le n $ be such that $\beta_j \ne 0$. The
previous identity implies that $ \xi \in \{ x_j \ne 0\} $.  Then
$$
\xi(u) = (\xi_j^{d-1} \, \xi_0: \cdots : \xi_j^{d-1} \, \xi_n) =
\left( \frac{\partial \Res_{n, d} }{\partial U_{i_0,  (d-1) \, e_j
+e_0}}(u):  \cdots : \frac{\partial \Res_{n, d} }{\partial U_{i_0,
(d-1) \, e_j +e_n }}(u) \right).
$$
\end{proof}

As an immediate  consequence of this result
and Proposition \ref{resdensa} we obtain:

\begin{cor}
Let notation be as in Corollary~\ref{resultante}.
Then the rational map $(\P^N)^{n+1} \dashrightarrow \P^n$,  $u \mapsto \xi(u) $ can
be (deterministically) computed within complexity
$(n d^{n})^{O(1)} $.
\end{cor}


%% file: ReferenciasChow.tex
\typeout{Referencias}


%% file: DireccionesChow.tex
\noindent {\sc Gabriela Jeronimo, Teresa Krick, and Juan Sabia: }
 Departamento de Matem{\'a}tica,
Universidad de Buenos Aires,
 Ciudad Universitaria,
1428 Buenos Aires, Argentina. \\
{\tt E-mail: jeronimo@dm.uba.ar (G.J.), krick@dm.uba.ar (T.K.),
jsabia@dm.uba.ar (J.S.)}

\vspace{4mm}

\noindent {\sc Mart{\'\i}n Sombra: }
Universit{\'e} de Paris 7, UFR de Math{\'e}matiques,
{\'E}quipe de G{\'e}om{\'e}trie et Dynamique, 2 place Jussieu, 75251 Paris Cedex 05,
France; and
Departamento de Matem{\'a}tica,
Universidad Nacional de La Plata,
Calle 50 y 115,
1900 La Plata, Argentina. \\
{\tt E-mail: sombra@math.jussieu.fr}